\documentclass[11pt]{article}

\usepackage{hyperref}

\usepackage{amsmath}
\usepackage{amssymb}
\usepackage{amsthm}
\usepackage{amscd}
\usepackage{amsfonts}
\usepackage[all]{xy}
\usepackage{ifthen}
\usepackage{a4}

\long\def\forget#1{}

\newdir^{ (}{{}*!/-3pt/\dir^{(}}    
\newdir^{  }{{}*!/-3pt/\dir^{}}    
\newdir_{ (}{{}*!/-3pt/\dir_{(}}    
\newdir_{  }{{}*!/-3pt/\dir_{}}

\def\quotes#1{{''#1''}} % Anführungsstriche

\def\INDENT{\hspace*{\parindent}}

\catcode`\ä=\active
\catcode`\ö=\active
\catcode`\ü=\active
\catcode`\Ä=\active
\catcode`\Ö=\active
\catcode`\Ü=\active
\catcode`\ß=\active
\catcode`\é=\active
\catcode`\è=\active
\catcode`\É=\active

\defä{\"a}
\defö{\"o}
\defü{\"u}
\defÄ{\"A}
\defÖ{\"O}
\defÜ{\"U}
\letß=\ss
\defé{\'{e}}
\defè{\`{e}}
\defÉ{\'{E}}

\DeclareMathOperator{\Quot}{Quot}
\DeclareMathOperator{\Spec}{Spec}

\DeclareMathOperator{\id}{id}
\DeclareMathOperator{\rank}{rank}
\DeclareMathOperator{\coker}{coker}
\DeclareMathOperator{\im}{im}

\newcommand{\sep}{{\rm sep}}
\newcommand{\alg}{{\rm alg}}
\DeclareMathOperator{\inv}{inv}
\newcommand{\op}{{\rm op}}
\DeclareMathOperator{\weight}{wt}
\DeclareMathOperator{\supp}{supp}

\DeclareMathOperator{\Hom}{Hom}
\DeclareMathOperator{\End}{End}
\DeclareMathOperator{\Isog}{Isog}
\DeclareMathOperator{\QHom}{QHom}
\DeclareMathOperator{\QEnd}{QEnd}
\DeclareMathOperator{\QIsog}{QIsog}
\DeclareMathOperator{\Gal}{Gal}

\newtheoremstyle{statement}%
{2\parskip}% space above
{\parskip}% space below
{\itshape}% body font
{}% indent amount
{\bfseries}% theorem head font
{}% punctuation after head
{\newline}% space after theorem head
{}% 'theorem head spec'

\newtheoremstyle{notice}%
{2\parskip}% space above
{\parskip}% space below
{}% body font
{}% indent amount
{\itshape}% theorem head font
{}% punctuation after head
{\newline}% space after theorem head
{}% 'theorem head spec'

\newtheorem{all}{all}[section]

\theoremstyle{plain}
\newtheorem{Definition}[all]{Definition}
\newtheorem{Lemma}[all]{Lemma}
\newtheorem{Proposition}[all]{Proposition}
\newtheorem{Corollary}[all]{Corollary}
\newtheorem{Theorem}[all]{Theorem}

\theoremstyle{remark}

\theoremstyle{definition}
\newtheorem{Remark}[all]{Remark}
\newtheorem{Example}[all]{Example}

\def\Proofof#1{Proof $($#1\,$)$.} % 'Proof (Proposition 1.x).'

\newenvironment{suchthat}
{\setlength{\parskip}{0ex}%
\begin{enumerate}\setlength{\parskip}{0ex}\setlength{\itemsep}{0ex}%
}{%
\end{enumerate}%
}%

\def\mal{^{\SSC\times}}

\def\II#1{{[\,#1\,]}}  % Polynomring Klammern, etwa 'A[t]'
\def\dual#1{{#1}^\vee} % dual (alternativ vielleicht '\widehat{#1}')

\def\Z{\mathbb{Z}} % Ganze Zahlen
\def\N{\mathbb{N}} % Natürliche Zahlen
\def\Q{\mathbb{Q}} % Natürliche Zahlen
\def\Ff{{\mathbb{F}}} % Körper F
\def\Fq{{\mathbb{F}_q}} % Körper Fq
\def\Fs{{\mathbb{F}_s}} % Körper Fa
\def\FqItI{{\Fq\II{t}}} % Fq[t]
\def\AA{\mathbb{A}} % Affiner Raum
\def\PP{\mathbb{P}} % Projektiver Raum

\def\O{{\cal O}} % Strukturgarbe
\def\I{{\cal J}} % Idealgarbe
\def\F{{\cal F}} % Modulgarbe
\def\G{{\cal G}} % Modulgarbe
\def\M{{\cal M}} % Hilfs-Modulgarbe
\def\HOM{{\cal H}\mbox{\it om}} % Hom-Garbe

\def\OC{{\O_C}} % Kurven und Strukturgarben...
\def\OS{{\O_S}}
\def\CS{{C_S}}
\def\OCS{{\O_\CS}}
\def\CL{{C_L}}
\def\OCL{{\O_\CL}}

\def\FF{{\underline{\F}}} % Abelsche Garbe
\def\GG{{\underline{\G}}} % Abelsche Garbe
\def\ZZ{{\underline{0}}} % Nullgarbe

\def\ulM{{\underline{M\!}\,}{}}
\def\ulN{{\underline{N\!}\,}{}}
\def\ulTM{{\underline{\,\,\wt{\!\!M}\!}\,}{}}
\def\ulTN{{\underline{\wt{N}\!}\,}{}}
\def\ulHM{{\underline{\hat M\!}\,}{}}
\def\ulHN{{\underline{\hat N\!}\,}{}}

\renewcommand{\mod}{\operatorname{mod}}

\def\P{{\mit\Pi}} % Pi's
\def\t{{\mit\tau}} % tau's
\def\s{{\sigma^\ast}} % sigma*

\def\TP{\widetilde{\P}} % Tilde...
\def\Tt{\tilde{\tau}}
\def\TF{\widetilde{\F}}
\def\TFF{\widetilde{\FF}}

\def\Tr{\tilde{r}}
\def\Td{\tilde{d}}

\def\HP{\widehat{\P}} % Hat...
\def\Ht{\hat{\tau}}
\def\HF{\widehat{\F}}

\def\Hr{\hat{r}}

\def\chr{\varepsilon} % Charakteristischer Punkt, im c.
\def\otimesidOCL#1{\!\otimes1} % '\otimes\id_\OCL(#1)'

\def\matr#1#2#3#4{\left(\genfrac{}{}{0pt}{}{#1}{#2}\,\genfrac{}{}{0pt}{}{#3}{#4}\right)}
\def\tmatr#1#2#3#4{{\textstyle\Big(\genfrac{}{}{0pt}{}{#1}{#2}\,\genfrac{}{}{0pt}{}{#3}{#4}\Big)}}
\def\smatr#1#2#3#4{{\scriptstyle\big(\genfrac{}{}{0pt}{}{#1}{#2}\,\genfrac{}{}{0pt}{}{#3}{#4}\big)}}
\def\vect#1#2{\left(\genfrac{}{}{0pt}{}{#1}{#2}\right)}
\def\tvect#1#2{{\textstyle\big(\genfrac{}{}{0pt}{}{#1}{#2}\big)}}

\def\TA{\tilde{A}}
\def\Lsep{{L^{\sep}}}
\def\liminv#1{\displaystyle\lim_{\stackrel{\textstyle\longleftarrow}{#1}}}
\def\CFs{{C_\Fs}} % CFs
\def\Ca{{C'}} % C affin
\def\CLa{{C'_L}} % CL affin
\def\Av{{A_v}} % Vervollständigung Av
\def\AvG{{\Av\II{G}}} % Av[G]
\def\Qv{{Q_v}} % Vervollständigung Qv
\def\QvG{{\Qv\II{G}}} % Qv[G]
\def\cl#1{{\overline{#1}}} % Algebraischer Abschluss
\def\FS{{\cl{\Ff}_s}} % FS
\def\AxL{{\tilde A\otimes_\Fq L}} % A x L
\def\VvFF{{V_v\FF}} % Tate module VvFF
\def\GalFSFs{{\Gal(\FS/\Fs)}} % Absolute Galois-Gruppe von Fs

\def\xQv{\otimes_Q\Qv}
\def\Fss{{\Ff_{s'}}}
\def\xFss{\otimes_\Fs\Fss}

\def\smallexact#1#2#3#4#5#6#7#8#9{%
{\,#1\rightarrow#3\rightarrow#5%
\ifthenelse{\equal{#7}{}}{}{\rightarrow#7% 
\ifthenelse{\equal{#9}{}}{}{\rightarrow#9}}%
\,}}

\def\exact#1#2#3#4#5#6#7#8#9{%
{#1\longrightarrow#3\longrightarrow#5%
\ifthenelse{\equal{#7}{}}{}{\longrightarrow#7% 
\ifthenelse{\equal{#9}{}}{}{\longrightarrow#9}}%
\,}}

\def\bigexact#1#2#3#4#5#6#7#8#9{%
\begin{CD}% 
#1 @>{#2}>> #3 @>{#4}>> #5 
\ifthenelse{\equal{#7}{}}{}{@>{#6}>> #7 
\ifthenelse{\equal{#9}{}}{}{@>{#8}>> #9 }}
\end{CD}%
}

\newcommand{\DS}{\displaystyle}
\newcommand{\TS}{\textstyle}
\newcommand{\SC}{\scriptstyle}
\newcommand{\SSC}{\scriptscriptstyle}
\let\setminus\smallsetminus
\DeclareMathOperator{\GL}{GL}
\DeclareMathOperator{\Graph}{Graph}
\DeclareMathOperator{\Id}{Id}
\DeclareMathOperator{\charakt}{char}
\def\isoto{\stackrel{}{\mbox{\hspace{1mm}\raisebox{+1.4mm}{$\SC\sim$}\hspace{-3.5mm}$\longrightarrow$}}}
\newcommand{\longto}{\longrightarrow}
\newcommand{\onto}{\mbox{\mathsurround=0pt \;$\longrightarrow \hspace{-0.7em} \to$\;}}
\newcommand{\shortonto}{\mbox{\mathsurround=0pt \;$\to \hspace{-0.8em} \to$\;}}
\newcommand{\into}{ \mbox{\mathsurround=0pt \;\raisebox{0.63ex}{\small $\subset$} \hspace{-1.07em} $\longrightarrow$\;}}
\newcommand{\es}{\enspace}
\newcommand{\et}{{\rm \acute{e}t}}
\newcommand{\nil}{{\rm nil}}
\DeclareMathOperator{\rk}{rk}
\newcommand{\wt}{\widetilde}
\newcommand{\wh}{\widehat}
\newcommand{\Fa}{{\mathfrak{a}}}
\newcommand{\ulK}{{\ul K}}
\newcommand{\BF}{{\mathbb{F}}}
\newcommand{\CM}{{\cal{M}}}
\newcommand{\ul}[1]{{\underline{#1}}}
\newcommand{\dbl}{{\mathchoice{\mbox{\rm [\hspace{-0.15em}[}}
                              {\mbox{\rm [\hspace{-0.15em}[}}
                              {\mbox{\scriptsize\rm [\hspace{-0.15em}[}}
                              {\mbox{\tiny\rm [\hspace{-0.15em}[}}}}
\newcommand{\dbr}{{\mathchoice{\mbox{\rm ]\hspace{-0.15em}]}}
                              {\mbox{\rm ]\hspace{-0.15em}]}}
                              {\mbox{\scriptsize\rm ]\hspace{-0.15em}]}}
                              {\mbox{\tiny\rm ]\hspace{-0.15em}]}}}}
\newcommand{\dpl}{{\mathchoice{\mbox{\rm (\hspace{-0.15em}(}}
                              {\mbox{\rm (\hspace{-0.15em}(}}
                              {\mbox{\scriptsize\rm (\hspace{-0.15em}(}}
                              {\mbox{\tiny\rm (\hspace{-0.15em}(}}}}
\newcommand{\dpr}{{\mathchoice{\mbox{\rm )\hspace{-0.15em})}}
                              {\mbox{\rm )\hspace{-0.15em})}}
                              {\mbox{\scriptsize\rm )\hspace{-0.15em})}}
                              {\mbox{\tiny\rm )\hspace{-0.15em})}}}}

\def\?{\ 
???\ \immediate\write16{}
\immediate\write16{Warning: There was still a question mark . . . }
\immediate\write16{}}

\begin{document}
% =============================================================================

\author{Matthias Bornhofen, Urs Hartl%
\footnote{Corresponding author: Urs Hartl, Institute of Mathematics, University of Muenster, Einsteinstr.\ 62, D--48149 Muenster, Germany, http:/\!/www.math.uni-muenster.de/u/urs.hartl/ }}

\title{Pure Anderson Motives and Abelian $\tau$-Sheaves over Finite Fields}

% In the final version I might want to fix the date:
%\date{August 30, 2004} 

\maketitle

\begin{abstract}
\noindent
Pure $t$-motives were introduced by G. Anderson as higher dimensional generalizations of Drinfeld modules, and as the appropriate analogs of abelian varieties in the arithmetic of function fields. In this article we develop their theory regarding morphisms, isogenies, Tate modules, and local shtukas. The later are the analog of $p$-divisible groups. We investigate which pure $t$-motives are semisimple, that is, isogenous to direct sums of simple ones. We give examples for pure $t$-motives which are not semisimple. Over finite fields the semisimplicity is equivalent to the semisimplicity of the endomorphism algebra, but also this fails over infinite fields. Still over finite fields we study the endomorphism rings of pure $t$-motives and obtain answers which are similar to Tate's famous results for abelian varieties. Finally we clarify the relation of pure $t$-motives to the abelian $\tau$-sheaves introduced by the second author for the purpose of constructing moduli spaces. We obtain an equivalence of the respective quasi-isogeny categories.

\end{abstract}

%%%%%%%%%%%%%%%%%%%%%%%%%%%%%%%%%%%%%%%%%%%%%%%%%%%%%%%%%%%%% 
%% 
%%     Introduction
%% 
%%%%%%%%%%%%%%%%%%%%%%%%%%%%%%%%%%%%%%%%%%%%%%%%%%%%%%%%%%%%% 

\thispagestyle{empty}

\section*{Introduction}

In the arithmetic of number fields, elliptic curves and abelian varieties are important objects to study and their theory has been vastly developed in the last two centuries. For the arithmetic of function fields Drinfeld \cite{Drinfeld,Drinfeld3} has invented the concepts of \emph{elliptic modules} (today called \emph{Drinfeld modules}) and \emph{elliptic sheaves} in the 1970's, both as the analogs of elliptic curves. Since then, the arithmetic of function fields has evolved into an equally rich parallel world to the arithmetic of number fields.
As for higher dimensional generalizations of elliptic modules or sheaves there are different notions, for instance Anderson's \emph{abelian $t$-modules} and \emph{$t$-motives} \cite{Anderson}, \emph{Drinfeld-Anderson shtukas} \cite{Drinfeld5,HH}, or \emph{abelian $\tau$-sheaves} which were introduced by the second author in \cite{Hl} in order to construct moduli spaces for pure $t$-motives. In the present article we advertise the point of view that pure $t$-motives and abelian $\tau$-sheaves are the appropriate analogs for abelian varieties. This is also supported by the results in \cite{Hl}. It is due to the fact that both structures have the feature of purity built in as opposed to general $t$-motives or Drinfeld-Anderson shtukas. For example non-zero morphisms exist only between pure $t$-motives or abelian $\tau$-sheaves of the same weight (see \ref{PROP.1} and \ref{Cor2.9b} in the body of the article).

There is a strong relation between pure $t$-motives, which we slightly generalize to \emph{pure Anderson motives}, and abelian $\tau$-sheaves. To give their definition let $C$ be a connected smooth projective curve over $\Ff_q$, let $\infty\in C(\Ff_q)$ be a fixed point, and let $A=\Gamma(C\setminus\{\infty\},\O_C)$. For a field $L\supset\Fq$ let $\s$ be the endomorphism of $A_L:=A\otimes_\Fq L$ sending $a\otimes b$ to $a\otimes b^q$ for $a\in A$ and $b\in L$. Let $c^\ast:A\to L$ be an $\Fq$-homomorphism and let $J=(a\otimes 1-1\otimes c^\ast(a):a\in A)\subset A_L$. A \emph{pure Anderson motive $\ulM=(M,\tau)$ of rank $r$, dimension $d$ and characteristic $c^\ast$} consists of a locally free $A_L$-module $M$ of rank $r$ and an $A_L$-homomorphism $\t:\s M:= M\otimes_{A_L,\s}A_L\to M$ with $\dim_L\coker\t=d$ and $J^d\cdot\coker\t=0$, such that $M$ possesses an extension to a locally free sheaf $\CM$ on $C\times_\Fq L$ on which $\t^l:(\s)^l\CM\to\CM(k\cdot\infty)$ is an isomorphism near $\infty$ for some positive integers $k$ and $l$. The last condition is the purity condition. The ratio $\frac{k}{l}$ equals $\frac{d}{r}$ (see \ref{Prop1'.1b}) and is called the \emph{weight of $\ulM$}. Anderson's definition of pure $t$-motives~\cite{Anderson} is recovered by setting $C=\PP^1_\Fq$ and $A=\Fq[t]$. 

In addition to this data an \emph{abelian $\tau$-sheaf} consists of a sequence of sheaves $\CM_i\subsetneq\CM_{i+1}$ lying between $\CM_0:=\CM$ and $\CM_l:=\CM(k\cdot\infty)$ whose stalks at $\infty$ are the images of $\t^i$ for $i=0,\ldots,l$ (see \ref{Def1.1}). The quasi-isogeny categories of pure Anderson motives and abelian $\tau$-sheaves are equivalent (\ref{PropX.1}, \ref{Cor2.9d}). An abelian $\tau$-sheaf of dimension $d=1$ is the same as an elliptic sheaf. In this sense abelian $\tau$-sheaves are higher dimensional elliptic sheaves. 
The concept of abelian $\tau$-sheaves was introduced by the second author~\cite{Hl} for the following reasons. In contrast to pure Anderson motives, abelian $\tau$-sheaves possess nice moduli spaces which are %
%FIX1: PROPERTIES OF THE STACKS
%irreducible and normal Deligne-Mumford stacks locally of finite type, separated, and flat over $C$;
Deligne-Mumford stacks locally of finite type and separated over $C$;
 see~\cite{Hl}. Moreover, let $c:\Spec L\to \Spec A\subset C$ be the morphism induced by $c^\ast$. The notion of abelian $\tau$-sheaves is still meaningful if $c:\Spec L\to C$ is not required to factor through $\Spec A$. Indeed, the possibility to have $\im(c)=\infty$ was crucial for the uniformization of the moduli spaces of abelian $\tau$-sheaves and the derived results on analytic uniformization of pure Anderson motives in \cite{Hl}. For these reasons we develop the theory of abelian $\tau$-sheaves and pure Anderson motives simultaneously in the present article.

Let $Q$ be the function field of $C$. Then the endomorphism algebra of a pure Anderson motive or an abelian $\tau$-sheaf is a finite dimensional $Q$-algebra (\ref{PropT.2}, \ref{ThmT.3}). In contrast the endomorphism algebra of an abelian variety is a finite dimensional algebra over the rational numbers. Through this fact pure Anderson motives and abelian $\tau$-sheaves belong to the arithmetic of function fields. We further investigate their \mbox{(quasi-)}isogenies. An isogeny $f:(M,\tau)\to(M',\tau')$ between pure Anderson motives of the same characteristic is an injective morphism $f:M\to M'$ with torsion cokernel such that $f\circ\tau=\tau'\circ\s f$. We show that in fact $\coker f$ is annihilated by an element of $A$ (as opposed to $A_L$); see \ref{Cor1.11b}. Therefore every isogeny possesses a dual (\ref{Cor1.11b}) and the group of quasi-isogenies equals the group of units in the endomorphism $Q$-algebra (\ref{QISOG-GROUP}). 
We give various other descriptions for \mbox{(quasi-)}isogenies (\ref{PROP.1.42A}, \ref{PropAltDescrQHom}). Also we prove that the existence of a separable isogeny defines an equivalence relation on pure Anderson motives over a finite field (\ref{ThmW5.2}), but not over an infinite field (\ref{Ex8.10}). 
For an isogeny $f$ between pure Anderson motives we define the degree of $f$ as an ideal of $A$ (\ref{Def1.7.6}) which annihilates $\coker f$ (\ref{Prop3.28a}). If $\ulM$ is a semisimple (see below) pure Anderson motive over a finite field, the degree of any isogeny $f:\ulM\to\ulM$ is a principal ideal and has a canonical generator (\ref{Prop3.4.1}). In particular $f$ has a canonical dual.

Next we address the question whether every pure Anderson motive is \emph{semisimple}, that is, isogenous to a direct sum of simple pure Anderson motives. A pure Anderson motive is called \emph{simple} if it has no non-trivial quotient motives. This question is the analog of the classical Theorem of Poincar\'e-Weil on the semisimplicity of abelian varieties. By giving a counterexample (\ref{Ex3.1}) we demonstrate that the answer to this question is negative in general. On the positive side we show that every pure Anderson motive over a finite base field becomes semisimple after a field extension whose degree is a power of $q$ (\ref{Thm3.8b}), and then stays semisimple after any further field extension (\ref{Cor3.8c}). The endomorphism $Q$-algebra of a semisimple pure Anderson motive is semisimple (\ref{QEND-DIVISION-MATRIX}) and over a finite field also the converse is true (\ref{Thm3.8}). This is false however over an infinite field (\ref{Ex3.10c}).

Like for abelian varieties the behavior of a pure Anderson motive $\ulM$ over a finite field is controlled by its Frobenius endomorphism $\pi$ (defined in \ref{Def2.19b}). If $\ulM$ is semisimple we determine the dimension and the local Hasse invariants of its endomorphism $Q$-algebra $E$ in terms of $\pi$ (\ref{Thm3.5a}, \ref{THEOREM-2}). We further give a criterion for two semisimple pure Anderson motives to be isogenous (\ref{THEOREM-1}). In the last section we sketch a few results for the question, which orders of $E$ occur as the endomorphism rings of pure Anderson motives (\ref{ThmW3.13}, \ref{ThmW6.1}). There is a relation between the breaking up of the isogeny class of a semisimple pure Anderson motive into isomorphism classes, and the arithmetic of $E$. We indicate this by treating the case of pure Anderson motives defined over the minimal field $\Fq$. In this case $E$ is commutative (\ref{ThmW6.1}).

Many of these results parallel Tate's celebrated article \cite{Tat} on abelian varieties over finite fields. To prove them a major tool are the Tate modules and local shtukas attached to pure Anderson motives and abelian $\tau$-sheaves, and the analog~\cite{Taguchi95b,Tam} of Tate's conjecture on endomorphisms. These local structures behave like in the classical case of abelian varieties, local shtukas playing the role of the $p$-divisible groups of the abelian varieties. The only difference is that $p$-divisible groups are only useful for abelian varieties in characteristic $p$, whereas the local shtukas at any place of $Q$ are important for the investigation of abelian $\tau$-sheaves and pure Anderson motives. One of the aims of this article is to demonstrate the utility of local shtukas. For instance we apply them in the computation of the Hasse invariants of $E$ in Theorem~\ref{THEOREM-2} and to reprove the standard fact that the set of morphisms between two pure Anderson motives is a projective $A$-module (\ref{ThmT.3}). Scattered in the text are several interesting examples displaying various phenomena 
% FIX VVVVVVVVVVVVVVVVVVVVVVVVVVVVVVVVVVVVVVVV
\forget{
(\ref{Ex1.8}, \ref{Ex1.18b}, \ref{Ex3.1}, \ref{Ex3.10c}, \ref{LAST-EXAMPLE}).
}
% FIX XXXXXXXXXXXXXXXXXXXXXXXXXXXXXXXXXXXXXXXX
(\ref{Ex1.8}, \ref{Ex1.8b}, \ref{Ex1.18b}, \ref{ExTRichter}, \ref{Ex3.1}, \ref{Ex3.10c}, \ref{LAST-EXAMPLE}, \ref{Ex3.15}).
% FIX ^^^^^^^^^^^^^^^^^^^^^^^^^^^^^^^^^^^^^^^^
% FIX VVVVVVVVVVVVVVVVVVVVVVVVVVVVVVVVVVVVVVVV
Note that this article has been split into two parts~\cite{BH_A,BH_B} for the purpose of publication.
% FIX ^^^^^^^^^^^^^^^^^^^^^^^^^^^^^^^^^^^^^^^^

% =============================================================================
\tableofcontents
% =============================================================================
%\include{diplom-1}
% =============================================================================

\subsection*{Notation}

In this article we denote by

\vspace{2mm}
\noindent
\begin{tabular}{@{}p{0.25\textwidth}@{}p{0.75\textwidth}@{}}
$\Fq$& the finite field with $q$ elements and characteristic $p$, \\
$C$& a smooth projective geometrically irreducible curve over $\Fq$, \\
$\infty\in C(\Fq)$& a fixed $\Fq$-rational point on $C$, \\
$\Ca = C\setminus\{\infty\}$&  \\
$A = \Gamma(\Ca,\OC)$& the ring of regular functions on $C$ outside $\infty$, \\
$Q = \Fq(C)$& the function field of $C$, 
\end{tabular}

\vspace{2mm}
\noindent 
thus we have $C' = \Spec A$ and $Q = \Quot A$. Furthermore, we denote by

\vspace{2mm}
\noindent
\begin{tabular}{@{}p{0.25\textwidth}@{}p{0.75\textwidth}@{}}
% FIX ^^^^^^^^^^^^^^^^^^^^^^^^^^^^^^^^^^^^^^^^^^^^^^
$Q_v$& the completion of $Q$ at the place $v\in C$, \\
$A_v$& the ring of integers in $Q_v$. For $v\ne\infty$ it is the completion of $A$ at $v$.\\
$\BF_v$ & the residue field of $A_v$. In particular $\BF_\infty\cong\Fq$.
\end{tabular}

\vspace{2mm}
\noindent
For a field $L$ containing $\Fq$ we write

\vspace{2mm}
\noindent
\begin{tabular}{@{}p{0.25\textwidth}@{}p{0.75\textwidth}@{}}
$C_L=C\times_{\Spec\Fq}\Spec L$,\\[1mm]
$A_L=A\otimes_\Fq L$,\\[1mm]
$Q_L=Q\otimes_\Fq L$,\\[1mm]
$A_{v,L}=A_v\wh\otimes_\Fq L$ & for the completion of $\O_{C_L}$ at the closed subscheme $v\times\Spec L$,\\[1mm]
$Q_{v,L}=A_{v,L}[\frac{1}{v}]$.& Note that this is not a field if $\BF_v\cap
L\supsetneq\Fq$. Nevertheless, it is always a finite product of fields.\\[1mm]
\end{tabular}
\begin{tabular}{@{}p{0.25\textwidth}@{}p{0.75\textwidth}@{}}
${\rm Frob}_q:L\to L$ & for the $q$-Frobenius endomorphism mapping $x$ to $x^q$,\\[1mm]
$\sigma = \id_C\times\Spec({\rm Frob}_q)$& for the endomorphism of $C_L$ which acts as the identity on the points and on $\O_C$ and as the $q$-Frobenius on $L$,\\
$\s$ & for the endomorphisms induced by $\sigma$ on all the above rings. For instance $\s(a\otimes b)=a\otimes b^q$ for $a\in A$ and $b\in L$.\\
$\s M=M\otimes_{A_L,\s}A_L$ & for an $A_L$-module $M$ and similarly for the other rings.
\end{tabular}

\forget{
\vspace{2mm}
\noindent
All schemes, as well as their products and morphisms, are supposed to be over $\Spec\Fq$. Let $S$ be a scheme. We denote by

\vspace{2mm}
\noindent
\begin{tabular}{@{}p{0.25\textwidth}@{}p{0.75\textwidth}@{}}
$\sigma_S: S\rightarrow S$& its $q$-Frobenius endomorphism which acts identically on points of $S$ and as the $q$-power map on the structure sheaf $\OS$, \\
$C_S = C\times S$& the base extension of $C$ from $\Spec\Fq$ to $S$, \\
$\sigma = \id_C\times\sigma_S$& the endomorphism on $C_S$ which acts identically on $C$ and as the $q$-Frobenius on $S$. 
\end{tabular}
}

\vspace{2mm}
\noindent
For a divisor $D$ on $C$ we denote by $\O_{C_L}(D)$ the invertible sheaf on $C_L$ whose sections $\varphi$ have divisor $(\varphi)\ge-D$. 
For a coherent sheaf $\F$ on $C_L$ we set $\F(D) := \F\otimes_{\O_{C_L}}\O_{C_L}(D)$. This notation applies in particular to the divisor $D = n\cdot\infty$ for $n\in\Z$.

%\pagebreak

\vspace{1cm}

% =============================================================================

\section{Pure Anderson Motives and Abelian $\tau$-Sheaves}

% =============================================================================

\subsection{Pure Anderson Motives}

Let $L$ be a field extension of $\Fq$ and fix an $\Fq$-homomorphism $c^\ast:A\to L$. Let $J\subset A_L$ be the ideal generated by $a\otimes 1-1\otimes c^\ast(a)$ for all $a\in A$, and let $r$ and $d$ be non-negative integers. Pure Anderson motives were introduced by G.\ Anderson~\cite{Anderson} and called \emph{pure $t$-motives} in the case where $A=\BF_p[t]$. In the slightly more general case we define:

\begin{Definition}[pure Anderson motives] \label{Def1'.1}
A \emph{pure Anderson motive $\ulM=(M,\tau)$ of rank $r$, dimension $d$, and characteristic $c^\ast$ over $L$} consists of a locally free $A_L$-module $M$ of rank $r$ and an injective $A_L$-module homomorphism $\tau:\s M\to M$ such that
\begin{suchthat}
\item 
the cokernel of $\t$ is an $L$-vector space of dimension $d$ and annihilated by $J^d$, and 
\item 
$M$ extends to a locally free sheaf $\CM$ of rank $r$ on $C_L$ such that for some positive integers $k,l$ the map $\tau^l:=\tau\circ\s(\tau)\circ\ldots\circ(\s)^{l-1}(\tau):(\s)^l M\to M$ induces an isomorphism $(\s)^l\CM_\infty\to\CM(k\cdot\infty)_\infty$ of the stalks at $\infty$. 
\end{suchthat}
We call $\chr:=\ker c^\ast\in \Spec A$ the \emph{characteristic point} of $\ulM$. We say that $\ulM$ has \emph{finite characteristic} (respectively \emph{generic characteristic}) if $\chr$ is a closed (respectively the generic) point. For $r>0$ we call $\weight(M,\tau):=\frac{k}{l}$ the \emph{weight} of $(M,\tau)$.
\end{Definition}

\noindent {\it Remark.} \nopagebreak 
Phrased in the language of modules our definition of purity is equivalent to the following due to Anderson~\cite[1.9]{Anderson}. Let $z$ be a uniformizing parameter of $A_{\infty,L}$. The Anderson motive is pure if and only if there exists an $A_{\infty,L}$-lattice $\hat\CM_\infty$ inside $M\otimes_{A_L}Q_{\infty,L}$ and positive integers $k,l$ such that $z^k\tau^l$ induces an isomorphism $(\s)^l\hat\CM_\infty\to\hat\CM_\infty$. This follows from the fact that $\hat\CM_\infty$ determines a unique extension $\CM$ of $M$ as above.

\begin{Proposition}\label{Prop1'.1b}
If $\ulM$ is a pure Anderson motive of rank $r>0$ then $\weight\ulM=\frac{d}{r}$\,. In particular $\dim\ulM>0$.
\end{Proposition}

\begin{proof}
Using \ref{DEGREE-LEMMA} below we compute
\[
k\,r \; = \; \deg \CM(k\cdot\infty)-\deg\CM
\; = \; \deg\CM(k\cdot\infty)-\deg(\s)^l\CM
\; = \;\dim_L\coker(\tau^l) \;=\; l\,d\,.\qedhere
\]
\end{proof}

\begin{Lemma}\label{DEGREE-LEMMA}
Let\/ $\G$ be a coherent sheaf on $C_L$. Then $\deg\s\G = \deg\G$.
\end{Lemma}

\begin{proof}
Let $pr:C_L\to \Spec L$ be the projection onto the second factor. If $\G=\OCL^{\!\oplus n}$ for some $n\in\N$, then $\s\G = \OCL^{\!\oplus n}$. If $\G=\OCL(D)$ for some divisor $D$ on $C$, then $\s\G=\OCL(\sigma^{-1}D)=\OCL(D)$. If $\G$ is a torsion sheaf we get
\[
\deg\s\G \,=\, \dim_L pr_\ast\s\G \,=\, \dim_L \s pr_\ast\G \,=\, \dim_L pr_\ast\G \,=\, \deg\G\,.
\]
Let now $\G$ be a locally free sheaf of rank $n$. Choose an embedding $f:\,\G\rightarrow\OCL(D)^{\oplus n}$ for some divisor $D$ on $C$ with $\coker f$ being a torsion sheaf. Since $\sigma=\id_C\times\sigma_L$ is flat being the base change of the flat morphism $\sigma_L:\Spec L\to\Spec L$, we have 
\[
\begin{CD} 
{0} @>>> {\G}   @>{f   }>> {\OCL(D)^{\oplus n}}   @>>> {\coker f}   @>>> {0} \\[1ex]
{0} @>>> {\s\G} @>{\s f}>> {\s\OCL(D)^{\oplus n}} @>>> {\s\coker f} @>>> {0} \\[1ex]
\end{CD}
\]
and therefore $\deg\s\G = \deg\G$ due to the additivity of the degree in exact sequences. Finally, if $\G$ is an arbitrary coherent sheaf, then
\[
\exact{0}{}{\G'}{}{\G}{}{\G''}{}{0}
\]
for some torsion sheaf $\G'$ and some locally free coherent sheaf $\G''$ because this sequence exists locally due to the fact that all local rings are principal ideal domains. Thus $\deg\s\G = \deg\G$, as desired.
\end{proof}

\begin{Proposition} \label{PropX.3}
If $(M,\tau)$ is a pure Anderson motive over $L$ then one can find an extension $\CM$ as above with $k$ and $l$ relatively prime.
\end{Proposition}

\begin{proof}
We let $z$ be a uniformizing parameter at $\infty$ and write $\frac{d}{r}=\frac{k}{l}$ with $k,l$ relatively prime positive
integers. Since $(M,\t)$ is pure it extends to a locally free sheaf $\CM'$ on $C_L$ on which $z^{k'}\t^{l'}$ is an isomorphism locally at $\infty$ for some positive integers $k',l'$ with $\frac{k'}{l'}=\frac{d}{r}=\frac{k}{l}$. We modify $\CM'$ to a locally free sheaf $\CM$ on $C_L$ by changing its stalk $\CM'_\infty$ at $\infty$ to
\[
\CM_\infty\;:=\;\sum_{j=0}^{\frac{l'}{l}-1}z^{kj}\t^{lj}\bigl((\s)^{lj}\CM'_\infty\bigr)\,.
\]
Then $M=\Gamma(C_L\setminus\{\infty\},\CM)$ and $z^{k}\tau^{l}:(\s)^{l}\CM_\infty\isoto\CM_\infty$ is an isomorphism at $\infty$ as desired. 
\forget{
We write $\frac{d}{r}=\frac{k'}{l'}$ with $k',l'$ relatively prime positive integers and set
\[
\ulN_\infty := (N_\infty,\tau):=(M,\tau)\otimes_{A_L}Q_{\infty,L}\,.
\]
In the terminology of \cite[Definition~7.3]{Hl} the pair $(N_\infty,\tau)$ is a \emph{Dieudonn\'e $\BF_q\dpl z\dpr$-module} over $L$ (also called a \emph{local isoshtuka} in Section~\ref{SectLS} below). Due to the purity condition there exists an $A_{\infty,L}$-lattice $\hat\CM_\infty$ inside $N_\infty$ and positive integers $k,l$ with $\frac{k}{l}=\frac{d}{r}$ such that $z^{k}\tau^{l}$ is an isomorphism $(\s)^l\hat\CM_\infty\to\hat\CM_\infty$. Therefore $\ulN_\infty$ is isoclinic of slope $\frac{k}{l}=\frac{k'}{l'}$. Consider the $A_{\infty,L}$-lattice
\[
\hat\CM'_\infty\;:=\;\sum_{j=0}^{\frac{l}{l'}-1}z^{k'j}\t^{l'j}\bigl((\s)^{l'j}\hat\CM_\infty\bigr)
\]
inside $\ulN_\infty$ on which $z^{k'}\tau^{l'}:(\s)^{l'}\hat\CM'_\infty\isoto\hat\CM'_\infty$ is an isomorphism. We glue $M$ with $\hat\CM'_\infty$ to the locally free sheaf $\CM'$ on $C_L$ on which now $(\tau')^{l'}:(\s)^{l'}\CM'\to \CM(k'\cdot\infty)$ is an isomorphism near $\infty$. 
}
\end{proof}

\begin{Definition} (compare \cite[4.5]{PT})\label{Def1'.2}
\begin{suchthat}
\item 
A \emph{morphism} $f:(M,\tau)\to (M',\tau')$ between Anderson motives of the same characteristic $c^\ast$ is a homomorphism $f:M\to M'$ of $A_L$-modules which satisfies $f\circ\tau=\tau'\circ\s(f)$.
\item 
If $f:\ulM\to\ulM'$ is surjective, $\ulM'$ is called a \emph{quotient} (or \emph{factor}) \emph{motive} of $\ulM$.
\item 
A morphism $f:\ulM\to \ulM'$ is called an \emph{isogeny} if $f$ is injective with torsion cokernel.
\item 
An isogeny is called \emph{separable} (respectively \emph{purely inseparable}) if the induced homomorphism $\t:\s\coker f\to\coker f$ is an isomorphism (respectively is \emph{nilpotent}, that is, if $\t^n=0$ for some $n$).
\end{suchthat}
\end{Definition}

\bigskip

\noindent {\it Remark.} \es
1. The set $\Hom(\ulM,\ulM')$ of morphisms is an $A$-module and $\End(\ulM)$ is an $A$-algebra. They are projective $A$-modules of rank $\le rr'$. This will be proved in Theorem~\ref{ThmT.3}.

\smallskip
\noindent
2. One has $\Hom(\ulM,\ulM')=\{0\}$ if $\ulM$ and $\ulM'$ are pure Anderson motives of different weights, justifying the terminology \emph{pure}. This can be derived from the Dieudonn\'e-Manin type classification \cite[Appendix B]{Laumon} of the local $\sigma$-isoshtuka $\ulM_\infty(\ulM):=\ulM\otimes_{A_L}Q_{\infty,L}$ of $\ulM$ at $\infty$; see Section~\ref{SectLS}. However, we will give a more elementary proof in Corollary~\ref{Cor2.9b} below.

\smallskip
\noindent
3. We will prove in Corollary~\ref{Cor1.11b} below that the cokernel of an isogeny is in fact annihilated by an element $a\in A$ (as opposed to $a\in A_L$). This was independently observed by N.~Stalder~\cite{Stalder} and also holds for non-pure Anderson motives.

\begin{Proposition}\label{Prop1.5b}
Let $(M,\t)$ be a pure Anderson motive and let $K$ be a finite torsion $A_L$-module equipped with an $A_L$-homomorphism $\t_K:\s K\to K$ such that both $\ker\t_K$ and $\coker\t_K$ are annihilated by a power of $J=(a\otimes1-1\otimes c^\ast(a):a\in A)\subset A_L$. Let further $\rho:M\shortonto K$ be a surjective morphism of $A_L$-modules with $\t_K\circ\s\rho=\rho\circ\t$. Then $(M',\t'):=(\ker\rho,\t|_{\s M'})$ is again a pure Anderson motive of the same rank and dimension and the inclusion $f:(M',\t')\to(M,\t)$ is an isogeny with $\coker f=(K,\t_K)$.
\end{Proposition}

\begin{proof}
Consider the diagram in which the bottom row is obtained from the snake lemma
\[
\xymatrix {& 0 \ar[r] & \s M' \ar[r]^{\s f} \ar@{^{ (}->}[d]_{\t'} & \s M \ar[r]^{\s\rho} \ar@{^{ (}->}[d]_{\t} & \s K \ar[r] \ar[d]_{\t_K} & 0 \\
& 0 \ar[r] & M' \ar[r]^f \ar@{->>}[d] & M \ar[r]^\rho \ar@{->>}[d] & K \ar[r] \ar@{->>}[d] & 0 \\
0 \ar[r] & \ker\t_K \ar[r] & \coker\t' \ar[r] & \coker\t \ar[r] & \coker\t_K \ar[r] & 0\;.}
\]
If follows that $\dim_L\coker\t'=\dim_L\coker\t=d$ and that also $\coker\t'$ is annihilated by $J^d$. The purity follows from the fact that one can extend $f$ to an isomorphism $\CM'_\infty\to\CM_\infty$ of the stalks at infinity.
\end{proof}

\noindent
{\it Remark.}
Note that without the requirement that a power of $J$ annihilates $\ker\tau_K$ and $\coker\tau_K$ the assertion  of the proposition is false as one can see from $A=\BF_q[t],\,M=A_L,\,\tau=t\otimes1-1\otimes c^\ast(t),\,K=\coker\tau,\,\tau_K=0$, when $c^\ast(t)^q\ne c^\ast(t)$.

\begin{Lemma}\label{Lemma1.5d}
Let $K$ be a finite torsion $A_L$-module and let $\t:\s K\to K$ be a morphism of $A_L$-modules. Then $\ulK=(K,\t)$ is an extension
\[
0\longto \ulK^\et\xrightarrow{\es f\;} \ulK\xrightarrow{\es g\;}\ulK^\nil\longto 0
\]
of $\ulK^\nil=(K^\nil,\t^\nil:\s K^\nil\to K^\nil)$ by $\ulK^\et=(K^\et,\t^\et:\s K^\et\to K^\et)$ where $\t^\nil$ is nilpotent and $\t^\et$ is an isomorphism, satisfying $f\circ\t^\et=\t\circ\s f$ and $g\circ\t=\t^\nil\circ\s g$. Moreover, if the base field $L$ is perfect the extension splits canonically.
\end{Lemma}

\begin{proof}
This was proved by Laumon~\cite[B.3.10]{Laumon}. He takes $ K^\et:=\bigcap_{n\ge1}\im\t^n$. If $L$ is perfect, $K^\et$ has the natural complement $\bigcup_{n\ge1}(\s)^{-n}(\ker\t^n)$ which is isomorphic to $\ulK^\nil=\ulK/\ulK^\et$.
\end{proof}

\begin{Proposition}\label{Prop1.5c}
Every isogeny $f:\ulM\to \ulM'$ can be factored $\ulM\xrightarrow{\es f_{\rm sep}\;}\ulM''\xrightarrow{\es f_{\rm insep}\;}\ulM'$ into a separable isogeny $f_{\rm sep}$ followed by a purely inseparable isogeny $f_{\rm insep}$. If the base field is perfect there exists also a (different) factorization $f=f'_{\rm sep}\circ f'_{\rm insep}$ as a purely inseparable isogeny followed by a separable one.
\end{Proposition}

\begin{proof}
Let $K:=\coker f$ and let $\t_K:\s K\to K$ be the induced morphism. By Lemma~\ref{Lemma1.5d} there is a surjective morphism $\rho:\ulM'\to\ulK\to\ulK^\nil$ and we define $\ulM''$ as the kernel of $\rho$. It is a pure Anderson motive by Proposition~\ref{Prop1.5b}. Clearly $f$ factors through $\ulM''$ and the isogeny $\ulM\to \ulM''$ has $\ulK^\et$ as cokernel, thus is separable. If $L$ is perfect we use the surjective morphism $\rho:\ulM'\to\ulK\to\ulK^\et$ instead.
\end{proof}

% =============================================================================

\bigskip
\subsection{Definition of abelian $\tau$-sheaves}

Since the purity condition 2 of Definition~\ref{Def1'.1} does not behave well
in families one has to rigidify $\ulM$ at $\infty$ in order to get moduli
spaces for pure Anderson motives. This was done in \cite{Hl}, where the rigidified objects are called \emph{abelian sheaves}. Over a field their definition is as follows.
Let $L\supset\Fq$ be a field and fix a morphism $c: \Spec L\rightarrow C$. Let $\I$ be the ideal sheaf on $C_L$ of the graph $\Graph(c)$ of $c$. Let $r$ and $d$ be non-negative integers.

\begin{Definition}[Abelian $\tau$-sheaf]\label{Def1.1}
An \emph{abelian $\tau$-sheaf $\FF = (\F_i,\P_i,\t_i)$ of rank $r$, dimension $d$ and characteristic $c$ over $L$} is a collection of locally free sheaves $\F_i$ on $C_L$ of rank $r$ together with injective morphisms $\P_i,\t_i$ of\/ $\O_{C_L}$-modules $(i\in\Z)$ of the form
\[
\begin{CD}
\cdots @>>> 
\F_{i-1} @>{\P_{i-1}}>> 
\F_{i}   @>{\P_{i}}>> 
\F_{i+1} @>{\P_{i+1}}>> 
\cdots 
\\ & & 
@AA{\t_{i-2}}A 
@AA{\t_{i-1}}A 
@AA{\t_{i}}A 
\\
\cdots @>>> 
\s\F_{i-2} @>{\s\P_{i-2}}>> 
\s\F_{i-1} @>{\s\P_{i-1}}>> 
\s\F_{i}   @>{\s\P_{i}}>> 
\cdots 
\\
\end{CD}
\]
subject to the following conditions:
\begin{suchthat}
\item the above diagram is commutative,
\item there exist integers\/ $k,l>0$ with\/ $ld=kr$ such that the morphism $\P_{i+l-1}\circ\cdots\circ\P_i$ identifies $\F_i$ with the subsheaf\/ $\F_{i+l}(-k\cdot\infty)$ of\/ $\F_{i+l}$ for all $i\in\Z$,
\item $\coker\P_i$ considered as an $L$-vector space has dimension $d$ for all $i\in\Z$,
\item $\coker\t_i$ is annihilated by $\I^d$ and as an $L$-vector space has dimension $d$ and for all $i\in\Z$.
\end{suchthat}
We call $\chr:=c(\Spec L)\in C$ the \emph{characteristic point} (or \emph{place}) and say that $\FF$ has \emph{finite} (respectively \emph{generic}) \emph{characteristic} if $\chr$ is a closed (respectively the generic) point. 
\end{Definition}

\begin{Remark} \label{Rem2.2}
1. By the second condition $\coker\P_i$ is only supported at $\infty$. Moreover, the periodicity condition implies $\F_{i+nl} = \F_i(nk\cdot\infty)$ and thus $\t_{i+nl}=\t_i\otimesidOCL{nk\cdot\infty}$ for all $n\in\Z$.

\smallskip
\noindent
2. The condition ``annihilated by $\I$'' in 4 can equivalently be reduced to ``supported on the graph of $c$ '', since the local ring of $C_L$ at the graph of $c$ is a principal ideal domain and the $d$-th power of a generator of $\I$ annihilates the $d$-dimensional $L$-vector space $\coker\t_i$.

\smallskip
\noindent
3. Trivially, $r=0$ implies $d=0$ since in this case we have all $\F_i=0$. Due to the second condition, the converse is also true because $d=0$ implies $r = \frac{l}{k}\,d = 0$ since the existence of such $k,l\not= 0$ is required. Without this the converse would in general not be true, because for example $\FF = (\OCS,\id_\OCS,\id_\OCS)$ has $\coker\id_\OCS = 0$ and therefore $d=0$, but $r=1$. This justifies the demand of the existence of such $k,l\not= 0$ since we do not want to consider the \quotes{degenerate} case $r>0$, $d=0$.

The case $r=0$, $d=0$ however is desired because it allows the \emph{zero sheaf} $\ZZ := (0,0,0)$ to be an abelian $\tau$-sheaf of rank $0$ and dimension $0$. Trivially, the zero sheaf satisfies the second condition for \emph{all} pairs $k,l > 0$.

\smallskip
\noindent
4. For $\FF\not= \ZZ$ one can ask whether the second condition is satisfied by the pair $k,l > 0$ with $ld = kr$ and $k,l$ relatively prime. This was required in the definition of abelian $\tau$-sheaves in \cite{Hl}. We will call those $\FF$ \emph{abelian $\tau$-sheaves with $k,l$ relatively prime}.
As a convention, an abelian $\tau$-sheaf $\FF$ without further specifications comes with all its parameters $\F_i,\P_i,\t_i$ $(i\in\Z)$ and $r,d,k,l$ with $k,l$ always chosen to be minimal. Similarly $\FF'$ carries a prime on its parameters, $\TFF$ a tilde on them, and so on. Note that the characteristic $c$ is fixed.

\smallskip
\noindent
5. Abelian $\tau$-sheaves of dimension $d=1$ are called \emph{elliptic sheaves} and were studied by Drinfeld~\cite{Drinfeld3}, Blum-Stuhler~\cite{BS} and others. The category of elliptic sheaves with $(k,l)=1$ over $L$ of rank $r$ with $\deg\F_0=1-r$ and whose characteristic does not meet $\infty$, that is, $\im(c)\subset C\setminus\{\infty\}$, is anti-equivalent to the category of Drinfeld-$A$-modules of rank $r$ over $L$, see \cite[Theorem~3.2.1]{BS} and Example~\ref{Ex1.8} below.
\end{Remark}

\smallskip

\begin{Definition}\label{DEFINITION-WEIGHT}
Let\/ $\FF$ be an abelian $\tau$-sheaf of rank $r$, dimension $d$ and characteristic $c$ over $S$. We set
\[
\weight(\FF) \;:=\; 
\left\{\;
\begin{array}{cl} 
\frac{d}{r} & \mbox{if\/ }\FF\not=\ZZ \\[1ex] 
0           & \mbox{if\/ }\FF=\ZZ 
\end{array} 
\;\right\}
\;\in\Q\,.
\]
We call\/ $\weight(\FF)$ the \emph{weight of $\FF$}.
\end{Definition}

\begin{Example} \label{Ex1.8}
Let  $C=\PP^1_\Fq$, $A=\FqItI$. Then $\CL=\PP^1_L$. Let $c: \Spec L\rightarrow \Spec\FqItI=\AA^1_\Fq$ such that $c^\ast: \FqItI\rightarrow L$ maps $t$ to $c^\ast(t)=:\vartheta$. Let $\O$ denote the structure sheaf of $\PP^1_L$ and let $a\in L$. Now consider the following diagram
\[
\begin{CD}
\cdots @>>> 
\O\oplus\O @>{\matr{1}{0}{0}{1}}>{\P_0}> 
\O\oplus\O(1\!\cdot\!\infty) @>{\matr{1}{0}{0}{1}}>{\P_1}> 
\O(1\!\cdot\!\infty)\oplus\O(1\!\cdot\!\infty) @>{\matr{1}{0}{0}{1}}>{\P_2}>
\cdots 
\\ & & 
@A{\textstyle\t_{\scriptscriptstyle-1}}\;AA
@A{\textstyle\t_{\scriptscriptstyle 0}}\;A\;{\matr{a}{t-\vartheta}{1}{0}}A
@A{\textstyle\t_{\scriptscriptstyle 1}}\;A\;{\matr{a}{t-\vartheta}{1}{0}}A 
\\
& &
\cdots @>{\matr{1}{0}{0}{1}}>{\s\P_{-1}}> 
\s(\O\oplus\O) @>{\matr{1}{0}{0}{1}}>{\s\P_0}> 
\s(\O\oplus\O(1\!\cdot\!\infty)) @>{\matr{1}{0}{0}{1}}>{\s\P_1}> 
\cdots 
\\
\end{CD}
\]
where the vectors in $\O\oplus\O$ are considered as column vectors.
This gives an example of an abelian $\tau$-sheaf $\FF$ of rank $2$, dimension $1$ with $(k,l)=1$, and characteristic $c$ over $\Spec L$ with $\weight(\FF)=\frac{1}{2}$.

Since the dimension is $1$, this abelian $\tau$-sheaf is an elliptic sheaf and comes from a Drinfeld module which can be recovered as follows; see \cite{BS}. Let $M := \Gamma(\AA^1_L,\O\oplus\O) = L[t]\oplus L[t]$ and let $\t := \P_0^{-1}\circ\t_0$. Since $\vect{1}{0}=\t\vect{0}{1}$, we have $M=L\{\t\}\cdot\vect{0}{1}$, and we calculate
\[
\t^2\tvect{0}{1} \,=\, \tvect{a}{t-\vartheta} \,=\, a\cdot\t\tvect{0}{1} + (t-\vartheta)\cdot\tvect{0}{1}
\quad \text{and} \quad
t\cdot\tvect{0}{1} = (\vartheta - a\t + \t^2)\tvect{0}{1}
\]
Let $\varphi: \FqItI\rightarrow L\{\t\}$ be the ring morphism mapping $t\mapsto \vartheta - a\t + \t^2$. Then we have back the Drinfeld Module $\varphi$ of rank 2 over $L$ which induces the abelian $\tau$-sheaf $\FF$.
\end{Example}

\begin{Example} \label{Ex1.8b}
We give another example which does not come from Drinfeld modules. Let $C=\PP^1_\Fq$ and $A=\Fq[t]$. Let 
\[
S=\Spec \Fq[\zeta,\alpha,\beta,\gamma,\delta]\,/\,\bigl((\alpha+\delta+2\zeta)^2\,,\, \alpha\delta-\beta\gamma-\zeta^2\bigr)
\]
and $c:S\to C\setminus V(t)\subset C$ be given by $c^\ast(\frac{1}{t})=\zeta$. Let $C_S:=C\times_{\Spec\Fq}S$, then
\[
\begin{CD}
\cdots @>>> 
\O_{C_S}^{\oplus2} @>{\quad\matr{1}{0}{0}{1}\quad}>{\P_0}> 
\O_{C_S}(1\!\cdot\!\infty)^{\oplus2} @>{\qquad\matr{1}{0}{0}{1}\qquad}>{\P_1}> 
\O_{C_S}(2\!\cdot\!\infty)^{\oplus2} @>{\quad\matr{1}{0}{0}{1}\quad}>{\P_2}>
\cdots 
\\ & & 
@A{\textstyle\t_{\scriptscriptstyle-1}}\;AA
@A{\textstyle\t_{\scriptscriptstyle 0}}\;A\;{\matr{1+\alpha t}{\gamma t}{\beta t}{1+\delta t}}A
@A{\textstyle\t_{\scriptscriptstyle 1}}\;A\;{\matr{1+\alpha t}{\gamma t}{\beta t}{1+\delta t}}A 
\\
& &
\cdots @>{\quad\matr{1}{0}{0}{1}\quad}>{\s\P_{-1}}> 
\s\O_{C_S}^{\oplus2} @>{\qquad\matr{1}{0}{0}{1}\qquad}>{\s\P_0}> 
\s\O_{C_S}(1\!\cdot\!\infty)^{\oplus2} @>{\quad\matr{1}{0}{0}{1}\quad}>{\s\P_1}> 
\cdots 
\\
\end{CD}
\]
is an abelian $\tau$-sheaf over $S$ of rank and dimension $2$ with $k=l=1$ since $(1-\zeta t)^2\cdot\coker\tau=(0)$. In fact $S$ is the (representable part of the) moduli space of abelian $\tau$-sheaves of rank and dimension $2$ with $(k,l)=1$ together with a level structure $\eta$ at $V(t)$ for which $(\F_0,\eta)$ is stable of degree zero. See~\cite[\S 4]{Hl} for the precise meaning of these terms, but note that in loc.\ cit.\ the exponent $2$ in $(\alpha+\delta+2\zeta)^2$ erroneously was missing, as was pointed out to us by M.~Molz. This illustrates the fact that abelian $\tau$-sheaves possess nice moduli spaces.
\end{Example}

\begin{Proposition}\label{PROP.2}
Let $\FF$ be an abelian $\tau$-sheaf and let $D$ be a divisor on $C$. Then the collection $\FF(D) := (\F_i(D),\P_i\otimesidOCL{D},\t_i\otimesidOCL{D})$ is an abelian $\tau$-sheaf of the same rank and dimension as $\FF$.
\end{Proposition}

\begin{proof}
Since the functor $\otimes_{\O_{C_L}}\O_{C_L}(D)$ is exact the proof is straightforward ones one notes that $\s(\F_i(D))=(\s\F_i)(D)$ because the divisor $D$ is $\sigma$-invariant.
\end{proof}

\begin{Definition}
Let $\FF$ be an abelian $\tau$-sheaf and let $n\in\Z$. We denote by
\[
\FF\II{n} := (\F_{i+n},\P_{i+n},\t_{i+n})
\]
the \emph{$n$-shifted abelian $\tau$-sheaf of\/ $\FF$} whose collection of\/ $\F$'s, $\P$'s and $\t$'s is just shifted by\/ $n$.
\end{Definition}

\noindent
It is trivial to see that $\FF\II{0}=\FF$ and that $\FF\II{n}$ is an abelian $\tau$-sheaf of the same rank and dimension as $\FF$ for every $n\in\Z$. Moreover, $\FF[l]\cong\FF(k\cdot\infty)$.

\bigskip

Next we come to the definition of morphisms in the category of abelian $\tau$-sheaves.

\begin{Definition}
A \emph{morphism} $f$ between two abelian $\tau$-sheaves $\FF = (\F_i,\P_i,\t_i)$ and $\FF' = (\F'_i,\P'_i,\t'_i)$ of the same characteristic $c:\Spec L\to C$ is a collection of morphisms $f_i: \F_i \rightarrow \F'_i$ $(i\in\Z)$ which commute with the $\P$'s and the $\t$'s, that is, $f_{i+1}\circ\P_i=\P'_i\circ f_i$ and $f_{i+1}\circ\t_i=\t'_i\circ\s f_i$. We denote the set of morphisms between $\FF$ and $\FF'$ by $\Hom(\FF,\FF')$. It is an $\Fq$-vector space.
\end{Definition}

For example, the collection of morphisms $(\P_i):\, \FF\rightarrow\FF\II{1}$ defines a morphism between the abelian $\tau$-sheaves $\FF$ and $\FF\II{1}$. 

\begin{Definition}
Let $\FF$ and $\FF'$ be abelian $\tau$-sheaves and let $f \in \Hom(\FF,\FF')$ be a morphism. 
\begin{suchthat}
\item $f$ is called \emph{injective}, if\/ $f_i$ is injective for all\/ $i\in\Z$.
\item $f$ is called \emph{surjective}, if\/ $f_i$ is surjective for all\/ $i\in\Z$.
\item $f$ is called an \emph{isomorphism}, if\/ $f$ is injective and surjective.
\end{suchthat}
We call $\FF$ an \emph{abelian quotient} (or \emph{factor}) \emph{$\tau$-sheaf} of\/ $\FF'$, if there is a surjective morphism from $\FF'$ onto $\FF$.
\end{Definition}

Abelian $\tau$-sheaves are pure in the following sense.

\begin{Proposition}\label{PROP.1}
Let $\FF$ and $\FF'$ be abelian $\tau$-sheaves. If\/ $\Hom(\FF,\FF')\not=\{0\}$, then\/ $\weight(\FF)=\weight(\FF')$. 
\end{Proposition}

\begin{proof}
Let $0\not=f\in\Hom(\F,\F')$ and let $i\in\Z$. Consider the sheaf $\HOM(\F^{\,}_i,\F'_i) = \F'_i\otimes_\OCL\dual{\F}_i$ and the set of all its locally free subsheaves $\M\subset\F'_i\otimes_\OCL\dual{\F}_i$. Then the set of their degrees $\deg\M$ is bounded above, say with upper bound $B$ by \cite[Lemma 1.I.3]{Se}.

Suppose $d'r<dr'$. Choose $n\in\Z$ with $ll'|nrr'$ such that $B+n(d'r-dr')<0$. Let $\P$ and $\P'$ be the identifying morphisms $\P_{i+nrr'-1}\circ\cdots\circ\P_i: \F_i \cong \F_{i+nrr'}(-ndr'\cdot\infty)$ and $\P'_{i+nrr'-1}\circ\cdots\circ\P'_i: \F'_i \cong \F'_{i+nrr'}(-nd'r\cdot\infty)$, respectively. Consider the following diagram
\[
\xymatrix{
\F_i\; \ar[r]^{\P\;\;}\ar[d]^{f_i}& \;\F_{i+nrr'} \makebox[0pt][l]{ $= \;\F_i(ndr'\cdot\infty)\;$} \ar[d]^{f_{i+nrr'}} &\qquad \\
\F'_i\; \ar[r]^{\P'\;\;}& \;\F'_{i+nrr'} \makebox[0pt][l]{ $= \;\F'_i(nd'r\cdot\infty)\;$\,.} &\qquad \\
}
\]
With $m:=n(d'r-dr')<0$, we conclude that $\HOM(\F_{i+nrr'},\F'_{i+nrr'}) \;=\; 
(\F'_i\otimes_\OCL\dual{\F}_i)(m\cdot\infty)$.
Now, considering the injective map $\varphi: \OCL\rightarrow(\F'_i\otimes_\OCL\dual{\F}_i)(m\cdot\infty)$, $1\mapsto f_{i+nrr'}$ we get a non-zero locally free subsheaf $\im\varphi=\M(m\cdot\infty)$ which is isomorphic to $\OCL$ with $\M\subset\F'_i\otimes_\OCL\dual{\F}_i$. Therefore 
\[
0\;=\;\deg\M(m\cdot\infty) \;=\; \deg\M+m\cdot\rank\M \;\le\; B + m \;<\; 0\,.
\]
This is a contradiction and shows $d'r\ge dr'$. The converse $d'r\le dr'$ follows analogously.
\end{proof}

\noindent
{\it Remark.} This result can also be proved by considering the local isoshtukas at $\infty$ of $\FF,\FF'$ (see Section~\ref{SectLS}) and using the Dieudonn\'e-Manin theory \cite[Appendix B]{Laumon} for local isoshtukas.

% -----------------------------------------------------------------------------

\bigskip

\subsection{Relation between pure Anderson motives and abelian $\tau$-sheaves} \label{SectRelation}

If $\FF=(\F_i,\P_i,\t_i)$ is an abelian $\tau$-sheaf of rank $r$, dimension $d$, and characteristic $c:\Spec L\to C$ with characteristic place $\chr=\im c\ne\infty$ then
\begin{equation}\label{Eq1.1}
\ulM(\FF)\;:=\;(M,\tau)\;:=\;\Bigl(\Gamma(C_L\setminus\{\infty\},\F_0)\,,\,\P_0^{-1}\circ\t_0\Bigr)
\end{equation}
is a pure Anderson motive of the same rank and dimension and of characteristic $c^\ast:A\to L$. We can take $\CM:=\F_0$ as the extension of $M$ to all of $C_L$. Conversely we have the following result.

\begin{Theorem}\label{PropX.1}
\begin{enumerate}
\item 
Let $(M,\tau)$ be a pure Anderson motive of rank $r$, dimension $d$, and characteristic $c^\ast:A\to L$ over $L$. Then $(M,\tau)=\ulM(\FF)$ for an abelian $\tau$-sheaf $\FF$ over $L$ of same rank and dimension with characteristic $c:=\Spec c^\ast:\Spec L\to \Spec A\subset C$. One can even find the abelian $\tau$-sheaf $\FF$ with $k,l$ relatively prime.
\item 
Let $\FF$ and $\FF'$ be abelian $\tau$-sheaves of the same weight and let $f_0:\ulM(\FF)\to \ulM(\FF')$ be a morphism. Then there exists an integer $m$ such that $f_0$ comes from a morphism $f:\FF\to\FF'(m\cdot\infty)$ as $f_0=\ulM(f)$.
\end{enumerate}
\end{Theorem}

\begin{proof}
1. Let $\CM$ be a locally free sheaf on $C_L$ with $M=\Gamma(C_L\setminus\{\infty\},\CM)$ as in Definition \ref{Def1'.1}. Let $k,l$ be positive integers with isomorphism
\[
\tau^l:(\s)^l\CM_\infty\isoto\CM(k\cdot\infty)_\infty\,.
\]
By Proposition~\ref{PropX.3} we may assume $(k,l)=1$.

For $i=0,\ldots,l$ let $\F_i$ be the locally free sheaf of rank $r$ on $C_L$ which coincides with $\CM$ on $C_L\setminus\{\infty\}$ and whose stalk at $\infty$ is the sum $(\im\t^i+\ldots+\im\t^{i+l-1})_\infty$ inside $\CM(2k\cdot\infty)_\infty$. Then $\t$ defines homomorphisms $\t_i:\s\F_i\to\F_{i+1}$ for $0\le i<l$ because $\s\im\t^i=\im\s\t^i$ due to the flatness of $\s=(\id_C\times\sigma_L)^\ast$. Since $\CM_\infty\subset\CM(k\cdot\infty)_\infty=(\im\t^l)_\infty$ there are natural inclusions $\P_i:\F_i\to\F_{i+1}$ satisfying $\P_{i+1}\circ\t_i=\t_{i+1}\circ\s\P_i$ and $\im(\P_{l-1}\circ\ldots\circ\P_0)=\F_l(-k\cdot\infty)\subset\F_l$. We now set $\F_{i+nl}:=\F_i(kn\cdot\infty), \P_{i+nl}:=\P_i\otimes\id,\t_{i+nl}:=\t_i\otimes\id$ for $0\le i<l$ and $n\in\Z$. Clearly $\coker\t_i$ is supported on $\Graph(c)$ for all $i$ and isomorphic to $\coker\t$ which is an $L$-vector space of dimension $d$. We compute 
\[
\dim_L\coker\P_i\;=\;\deg\F_{i+1}-\deg\F_i\;=\;\deg\F_{i+1}-\deg\s\F_i\;=\;\dim_L\coker\t_i
\]
for all $i$. Hence $\FF=(\F_i,\P_i,\t_i)$ is an abelian $\tau$-sheaf over $L$ with $\ulM(\FF)=(M,\t)$.

\medskip
\noindent 
2. Let $l$ be an integer satisfying condition 2 of Definition~\ref{Def1.1}. For $0<i$ set 
\[
f_i:=\P'_{i-1}\circ\ldots\circ\P'_0\circ f_0\circ\P_0^{-1}\circ\ldots\circ\P_{i-1}^{-1}
\]
and similarly for $i<0$. Since the $\P_j,\P'_j$ are isomorphisms outside $\infty$ there exists an integer $m$ such that $f_i$ is a morphism $f_i:\F_i\to\F_i'(m\cdot\infty)$ for all $0\le i\le l$. Now the periodicity condition 2 of Definition~\ref{Def1.1} shows that the latter is a morphism for all $i\in\Z$. Finally $(\P'_0)^{-1}\t_0\,\s(f_0)=f_0\P_0^{-1}\t_0$ implies that $f=(f_i)_i:\FF\to\FF(m\cdot\infty)$ is a morphism with $\ulM(f)=f_0$ as desired.
\end{proof}

\medskip

The aforementioned relation can more generally be described by the following terminology. Let $\Spec\TA\subset C$ be an affine open subscheme.

\begin{Definition}\label{Def1.16}
A\/ \emph{$\t$-module on $\TA$ over $L$ of rank $r$} is a pair $\ulM=(M,\t)$, where
\begin{suchthat}
\item $M$ is a locally free $\TA\otimes_\Fq L$-module of rank $r$,
\item $\t: \s M\rightarrow M$ is injective.
\end{suchthat}
A\/ \emph{morphism} between $(M,\t)$ and $(M',\t')$ is a homomorphism $f:\, M\rightarrow M'$ of $\TA\otimes_\Fq L$-modules which respects $\t'\circ \s f = f\circ\t$. We denote the set of morphisms between $\ulM$ and $\ulM'$ by $\Hom(\ulM,\ulM')$.
\end{Definition}

Let $\FF$ be an abelian $\tau$-sheaf. 
Consider a finite closed subset $D\subset C$ such that either $\infty\in D$ or there exists a uniformizing parameter $z$ at infinity inside $\TA:=\Gamma(C\setminus D,\O_C)$. Note that by enlarging $D$ it will always be possible to find such a $z\in\TA$ in the case $\infty\not\in D$. 

If $\infty\in D$ we have by the $\P$'s a chain of isomorphisms
\[
\bigexact{\rule[-1.5ex]{0pt}{5ex}\cdots}%
%{\mbox{\raisebox{-0.5ex}[0pt][0pt]{$\scriptstyle\sim$}}\;\;}%\P_{-1,C_L\setminus D}\;}%
{\sim\quad\;\;\;}
{\Gamma(C_L\setminus D,\F_{-1})}{\sim\quad\;\;\;}
{\Gamma(C_L\setminus D,\F_{ 0})}{\sim\quad\;\;\;}
{\Gamma(C_L\setminus D,\F_{ 1})}{\sim\quad\;\;\;}
{\cdots}
\]
since $\coker\P_i$ is only supported at $\infty$ for all $i\in\Z$. So we set $M:=\Gamma(C_L\setminus D,\F_{0})$ and $\t:=(\P_0|_{C_L\setminus D})^{-1}\circ\t_0|_{C_L\setminus D}$, and we define $\ulM^{(D)}(\FF):=(M,\t)$. Obviously, $\ulM^{(D)}(\FF)$ is a $\t$-module on $\TA$ and $\ulM^{(\infty)}(\FF)$ is the pure Anderson motive $\ulM(\FF)$ studied above.

If $\infty\notin D$ fix $z$ as above. Set $M_i:=\Gamma(C_L\setminus D,\F_i)$ and define
\begin{equation} \label{EQ.Phi}
\ulM^{(D)}(\FF) \,:=\, M_0\oplus\dots\oplus M_{l-1}\qquad\text{with} 
\qquad 
\t \,:=\, \left(
\raisebox{5ex}{$
\xymatrix @=0pc {
0 \ar@{.}[rrr] \ar@{.}[ddddrrrr] & & & 0 & **{!L !<0.8pc,0pc> =<1.5pc,1.5pc>}\objectbox{\TP^{-1}\circ z^k\t_{l-1}} \\
\t_0 \ar@{.}[dddrrr] & & & & 0 \ar@{.}[ddd]\\
0 \ar@{.}[dd] \ar@{.}[ddrr] \\
\\
0 \ar@{.}[rr] & & 0 & \t_{l-2} & 0
}$}
\qquad\qquad\right)
\end{equation}
and $\TP:=\P_{l-1}\circ\dots\circ\P_0$. Clearly $\t$ depends on the choice of $k$, $l$, and $z$. Again $\ulM^{(D)}(\FF)$ is a $\t$-module on $\TA:=\Gamma(C\setminus D,\O_C)$. Notice that $\coker\t$ is supported on $\Graph(c)\cap(C_L\setminus D)$.

\begin{Definition}\label{Def1.17}
We call\/ $\ulM^{(D)}(\FF)$ the \emph{$\t$-module on $\TA$ associated to $\FF$}. We abbreviate $\ulM^{(\{\infty\})}(\FF)$ to $\ulM(\FF)$.
\end{Definition}

If $\infty\notin D$ the $\t$-module $\ulM^{(D)}(\FF)$ is equipped with the endomorphisms
\begin{equation} \label{EQ.Pi+Lambda}
\P \,:=\, \left(
\raisebox{5.5ex}{$
\xymatrix @=0pc {
0 \ar@{.}[rrr] \ar@{.}[ddddrrrr] & & & 0 & **{!L !<0.8pc,0pc> =<1.5pc,1.5pc>}\objectbox{\TP^{-1}\circ z^k\P_{l-1}} \\
\P_0 \ar@{.}[dddrrr] & & & & 0 \ar@{.}[ddd]\\
0 \ar@{.}[dd] \ar@{.}[ddrr] \\
\\
0 \ar@{.}[rr] & & 0 & \P_{l-2} & 0
}$}
\qquad\qquad\right) , \enspace
\Lambda(\lambda) \,:=\, \left(
\raisebox{5.5ex}{$
\xymatrix @=0pc {
\!\!\lambda \cdot\id_{M_0}\!\!  \\
& \!\!\lambda^q \cdot\id_{M_1}\!\! \ar@{.}[ddrr] \\ 
\\
& & & \!\lambda^{q^{l-1}}\cdot\id_{M_{l-1}}\!\!
}$}
\right)
\end{equation}
for all $\lambda\in \Ff_{q^l}\cap L$. Actually $\Lambda(\lambda)$ is even an automorphism and the same holds for $\P$ if $z$ has no zeroes on $C\setminus(D\cup\{\infty\})$. They satisfy the relations $\P^l=z^k$ and $\P\circ\Lambda(\lambda^q)=\Lambda(\lambda)\circ\P$.

\begin{Lemma}\label{Lemma2.9a}
Assume that $\chr\notin D$ or that $\chr\ne\infty\in D$. If $\FF$ and $\FF'$ are abelian $\tau$-sheaves of different weights, then $\Hom(\ulM^{(D)}(\FF),\ulM^{(D)}(\FF'))=\{0\}$ (for any choice of $k$, $l$, $k'$, $l'$ and $z$ if $\infty\notin D$).
\end{Lemma}

Before proving the lemma we note a direct consequence of its interaction with Theorem~\ref{PropX.1}.

\begin{Corollary}\label{Cor2.9b}
If $\ulM$ and $\ulM'$ are pure Anderson motives of different weights, then $\Hom(\ulM,\ulM')=\{0\}$.
\qed
\end{Corollary}

\noindent
{\it Remark.}
Again this follows alternatively from the Dieudonn\'e-Manin classification of the local isoshtuka at $\infty$ associated with $\ulM,\ulM'$; see Section~\ref{SectLS}.

\begin{proof}[Proof of lemma~\ref{Lemma2.9a}]
Let $f\in \Hom(\ulM^{(D)}(\FF),\ulM^{(D)}(\FF'))$. If $\infty\in D$ set $\CM:=\F_0$ and $\CM':=\F'_0$. If $\infty\notin D$ set $\CM:=\bigoplus_{i=0}^{l-1}\F_i$ and $\CM':=\bigoplus_{i=0}^{l'-1}\F'_i$. Then $\Gamma(C_L\setminus D,\CM)=\ulM^{(D)}(\FF)$ and likewise for $\FF'$. Thus $f$ extends to a homomorphism $f:\CM\to\CM'(m\cdot D)$ for a suitable  $m\in\N$. We abbreviate $\t^i:=\t\circ\s(\t)\circ\ldots\circ(\s)^{i-1}(\t)$. Let $z\in Q$ be a uniformizing parameter at $\infty$. If $\infty\in D$ and $\chr\ne\infty$ then
\[
z^k\t^l:(\s)^l\CM_\infty \isoto\CM_\infty\qquad\text{and}\quad z^{k'}(\t')^{l'}:(\s)^{l'}\CM'_\infty \isoto\CM'_\infty
\]
are isomorphisms on the stalks at $\infty$. So for any $n\in\N$ we have for the stalk of $f$ at $\infty$
\[
f_\infty\;=\;\bigl(z^{k'}(\t')^{l'}\bigr)^{nl}\circ(\s)^{nll'}(z^{n(kl'-k'l)}f_\infty)\circ(z^k\t^l)^{-nl'}:\es\CM_\infty\to\CM'(m\cdot D)_\infty\,.
\]
In particular if $\frac{k}{l}>\frac{k'}{l'}$ (and similarly for $\frac{k}{l}<\frac{k'}{l'}$), $f_\infty\equiv0\mod z^{n(kl'-k'l)}$ for all $n$, whence $f_\infty=0$. Thus $f=0$ since $\CM$ is locally free.

If $\infty\notin D$ and $\chr\notin D$ then with the $\tau$ from (\ref{EQ.Phi}) the homomorphisms on the stalks at every point $v\in D$
\[
z^{-k}\t^l:(\s)^l\CM_v \isoto\CM_v\qquad\text{and}\quad z^{-k'}(\t')^{l'}:(\s)^{l'}\CM'_v \isoto\CM'_v
\]
are isomorphisms since $\chr\notin D$. So again for any $n\in \N$ the stalk $f_v$ satisfies
\[
f_v\;=\;\bigl(z^{-k'}(\t')^{l'}\bigr)^{nl}\circ(\s)^{nll'}(z^{-n(kl'-k'l)}f_v)\circ(z^{-k}\t^l)^{-nl'}:\es\CM_v\to\CM'(m\cdot D)_v\,.
\]
There exists a pole $v\in D$ of $z$. Then for $\frac{k}{l}>\frac{k'}{l'}$ (and similarly for $\frac{k}{l}<\frac{k'}{l'}$), $f_v=0$, whence $f=0$ as desired.
\end{proof}

\begin{Example}\label{Ex1.18b}
We give an example showing that the assertion of the lemma is false in case $\chr=\infty\in D$. Let $C=\PP^1_\Fq$, $\F_i=\O_{C_L}(i\cdot\infty)$, $\F'_i=\O_{C_L}(2i\cdot\infty)$ and let $\P_i$ and $\t_i$ be the natural inclusions $\F_i\subset\F_{i+1}$ and $\s\F_i\subset\F_{i+1}$ and likewise for $\F'_i$. Then $\FF=(\F_i,\P_i,\t_i)$ is an abelian $\tau$-sheaf of weight $1$ and $\FF'=(\F'_i,\P'_i,\t'_i)$ is an abelian $\tau$-sheaf of weight $2$ both of characteristic $\infty$. Clearly $\ulM^{(\infty)}(\FF)=(A_L,\id)=\ulM^{(\infty)}(\FF')$ contradicting the assertion of the lemma.
\end{Example}

% -----------------------------------------------------------------------------

\bigskip

\subsection{Kernel sheaf and image sheaf}

In this section we show that the kernel and the image of a morphism of pure Anderson motives are themselves pure Anderson motives and likewise for abelian $\tau$-sheaves provided that the characteristic point $\chr=c(\Spec L)$ is different from $\infty$. 

\begin{Proposition}\label{Prop1.9a}
Let $\ulM$ and $\ulM'$ be pure Anderson motives and let $f\in\Hom(\ulM,\ulM')$. Then
\[
\ker f \;:=\;(\ker f,\tau|_{\s\ker f})\qquad\text{and}\qquad \im f\;:=\;(\im f,\t'|_{\s\im f})
\]
are again pure Anderson motives with $\weight(\ker f)=\weight(\im f)=\weight(\ulM)$.
\end{Proposition}

\begin{proof}
Let $\CM,\CM',k,l,k',l',z$ be as in definition~\ref{Def1'.1} and the subsequent remark. After replacing $\CM'_\infty$ by $z^{-n}\CM'_\infty$ for an integer $n$ we may assume that $f$ extends to a morphism $\CM\to\CM'$. Since all local rings of $C_L$ are discrete valuation rings the subsheaves $\wt\CM:=\ker f$ and $\wh\CM:=\im f$ are themselves locally free by the elementary divisor theorem. Set $\wt M:=\Gamma(C_L\setminus\{\infty\},\wt\CM)$ and $\wh M:=\Gamma(C_L\setminus\{\infty\},\wh\CM)$. Moreover the restrictions $\tilde\t:=\t|_{\s\wt M}$ and $\hat\t:=\t'|_{\s\wh M}$ are clearly injective. If $f\ne0$ then $\weight(\ulM)=\weight(\ulM')$ by corollary~\ref{Cor2.9b}. Let $\tilde l=\hat l$ be the least common multiple of $l$ and $l'$ and let $\tilde k=\hat k=\weight(\ulM)\cdot\tilde l$. Then
\[
z^{\tilde k}\t^{\tilde l}:(\s)^{\tilde l}\CM_\infty\isoto\CM_\infty \qquad\text{and}\qquad
z^{\tilde k}(\t')^{\tilde l}:(\s)^{\tilde l}\CM'_\infty\isoto\CM'_\infty
\]
are isomorphisms. Thus also
\[
z^{\tilde k}\tilde\t^{\tilde l}:(\s)^{\tilde l}\wt\CM_\infty\isoto\wt\CM_\infty \qquad\text{and}\qquad
z^{\hat k}\hat\t^{\hat l}:(\s)^{\hat l}\wh\CM_\infty\isoto\wh\CM_\infty
\]
are isomorphisms. Since $\tau$ and $\tau'$ are isomorphism outside $\Graph(c)$ the same is true for $\tau|_{\s\ker f}$ and $\tau'|_{\s\im f}$. So the cokernels of the latter two are supported on $\Graph(c)$. This proves the proposition by Remark~\ref{Rem2.2}.
\end{proof}

\begin{Proposition}\label{KERISABELIAN}\label{IMISABELIAN}%
Let $\FF$ and $\FF'$ be abelian $\tau$-sheaves of characteristic different from $\infty$ and let $f\in\Hom(\FF,\FF')$. Then the \emph{kernel $\tau$-sheaf} and the \emph{image $\tau$-sheaf}
\[
\begin{array}{r@{\;}c@{\;}l}
\ker f &:=& (\ker f_i, \P_i|_{\ker f_i}, \t_i|_{\s\ker f_i}) \\[1ex]
\im f &:=& (\im f_i, \P'_i|_{\im f_i}, \t'_i|_{\s\im f_i})
\end{array}
\]
are abelian $\tau$-sheaves of the same characteristic as $\FF$ and $\FF'$.
\end{Proposition}

\begin{proof}\def\kerf{\ker f}\def\imf{\im f}%
We will conduct the proof for $\ker f$ and $\im f$ simultaneously. If $f=0$, then $\ker f=\FF$ and $\im f=\ZZ$, and we are done. Otherwise, we have a non-zero morphism between $\FF$ and $\FF'$, and by proposition \ref{PROP.1} both abelian $\tau$-sheaves $\FF$ and $\FF'$ have the same weight. We choose an integer $l$ that satisfies condition 2 of \ref{Def1.1} for both $\FF$ and $\FF'$ and we set $k=l\cdot\weight(\FF)$.

\def\kerf{\TF}\def\imf{\HF}%
Let $i\in\Z$. Since all local rings of $C_L$ are principal ideal domains the elementary divisor theorem yields that $\TF_i:=\ker f_i\subset\F_i$ and $\HF_i:=\im f_i\subset\F'_i$ are locally free coherent sheaves. 
The induced morphisms $\TP_i := \P_i|_{\kerf_i}$ and $\Tt_i := \t_i|_{\s\kerf_i}$ map injectively into $\kerf_{i+1}$ since $\s\ker f_i=\ker \s f_i$ due to the flatness of $\sigma$.
Similarly, we get this for $\HP_i := \P'_i|_{\imf_i}$ and $\Ht_i := \t'_i|_{\s\imf_i}$.
To examine $\coker\TP_i$ and $\coker\HP_i$ consider the diagram with exact
rows and columns in which the last column is exact by the 9-lemma
\[
\xymatrix{
& 0 \ar[d] & 0 \ar[d] & 0 \ar[d] & \\
0 \ar[r]& \kerf_i \ar[r]^{\TP_i}\ar[d]      & \kerf_{i+1} \ar[r]\ar[d]          & \coker\TP_i \ar[r]\ar[d]& 0 \\
0 \ar[r]& \F_i    \ar[r]^{\P_i} \ar[d]^{f_i}& \F_{i+1}    \ar[r]\ar[d]^{f_{i+1}}& \coker\P_i  \ar[r]\ar[d]& 0 \\
0 \ar[r]& \imf_i  \ar[r]^{\HP_i}\ar[d]      & \imf_{i+1}  \ar[r]\ar[d]          & \coker\HP_i \ar[r]\ar[d]& 0 \\
& 0 & 0 & 0 & \\
}
\]
Thus $\coker\TP_i$ and $\coker\HP_i$ are torsion sheaves like $\coker\P_i$, and we conclude that
the ranks $\Tr:= \rank\kerf_i$ and $\Hr:=\rank\imf_i$ are independent of $i$.

To show that $\underline\kerf$ and $\underline\imf$ are abelian $\tau$-sheaves let $\P$ and $\P'$ be the identifying morphisms $\P_{i+l-1}\circ\cdots\circ\P_i: \F_i \isoto \F_{i+l}(-k\cdot\infty)$ and $\P'_{i+l-1}\circ\cdots\circ\P'_i: \F'_i \isoto\F'_{i+l}(-k\cdot\infty)$, respectively. Since $\P$ and $\P'$ are isomorphisms we obtain the same for $\TP_{i+l-1}\circ\ldots\circ\TP_i$ and $\HP_{i+l-1}\circ\ldots\circ\HP_i$, whence the periodicity condition 2.

To establish conditions 3 and 4 we need that the characteristic is different from $\infty$. Let $c: \Spec L\rightarrow C':=C\setminus\{\infty\}$ and let
\[
\xymatrix{ & & &
\makebox[0pt][r]{$M \::=\:$ }%
\kerf_0|_{\CLa}\: \ar[r]^<<{\sim}^{\:\:\TP_0} &
\:\kerf_1|_{\CLa}\: \ar[r]^<<{\sim}^{\TP_1} & 
\:\cdots\makebox[2.4em][l]{}
}
\]
Set $\Tt := {\TP_0}^{-1}\!\circ\Tt_0:\; \s M\rightarrow M$ and set $\Td := \dim_L\coker\Tt$. Similar to the diagram chase for the $\coker\TP_i$, we get an injective morphism $\coker\Tt_i\rightarrow\coker\t_i$. Hence the support of $\coker\Tt_i$ lies outside $\infty$, and we have $\coker\Tt_i \:=\: \coker\Tt_i|_{\CLa} \:\cong\: \coker\Tt$ for all $i\in\Z$.
Now the exact sequences
\[
\begin{CD} 
{0} @>>> {\kerf_i}   @>{\TP_i}>> {\kerf_{i+1}} @>>> {\coker\TP_i} @>>> {0} \\[1ex]
{0} @>>> {\s\kerf_i} @>{\Tt_i}>> {\kerf_{i+1}} @>>> {\coker\Tt_i} @>>> {0} 
\end{CD}
\]
yield 
\[
\dim_L\coker\TP_i\:=\:\deg\kerf_{i+1}- \deg\kerf_i\:=\:\deg\kerf_{i+1}-\deg\s\kerf_i \:=\:\dim_L\coker \Tt_i\:=\:\dim_L\coker\Tt \:=\: \Td
\]
for all $i\in\Z$.
Clearly, $\coker\Tt_i$ is supported on the graph of $c$ due to its injection into $\coker\t_i$. Again, this argument adapts to $\coker\HP_i$ and $\coker\Ht_i$, as well.
\end{proof}

\begin{Corollary}
Let\/ $\FF$ and $\FF'$ be abelian $\tau$-sheaves of characteristic different from $\infty$ and let\/ $f \in \Hom(\FF,\FF')$ be a morphism. 
\begin{suchthat}
\item $f$ is injective if and only if\/ $\ker f=\ZZ$\,.
\item $f$ is surjective if and only if\/ $\im f=\FF'$\,. \qed
\end{suchthat}
\end{Corollary}

\begin{Example} \label{ExTRichter}
As was pointed out to us by T.\ Richter the assumption $\chr\ne\infty$ cannot be dropped. For instance consider the abelian $\tau$-sheaf on $C_L=\PP^1_L$ with $\F_i=\O_{\PP^1_L}\bigl(\bigl\lceil\frac{i-1}{2}\bigr\rceil\bigr)\oplus\O_{\PP^1_L}\bigl(\bigl\lceil\frac{i}{2}\bigr\rceil\bigr)$,
where $\bigl\lceil\frac{i}{2}\bigr\rceil$ denotes the smallest integer $\ge\frac{i}{2}$. Let $\P_i=\left(\begin{array}{cc} 1 & 0 \\ 0 & 1 \end{array}\right)$ and $\t_i=\left(\begin{array}{cc} 0 & 1 \\ 1 & 0 \end{array}\right)$. Then $\FF=(\F_i,\P_i,\t_i)$ is an abelian $\tau$-sheaf with $r=l=2,d=k=1$, and characteristic $\infty$. We rewrite everything in terms of the bases $\biggl(\!\!\!\begin{array}{c}z^{-\lceil\frac{i-1}{2}\rceil}\\ 0 \end{array}\!\!\!\biggr),\biggl(\!\!\begin{array}{c}0\\ z^{-\lceil\frac{i}{2}\rceil}\end{array}\!\!\biggr)$ of $\F_i|_{\PP^1_L\setminus\{0\}}$, where $\PP^1_L\setminus\{0\}=\Spec L[z]$. With respect to these bases $\P_i$ and $\t_i$ are described by the matrices
\[
\P_i=\left(\begin{array}{cc} 1 & 0 \\ 0 & z \end{array}\right)\text{ for }2|i\,,\quad\P_i=\left(\begin{array}{cc} z & 0 \\ 0 & 1 \end{array}\right)\text{ for }2\nmid i\,,\qquad \text{and} \qquad \t_i=\left(\begin{array}{cc} 0 & 1 \\ z & 0 \end{array}\right)\text{ for all }i\,.
\]
There is an endomorphism $f$ of $\FF$ given by $f_i=\left(\begin{array}{cc} z & z \\ z & z \end{array}\right)\text{ for }2|i\,$ and $f_i=\left(\begin{array}{cc} z & 1 \\ z^2 & z \end{array}\right)\text{ for }2\nmid i$. We compute
\[
\begin{array}{rcll}
\TS\ker f_i&=&\TS{-1\choose 1}\cdot\O_{\PP^1_L}(\frac{i}{2}\cdot\infty) & \quad \text{ if }2|i \text{ and}\\[2mm]
\TS\ker f_i&=&\TS{-1\choose z}\cdot\O_{\PP^1_L}(\frac{i-1}{2}\cdot\infty) & \quad \text{ if }2\nmid i\,.
\end{array}
\]
Therefore $\P_i|_{\ker f_i}$ is an isomorphism if $2|i$ and has cokernel of $L$-dimension $1$ for $2\nmid i$. Thus $\ker f$ is not an abelian $\tau$-sheaf.
\end{Example}

% -----------------------------------------------------------------------------

\bigskip

\subsection{Isogenies between Abelian $\tau$-Sheaves and Pure Anderson Motives}

In the theory of abelian varieties the concept of \emph{isogenies} is central, defining an equivalence relation which allows a classification of abelian varieties into isogeny classes that are larger than isomorphism classes. In the following, we adapt the idea of isogenies to abelian $\tau$-sheaves. They are defined by the following conditions.

\begin{Proposition}\label{PROP.1.42A}% 
  Let $f:\FF\to\FF'$ be a morphism between two abelian $\tau$-sheaves $\FF = (\F_i,\P_i,\t_i)$ and\/ $\FF' = (\F'_i,\P'_i,\t'_i)$. Then the following assertions are equivalent:
\begin{suchthat}
\item 
$f$ is injective and the support of all\/ $\coker f_i$ is contained in $D\times\Spec L$ for a finite closed subscheme $D\subset C$,
\item 
$f$ is injective and $\FF$ and $\FF'$ have the same rank and dimension,
\item 
$\FF$ and $\FF'$ have the same weight and the fiber $f_{i,\eta}$ at the generic point $\eta$ of $C_L$ is an isomorphism for some (any) $i\in\Z$.
\end{suchthat}
\end{Proposition}

\begin{proof}
$1\Rightarrow 3$ follows from \ref{PROP.1} and the fact that $\P_{i,\eta}$ and $\P'_{i,\eta}$ are isomorphisms for all $i$.
Since $3\Rightarrow 2$ is evident it remains to establish $2\Rightarrow 1$.

We will first reduce to the case $A=\BF_q[t]$. By the theorem of Riemann-Roch there exists a rational function $t\in Q$ on $C$ with poles only at $\infty$ and whose zeroes do not meet the characteristic point, nor the support of the $\coker f_i$. This function defines an inclusion of function fields $\Fq(t)\subset Q$ and thus a finite flat morphism between the respective curves $\varphi:\, C\rightarrow \PP^1_\Fq$ with $\varphi^{-1}(\infty_{\PP^1})=\{\infty\}$. The direct images $\GG:=\varphi_\ast\FF$ and $\GG':=\varphi_\ast\FF'$ under $\varphi$ are abelian $\tau$-sheaves on $\PP^1_\Fq$ of rank $r\cdot\deg\varphi$, dimension $d$, and characteristic $\varphi\circ c$ by \cite[Proposition~1.6]{HH}. We define $\TA:=\Gamma(\PP^1_\Fq\setminus\{0\},\O_{\PP^1})$ such that $\TA=\Fq\II{z}$ for some $z\in\TA$ with a simple pole at $0$ and a simple zero at $\infty$. We choose an integer $l$ that satisfies condition 2 of \ref{Def1.1} for both $\GG$ and $\GG'$. Consider $\ulM^{(0)}(\GG)=(M,\t)$ and $\ulM^{(0)}(\GG')=(M',\t')$; see Definition~\ref{Def1.17}. Set $s:=lr\deg\varphi=\rk M$ and $e:=ld=s\cdot\weight(\GG)$.

Now choose $\TA_L$-bases of $M$ and $M'$. This is possible since $\TA_L$ is a principal ideal domain and that was the reason why we constructed $\varphi$. With respect to these bases, the endomorphisms $\t$ and $\t'$ and the $\TA$-morphism $g=\ulM^{(0)}(\varphi_\ast f):\, M\rightarrow M'$ which is induced by $f$ can be described by quadratic matrices $T$, $T'$ and $H$, and we have the formula $T'\s\!H = HT$.

Let $\zeta:=c^\ast(z)\in L$. By the elementary divisor theorem we find matrices $U,V\in GL_s(\AxL)$ with 
\[
UTV \,=\, \left(
\begin{array}{ccc} 
\makebox[2.5em][l]{$(z-\zeta)^{n_1}$} & & \\
& \ddots & \\
& & \makebox[2.5em][r]{$(z-\zeta)^{n_s}$}
\end{array}
\right)
\]
for some integers $n_1\le\dots\le n_s$. Thus $\coker\t\cong \bigoplus_{i=1}^s \AxL\,/\,(z-\zeta)^{n_i}$ and  $e=\sum_{i=1}^s n_i$. Since
\[
\det T\cdot\det UV = \det \,UTV = (z-\zeta)^e 
\]
we calculate $\det T=b\cdot(z-\zeta)^e$ for some $b\in(\AxL)^{\times}={L\II{z}}^{\!\times}=L^{\!\times}$. Analogously, we have $\det T'=b'\cdot(z-\zeta)^e$ for some $b'\in L^{\!\times}$ as $\GG$ and $\GG'$ have the same dimension $d$. Since $\det H\ne0$ due to the injectivity of $f$ we conclude
\[
\det T'\cdot\det\s\!H \,=\, \det H\cdot\det T 
\quad\Rightarrow\quad
\frac{\displaystyle \det\s\!H}{\displaystyle \det H} \,=\, \frac{\displaystyle b}{\displaystyle b'\!\!} \;\in L^{\!\times}\,.
\]
In an algebraic closure $L^\alg$ of $L$ there exists a $\lambda$ with $\lambda^{q-1}=\frac{b'\!\!}{b}$. So we have 
\[
a:=\lambda\cdot\det H=\s(\lambda\cdot\det H)\in L^\alg[z]
\]
and, due to the $\sigma$-invariance, even $a\in\Fq\II{z}=\TA$ (and hence $\lambda\in L$).
Again using the elementary divisor theorem one sees that $a$ annihilates $\coker g$. Now our proof is complete as the support of $\coker f_i$ is contained in the divisor of zeroes $(\varphi^\ast(a))_0\times\Spec L$ on $C_L$ for $0\le i<l$ by construction (for this purpose we used $g=\ulM^{(0)}(\varphi_\ast f)$ which captures all these $f_i$) and for the remaining $i\in\Z$ by periodicity.
\end{proof}

\begin{Definition}[Isogeny]
\begin{suchthat}
\item 
A morphism $f:\FF\to \FF'$ satisfying the equivalent conditions of proposition~\ref{PROP.1.42A} is called an \emph{isogeny}.
We denote the set of isogenies between $\FF$ and $\FF'$ by $\Isog(\FF,\FF')$.
\item 
An isogeny $f:\FF\to\FF'$ is called \emph{separable} (respectively \emph{purely inseparable}) if for all $i$ the induced morphism $\t_i:\s\coker f_i\to\coker f_{i+1}$ is an isomorphism (respectively is \emph{nilpotent}, that is, $\t_{i}\circ\s\t_{i-1}\circ\ldots\circ(\s)^n\t_{i-n}=0$ for some $n$).
\end{suchthat}
\end{Definition}

\medskip

Proposition~\ref{PROP.1.42A} has  important consequences for pure Anderson motives.

\begin{Corollary}\label{Cor1.26a}
Let $f:\FF\to\FF'$ be a morphism between abelian $\tau$-sheaves of characteristic different from $\infty$. Then $f$ is an isogeny if and only if $\ulM(f):\ulM(\FF)\to\ulM(\FF')$ is an isogeny between the associated pure Anderson motives.\qed
\end{Corollary}

\begin{Corollary}\label{Cor1.11b}
Let $f:\ulM\to\ulM'$ be an isogeny between pure Anderson motives (Definition~\ref{Def1'.2}). Then
\begin{suchthat}
\item
there exists an element $a\in A$ which annihilates $\coker f$,
\item 
there exists a dual isogeny $\dual{f}:\ulM'\to\ulM$ such that $f\circ\dual{f}=a\cdot\id_{\ulM'}$ and $\dual{f}\circ f=\id_\ulM$.
\end{suchthat}
\end{Corollary}

\begin{proof}
1 follows from Corollary~\ref{Cor2.9b}, Theorem~\ref{PropX.1}, and Proposition~\ref{PROP.1.42A} by noting that $D$ is contained in the zero locus of a suitable $a\in A$ by the Riemann-Roch theorem.

For 2  consider the diagram 
\[
\xymatrix{
0 \ar[r] &
M \ar[r]^{f}\ar[d]_{a} &
M' \ar[r]\ar[d]_{a} \ar@{-->}[dl]_{\dual{f}} &
\coker f \ar[d]_{a}^{\;(=0)} \ar[r] & 
0
\\
0 \ar[r] &
M \ar[r]^{f} &
M' \ar[r] &
\coker f \ar[r] & 
0\,.
}
\]
By diagram chase, we get a morphism $\dual{f}: M'\rightarrow M$ which is \emph{dual} to $f$ in the sense that $\dual{f}\circ f = a$ and $f\circ\dual{f} = a$. 
\end{proof}

\begin{Remark}\label{Rem1.26'}
The dual isogeny $\dual{f}$ clearly depends on $a$ and there rarely is a canonical choice for $a$. If $C=\PP^1$ and $A=\Fq[t]$ we obtain from the elementary divisor theorem a unique minimal monic element $a\in A$ (which still depends on the choice of the isomorphism $A\cong\Fq[t]$, though) that annihilates $\coker f$. Also if $f\in\End(\ulM)$ is an isogeny of a semisimple pure Anderson motive over a finite field we will exhibit in Theorem~\ref{Prop3.4.1} below a canonical $a$. But in general there is no canonical choice.

Nevertheless, since $A$ is a Dedekind domain, a power of the ideal annihilating $\coker f$ will be principal and one may take $a$ as a generator. This has the advantage that the support of $\coker f$ equals $V(a)\subset\Spec A$. In particular if the characteristic point $\chr$ is not contained in the support of $\coker f$ and in $V(a)$, also $\dual{f}$ will be separable. On the other hand, if $f\in \End(\ulM)$ then $f$ is integral over $A$, since $\End(\ulM)$ is a finite $A$-module by Proposition~\ref{PropT.1} below. Then $f$ generates a finite commutative $A$-algebra $A[f]$. Our discussion of the choice of $a$ shows that the set $V(f)\subset\Spec A[f]$ of zeroes of $f$ lies above $\supp(\coker f)\subset\Spec A$.
\end{Remark}

\begin{Lemma}\label{ISOGENYCOMPOSE}
Let $f\in\Hom(\FF,\FF')$ and $f'\in\Hom(\FF',\FF'')$ be morphisms between abelian $\tau$-sheaves and let $D$ be a divisor on $C$. 
\begin{suchthat}
\item If two of $f$, $f'$, and $f'\circ f$ are isogenies then so is the third.
\item If $f$ is an isogeny then also $f\otimesidOCL{D}\in\Isog(\FF(D),\FF'(D))$ is an isogeny.
\item If $D$ is effective then the natural inclusion $\FF\subset\FF(D)$ is an isogeny.
\end{suchthat}
\end{Lemma}

\begin{proof}
1 is obvious.

\smallskip
\noindent 
2. Clearly the tensored morphisms $f_i\otimes 1:\F_i(D)\rightarrow\F'_i(D)$ remain injective and the support of $\coker f_i\otimes 1$ equals the support of $\coker f_i$.

\smallskip
\noindent
3. The inclusion $\FF\subset\FF(D)$ is a morphism because the divisor $D$ is $\sigma$-invariant.
\end{proof}

\medskip

In the following we want to define the \emph{degree} of an isogeny which should be an ideal of $A$ since in the function field case we have substituted $A$ for $\Z$.
Let $f:\ulM\to\ulM'$ be an isogeny between pure Anderson motives. Then the $A_L$-module $\coker f$ is a finite $L$-vector space equipped with a morphism of $A_L$-modules $\t':\s\coker f\to\coker f$. Since $\coker f$ is annihilated by an element of $A$ it decomposes by the Chinese remainder theorem
\[
(\coker f,\t')\;=\;\bigoplus_{v\in\supp(\coker f)}(\coker f,\t')\otimes_A A_v\;=:\;\bigoplus_{v\in\supp(\coker f)}\ulK_v\,.
\]
If $v\ne\chr$ the morphism $\t'$ on $\ulK_v$ is an isomorphism and so Lang's theorem implies that
\[
(\ulK_v\otimes_L L^\sep)^\t\otimes_{\Fq} L^\sep \isoto \ulK_v\otimes_L L^\sep
\]
is an isomorphism; see for instance \cite[Lemma 1.8.2]{Anderson}. In particular
\[
[\BF_v:\Fq]\cdot\dim_{\BF_v}(\ulK_v\otimes_L L^\sep)^\t\;=\;\dim_\Fq(\ulK_v\otimes_L L^\sep)^\t\;=\;\dim_{L^\sep}(\ulK_v\otimes_L L^\sep)\;=\;\dim_L \ulK_v\,.
\]

\smallskip

On the other hand if the characteristic is finite and $v=\chr$, the characteristic morphism $c^\ast:A\to L$ yields $\BF_\chr\subset L$ and determines the distinguished prime ideal 
\[
\Fa_0\;:=\;(a\otimes 1-a\otimes c^\ast(a):a\in\BF_\chr)\subset A_{\chr,L}\,. 
\]
If we set $n:=[\BF_\chr:\Fq]$ and $\Fa_i:=(\s)^i\Fa_0=(a\otimes 1-a\otimes c^\ast(a)^{q^i}:a\in\BF_\chr)$, then we can decompose $A_{\chr,L}=\bigoplus_{i\in\Z/n\Z}A_{\chr,L}/\Fa_i$ and $\t$ is an isomorphism
\[
\s(\ulK_\chr/\Fa_{i-1}\ulK_\chr)\isoto \ulK_\chr/\Fa_i
\]
for $i\ne0$ since $\t$ is an isomorphism on $\ulM$ and $\ulM'$ outside the graph of $c^\ast$. 
% FIX 2: THIS SENTENCE HAS TO BE MODIFIED IF THE ``DEGREE'' MOVES BACK INTO CHAPTER 3
 (This argument will be used again in Proposition~\ref{PropLS4}.)
% FIX 
 In particular 
\[
[\BF_\chr:\Fq]\cdot\dim_L(\ulK_\chr/\Fa_0\ulK_\chr)\;=\;\dim_L\ulK_\chr\,.
\]

\begin{Definition}\label{Def1.7.6}
We assign to the isogeny $f$ the ideal
\[
\deg(f)\;:=\;\prod_{v\in\supp(\coker f)}v^{(\dim_L\ulK_v)/[\BF_v:\Fq]}\;=\;\chr^{\dim_L(\ulK_\chr/\Fa_0\ulK_\chr)}\cdot\prod_{v\ne\chr}v^{\dim_{\BF_v}(\ulK_v\otimes_L L^\sep)^\t}
\]
of $A$ and call it the \emph{degree of $f$}. We call $\chr^{\dim_L(\ulK_\chr/\Fa_0\ulK_\chr)}$ the \emph{inseparability degree of $f$} and $\prod_{v\ne\chr}v^{\dim_{\BF_v}(\ulK_v\otimes_L L^\sep)^\t}$ the \emph{separability degree of $f$}.
\end{Definition}

\noindent
{\it Remark.} The separability degree of $f$ is the Euler-Poincar\'e characteristic $EP\bigl(\bigoplus_{v\ne\chr}\ulK_v\otimes_L L^\sep\bigr)^\t$; see Gekeler~\cite[3.9]{Gekeler} or Pink-Traulsen~\cite[4.6]{PT}. Recall that the Euler-Poincar\'e characteristic of a finite torsion $A$-module is the ideal of $A$ defined by requiring that $EP$ is multiplicative in short exact sequences, and that $EP(A/v):=v$ for any maximal ideal $v$ of $A$.

\begin{Lemma}\label{Lemma1.7.7}
\begin{suchthat}
\item 
If $f:\ulM\to\ulM'$ and $g:\ulM'\to\ulM''$ are isogenies then $\deg(gf)=\deg(f)\cdot\deg(g)$.
\item 
$\dim_{\Fq}A/\deg(f)\;=\;\dim_L\coker f$.
\end{suchthat}
\end{Lemma}

\begin{proof}
1 is immediate from the short exact sequence
\[
\xymatrix{ 0 \ar[r] & \coker f \ar[r]^g & \coker(gf) \ar[r] & \coker g \ar[r] & 0 }
\]
and 2 is obvious.
\end{proof}

\begin{Proposition}\label{Prop3.28a}
The ideal $\deg(f)$ annihilates $\coker f$.
\end{Proposition}

\begin{proof}
If $v=\chr$ and $a$ is a uniformizer at $\chr$, then multiplication with $a$ is nilpotent on the $L$-vector space $\ulK_\chr/\Fa_0\ulK_\chr$. In particular $a^{\dim_L(\ulK_\chr/\Fa_0\ulK_\chr)}$ annihilates $\ulK_\chr/\Fa_0\ulK_\chr$, and hence also $\ulK_\chr$.

If $v\ne\chr$ and $a$ is a uniformizer at $v$, we obtain analogously that $a^{\dim_{\BF_v}(\ulK_v\otimes_L L^\sep)^\t}$ annihilates the $\BF_v$-vector space $(\ulK_v\otimes_L L^\sep)^\t$ and therefore also the $L$-vector space $\ulK_v$.
\end{proof}

\begin{Proposition}\label{Prop1.7.8}
Let $f:\ulM\to \ulM'$ be an isogeny such that $\deg(f)=aA$ is principal (for example this is the case if $C=\PP^1$ and $A=\Fq[t]$). Then there is a uniquely determined dual isogeny $\dual{f}:\ulM'\to \ulM$ satisfying $f\circ\dual{f}=a\cdot\id_{\ulM'}$ and $\dual{f}\circ f=a\cdot\id_\ulM$.
\end{Proposition}

\begin{proof}
Since $\deg(f)$ annihilates $\coker f$ the proposition is immediate.
\end{proof}

In Theorem~\ref{Prop3.4.1} we will see that any isogeny $f\in\End(\ulM)$ of a semisimple pure Anderson motive over a finite field satisfies the assumption that $\deg(f)$ is principal.

% -----------------------------------------------------------------------------

\bigskip

\subsection{Quasi-morphisms and  quasi-isogenies}

We want to establish the existence of dual isogenies also for abelian $\tau$-sheaves.
But if we follow the proof of Corollary~\ref{Cor1.11b}, the problem is that multiplication with $a$ is not an endomorphism of an abelian $\tau$-sheaf, since it produces poles. We remedy this by defining \emph{quasi-morphisms} and \emph{quasi-isogenies} between $\FF$ and $\FF'$ which allow the maps to have finite sets of poles.

\begin{Definition}[Quasi-morphism and quasi-isogeny]
Let $\FF$ and $\FF'$ be abelian $\tau$-sheaves.
\begin{suchthat}
\item A \emph{quasi-morphism} $f$ between $\FF$ and $\FF'$ is a morphism $f\in\Hom(\FF,\FF'(D))$ for some effective divisor $D$ on $C$.
\item A \emph{quasi-isogeny} $f$ between $\FF$ and $\FF'$ is an isogeny $f\in\Isog(\FF,\FF'(D))$ for some effective divisor $D$ on $C$. 
\end{suchthat}
We call two quasi-morphisms $f_1\in\Hom(\FF,\FF'(D_1))$ and $f_2\in\Hom(\FF,\FF'(D_2))$ \emph{equivalent} $(\text{denoted }f_1\sim f_2)$, if the diagram
\[
\xymatrix@=0pt@C=3em{ 
& \FF'(D_1) \ar[rd] & \\
\FF \ar[ru]^{f_1}\ar[rd]_{f_2} & & \FF'\makebox[0pt][l]{$(D_1\!+\!D_2)$} \\
& \FF'(D_2) \ar[ru] &
}
\]
commutes where the two arrows on the right are the natural inclusions.
\end{Definition}

Clearly, the relation $\sim$ defines an equivalence relation on the set of quasi-morphisms between $\FF$ and $\FF'$ where the transitivity is seen from
\[
\xymatrix@=3ex@C=3em{ 
& \FF'(D_1) \ar[r]\ar@{-->}[rd] & \FF'\makebox[0pt][l]{$(D_1\!+\!D_2)$} &\quad\; \ar[rd] & \\
\FF \ar[ru]^{f_1}\ar[r]^{f_2}\ar[rd]_{f_3} & \FF'(D_2) \ar[ru]\ar[rd] & \FF'\makebox[0pt][l]{$(D_1\!+\!D_3)$} &\; \ar@{.>}[r] & \FF'\makebox[0pt][l]{$(D_1\!+\!D_2\!+\!D_3)$} \\
& \FF'(D_3) \ar[r]\ar@{-->}[ru] & \FF'\makebox[0pt][l]{$(D_2\!+\!D_3)$} &\quad\; \ar[ru] &
}
\]
by canceling the dotted arrow due to injectivity. Since the divisors of quasi-morphisms are not particularly interesting, we fade them out by forming equivalence classes of quasi-morphisms according to this equivalence relation.

\medskip

\begin{Definition}
Let $\FF$ and $\FF'$ be abelian $\tau$-sheaves. \nopagebreak
\begin{suchthat}
\nopagebreak
\item We set $\QHom(\FF,\FF')$ to be the set of quasi-morphisms between $\FF$ and $\FF'$ modulo $\sim$. 
\item The equivalence class of a quasi-morphism $f$ between $\FF$ and $\FF'$ modulo $\sim$ is denoted by $\II{f}$, and we call it a \emph{quasi-morphism} between $\FF$ and $\FF'$ as well.
\item We set $\QIsog(\FF,\FF')$ to be the subset of\/ $\QHom(\FF,\FF')$ whose elements $\II{f}$ are represented by isogenies $f$.
\end{suchthat}
\vspace{1mm}
We write $\QEnd(\FF):=\QHom(\FF,\FF)$ and $\QIsog(\FF):=\QIsog(\FF,\FF)$.
\end{Definition}

\noindent {\it Remark.}\INDENT
1. By Lemma~\ref{ISOGENYCOMPOSE}, it holds for $f_1\sim f_2$, that $f_1$ is a quasi-isogeny if and only if $f_2$ is a quasi-isogeny. This justifies our definition of $\QIsog(\FF,\FF')$.

2. Proposition \ref{PROP.1.42A} and Lemma~\ref{ISOGENYCOMPOSE} hold analogously for quasi-morphisms and quasi-isogenies, since every quasi-morphism $f\in\QHom(\FF,\FF')$ is a morphism $f\in\Hom(\FF,\FF'(D))$ for some effective divisor $D$ on $C$.

3. Every pair of quasi-morphisms $\II{f_1},\II{f_2}\in\QHom(\FF,\FF')$ can be represented by morphisms $f_1,f_2\in\Hom(\FF,\FF'(D))$ with the same divisor $D$ on $C$. 
In particular we can form the sum
\[
\II{f_1}+\II{f_2} := \II{f_1+f_2} \:\in\QHom(\FF,\FF')\,.
\]
Since poles are negligible, we can extend this structure to a $Q$-vector space by now being able to admit multiplication by elements of $Q$. Let $\II{f}\in\QHom(\FF,\FF')$ be represented by $f\in\Hom(\FF,\FF'(D))$ and let $a\in Q$. Then $a\cdot f\in\Hom(\FF,\FF'(D\!+\!(a)_\infty))$ where $(a)_\infty$ denotes the divisor of poles of $a$, and we define
\[
a\cdot\II{f} := \II{a\cdot f} \:\in\QHom(\FF,\FF')\,.
\]
Moreover, Quasi-morphisms can be composed. Let $\FF$, $\FF'$ and $\FF''$ be abelian $\tau$-sheaves and let $\II{f}\in\QHom(\FF,\FF')$ and $\II{f'}\in\QHom(\FF',\FF'')$ be quasi-morphisms between $\FF$, $\FF'$ and $\FF''$, which are represented by $f\in\Hom(\FF,\FF'(D))$ and $f'\in\Hom(\FF',\FF''(D'))$, respectively. In order to compose $f'$ and $f$, we have to raise $f'$ to be a morphism from $\FF'(D)$. We achieve this by simply tensoring with $\otimes_\OCL\OCL(D)$. Now $(f'\otimesidOCL{D})\circ f\in\Hom(\FF,\FF''(D\!+\!D'))$, and we can define the composition
\[
\II{f'}\circ\II{f} := \II{(f'\otimesidOCL{D})\circ f} \:\in\QHom(\FF,\FF'')\,.
\]
All these operations are well-defined which can easily be seen by diagram arguments similar to the one we presented for the transitivity of $\sim$. Altogether we obtain

\begin{Corollary}\label{QEND-Q-ALGEBRA}\label{COMPOSITIONISISOGENY}
Let $\FF$ and $\FF'$ be abelian $\tau$-sheaves. With the above given structure, we have
\begin{suchthat}
\item the composition of quasi-isogenies is again a quasi-isogeny,
\item $\QHom(\FF,\FF')$ is a $Q$-vector space,
\item $\QEnd(\FF)$ is a $Q$-algebra. \qed
\end{suchthat}
\end{Corollary}

\smallskip

\noindent {\it Remark.}
The $Q$-vector spaces $\QHom(\FF,\FF')$ and $\QEnd(\FF)$ are finite dimensional. We will prove this in Proposition~\ref{PropT.2} below.

\bigskip

As an abuse of notation, we will write $f\in\QHom(\FF,\FF')$ instead of $\II{f}$ to denote the quasi-morphism represented by $f\in\Hom(\FF,\FF'(D))$.

\begin{Remark} \label{MULTISISOGENY}
For every $a\in Q\mal$, the multiplication by $a$ is a quasi-isogeny on $\FF$. Since $a$ injects $\F_i$ into $\F_i((a)_\infty)$ and commutes with the $\P_i$ and the $\t_i$, it is a morphism of abelian $\tau$-sheaves. Additionally, its cokernels are supported on $(a)_0$, the divisor of zeroes of $a$.
\end{Remark}

Now we come back to the idea of defining a dual isogeny. As already mentioned, a global definition fails because the annihilating multiplication by $a$ is not a morphism between $\F_i$ and $\F'_i$. This problem will now be solved by using quasi-morphisms and quasi-isogenies.

Let $\FF$ and $\FF'$ be abelian $\tau$-sheaves and let $f\in\QIsog(\FF,\FF')$ be a quasi-isogeny represented by an isogeny $f:\FF\to\FF'(D)$ for an effective divisor $D$ on $C$. By the annihilating property of the support, we can find $a\in Q\mal$ with $a\cdot\coker f_i=0$ for all $i\in\Z$. Now consider the following diagram.
\[
\xymatrix@M=0.5em{
0 \ar[r] &
\F_i \ar[rr]^{f_i}\ar[d]_{a} &&
\F'_i(D) \ar[rr]\ar[d]_{a} \ar@{-->}[dll]_{\dual{f_i}} &&
\coker f_i \ar[d]_{a}^{\;(=0)} \ar[r] & 
0
\\
0 \ar[r] &
\F_i\makebox[0pt][l]{$((a)_\infty)$} & \ar[r]^{f_i\quad} &
\F'_i\makebox[19pt][l]{$(D+(a)_\infty)$} & \ar[r] &
\coker f_i \ar[r] & 
0\,.
}
\]
As in \ref{Cor1.11b}, we get a morphism $\dual{f_i}: \F'_i(D)\rightarrow\F_i\bigl((a)_\infty\bigr)$ satisfying $\dual{f_i}\circ f_i = a$ and $f_i\circ\dual{f_i} = a$. Collecting these $\dual{f_i}$ together, we obtain a \emph{dual} morphism of abelian $\tau$-sheaves $\dual{f}\in\Hom(\FF'(D),\FF((a)_\infty))$ which is a quasi-morphism between $\FF'$ and $\FF$. 

\begin{Proposition} \label{QISOG-GROUP}%
Let $\FF$ and $\FF'$ be abelian $\tau$-sheaves.
\begin{suchthat}
\item
If $f\in \QIsog(\FF,\FF')$ is a quasi-isogeny then the dual $\dual{f}\in\QHom(\FF',\FF)$ of\/ $f$ is a quasi-isogeny and $f^{-1}:=a^{-1}\cdot\dual{f}$ is the inverse of $f$ in $\QHom(\FF',\FF)$.
\item 
$\QIsog(\FF)$ is the group of units in the $Q$-algebra $\QEnd(\FF)$.
\end{suchthat}
\end{Proposition}

\begin{proof}
Since the $f_i$ and the multiplication by $a\not=0$ are isomorphisms at the generic fiber the lemma follows from proposition~\ref{PROP.1.42A}.
\end{proof}

\noindent {\it Remark.}
The dual morphism $\dual{f}$ clearly depends on the choice of $a$ and again there is in general no canonical choice of $a$.

\bigskip

\begin{Definition}
Let $\FF$ and $\FF'$ be abelian $\tau$-sheaves. We call $\FF$ and $\FF'$ \emph{quasi-isogenous} $(\FF\approx\FF')$, if there exists a quasi-isogeny between $\FF$ and $\FF'$.
\end{Definition}

\begin{Corollary}\label{QISOG-EQREL}
The relation $\approx$ is an equivalence relation. \qed
\end{Corollary}

\begin{Proposition}\label{QHOM-QEND-ISOMORPHIC}
Let $\FF$ and $\FF'$ be abelian $\tau$-sheaves. If $\FF\approx\FF'$, then 
\begin{suchthat}
\item the $Q$-algebras $\QEnd(\FF)$ and $\QEnd(\FF')$ are isomorphic,
\item $\QHom(\FF,\FF')$ is free of rank $1$ both as a left module over $\QEnd(\FF')$ and as a right module over $\QEnd(\FF)$. \qed
\end{suchthat}
\end{Proposition}

Next we want to give an alternative description of $\QHom(\FF,\FF')$ similar to the description \cite[Proposition 3.4.5]{Papanikolas} in the case of ``dual $t$-motives''.

\begin{Proposition}\label{PropAltDescrQHom}
Let $\FF$ and $\FF'$ be abelian $\tau$-sheaves of the same weight and characteristic. Then the $Q$-vector space $\QHom(\FF,\FF')$ is canonically isomorphic to the space of morphisms between the fibers at the generic point $\eta$ of $C_L$
\[
\bigl\{\,f_{0,\eta}:\F_{0,\eta}\to\F'_{0,\eta}\enspace \text{with}\enspace f_{0,\eta}\circ\P_{0,\eta}^{-1}\circ\t_{0,\eta}\;=\;(\P'_{0,\eta})^{-1}\circ\t'_{0,\eta}\circ\sigma^\ast(f_{0,\eta})\,\bigr\}\,.
\]
This isomorphism is compatible with composition of quasi-morphisms.
\end{Proposition}

\begin{proof}
Clearly if $f\in \QHom(\FF,\FF')$ the map $f\mapsto f_{0,\eta}$ is a monomorphism of $Q$-vector spaces. To show that it is surjective let $f_{0,\eta}$ as above be given. As in the proof of \ref{PROP.1.42A} choose a finite flat morphism $\varphi:C\to\PP^1_{\Fq}$ with $\varphi^{-1}(\infty_{\PP^1})=\{\infty\}$, set $\Fq[t]=\Gamma(\PP^1_\Fq\setminus\{\infty_{\PP^1}\},\O_{\PP^1})$, and replace $\FF$ and $\FF'$ by $\varphi_\ast\FF$ and $\varphi_\ast\FF'$. Choose $L[t]$-bases of $M=\Gamma(\PP^1_L\setminus\{\infty_{\PP^1}\},\F_0)$ and $M'=\Gamma(\PP^1_L\setminus\{\infty_{\PP^1}\},\F'_0)$, and write $\P_0^{-1}\circ\t_0$ and $(\P'_0)^{-1}\circ\t'_0$ with respect to these bases as matrices $T$ and $T'$ with coefficients in $L[t]$. If $\chr\ne\infty$ let $\theta:=c^\ast(t)$ and set $e:=d$ and $e':=d'$. If $\chr=\infty$ set $e=e'=0$ and $\theta:=0$ (the choice of $\theta$ will not play a role in this case).

Then in both cases $\det T=b\cdot(t-\theta)^e$ and $\det T'=b'\cdot(t-\theta)^{e'}$ for $b,b'\in L^{\times}$. By considering the adjoint matrices we find in particular that $(t-\theta)^eT^{-1}$ and $(t-\theta)^{e'}(T')^{-1}$ have all their coefficients in $L[t]$. Write $f_{0,\eta}$ with respect to these bases as a matrix $F\in M_{r'\times r}(L(t))$. It satisfies $FT=T'\s F$.

Consider the ideals of $L^\alg[t]$ where $L^\alg$ is an algebraic closure of $L$
\[
I:=\bigl\{\,h\in L^\alg[t]:\enspace hF\in M_{r'\times r}(L^\alg[t])\,\bigr\}
\]
and $I^\sigma:=\{\s(h):h\in I\}$. Note that $I\ne(0)$. We claim that
\begin{eqnarray*}
h\in I \enspace& \Longrightarrow & (t-\theta)^{e'}h\in I^\sigma \quad\text{and}\\
h\in I^\sigma &\Longrightarrow & (t-\theta)^e h\in I\,.
\end{eqnarray*}
Indeed, let $h\in I$ and set $g:=(\s)^{-1}((t-\theta)^{e'}h)$. Then 
\[
\s(gF)\;=\;(t-\theta)^{e'}h\,\s F\;=\;(t-\theta)^{e'}(T')^{-1}\cdot hFT\;\in\; M_{r'\times r}(L^\alg[t])\,.
\]
Hence $g\in I$ and $(t-\theta)^{e'}h\in I^\sigma$. Conversely let $h\in I^\sigma$, that is, $h=\s(g)$ for $g\in I$. Then
\[
(t-\theta)^ehF\;=\;T'\s(gF)\cdot(t-\theta)^eT^{-1}\;\in\;M_{r'\times r}(L^\alg[t])
\]
proving the claim.

Since $L^\alg[t]$ is a principal ideal domain we find $I=(h)$ and $I^\sigma=(\s (h))$ for some $h\in I$. In particular $(t-\theta)^{e'}h=g\cdot\s(h)$ and $(t-\theta)^e \s(h)=f\cdot h$ for suitable $f,g\in L^\alg[t]$. We conclude $(t-\theta)^{e+e'}h=fg\,h$ and since the polynomials $h$ and $\s(h)$ are non-zero and have the same degree, $g=\beta\cdot(t-\theta)^{e'}$ for some $\beta\in (L^\alg)^{\times}$. Choose an element $\gamma\in (L^\alg)^{\times}$ with $\gamma^{q-1}=\beta$. Then
\[
a:=\gamma h=\s(\gamma h)\in \Fq[t]
\]
and $aF\in M_{r'\times r}(L[t])$. Thus $f_{0,\eta}$ defines a morphism $f_0:\F_0\to\F'_0(D_0)$ for $D_0:=(\varphi^\ast a)_0+m_0\cdot\infty$ with appropriate $m_0\in \N_0$. Here $(\varphi^\ast a)_0$ is the zero divisor of the element $\varphi^\ast a\in A$.

Now we define inductively on $C'_L:=C_L\setminus\{\infty\}$
\[
f^{\,}_i := \P'_{i-1}\circ f^{\,}_{i-1}\circ\P_{i-1}^{-1}:\; \F_i|_\CLa\rightarrow\F'_i((\varphi^\ast a)_0)|_\CLa \qquad (i>0)
\]
and analogously for $i<0$. To pass to the projective closure, we allow divisors $D_i=(\varphi^\ast a)_0+m_i\cdot\infty$ for sufficiently large $m_i>0$ such that $f_i:\,\F_i\rightarrow\F'_i(D_i)$ for all $i\in\Z$. Since $\FF$ and $\FF'$ have the same weight, we have the periodical identification if $l$ satisfies condition 2 of Definition~\ref{Def1.1} for both $\FF$ and $\FF'$
\[
\xymatrix{
& \F_{i+nl} \ar[r]^<<{\sim}\ar[d]^{f_{i+nl}} & \F_i(nk\cdot\infty) \ar[d]^{f_i\,\otimesidOCL{nk\cdot\infty}} & \ar@{}[d]^{\displaystyle(i,n\in\Z)\,.} \\
& \F'_{i+nl}(D_i) \ar[r]^<<{\sim} & \F'_i(D_i+nk\cdot\infty) &
}
\]
Take $m:=\max\{m_0,\dots,m_{l-1}\}$ and set $D:=(\varphi^\ast a)_0+m\cdot\infty$. Then $f_i:\,\F_i\rightarrow\F'_i(D)$ for all $i\in\Z$. Since the commutation with the $\P$'s and the $\t$'s holds by construction, the collection of the $f_i$ is the desired quasi-morphism $f\in\QHom(\FF,\FF')$. 
\end{proof}

\bigskip

The following proposition connects the theory of quasi-morphisms of abelian $\tau$-sheaves to the theory of morphisms of their associated pure Anderson motives and $\t$-modules.

\begin{Proposition}\label{CONNECTION}
Let $\FF$ and $\FF'$ be two abelian $\tau$-sheaves of the same weight and let $D\subset C$ be a finite closed subscheme as in Section~\ref{SectRelation}. 
\begin{suchthat}
\item
If $\infty\in D$ we have a canonical isomorphism of $Q$-vector spaces
\[
\QHom(\FF,\FF') \;=\; \Hom(\ulM^{(D)}(\FF),\ulM^{(D)}(\FF'))\otimes_{\TA} Q\,.
\]
\item
If $\infty\notin D$ choose an integer $l$ which satisfies condition 2 of \ref{Def1.1} for both $\FF$ and $\FF'$ and assume $\Ff_{q^l}\subset L$. Then we have a canonical isomorphism of $Q$-vector spaces
\[
\QHom(\FF,\FF') \;=\; \Hom_{\P,\Lambda}(\ulM^{(D)}(\FF),\ulM^{(D)}(\FF'))\otimes_{\TA} Q
\]
where the later is the space of all morphisms commuting with $\P$ and $\Lambda(\lambda)$ from (\ref{EQ.Pi+Lambda}) for all $\lambda\in\Ff_{q^l}$.
\end{suchthat}
By lemma~\ref{Lemma2.9a} the condition on the weights can be dropped if $\chr\notin D$ or if $\chr\ne\infty$ and $\infty\in D$.
\end{Proposition}

\begin{proof}
Let $\ulM:=\ulM^{(D)}(\FF)$ and $\ulM':=\ulM^{(D)}(\FF')$. We exhibit a monomorphism of $Q$-vector spaces from $\QHom(\FF,\FF')$ to $\Hom(\ulM,\ulM')\otimes_{\TA}Q$ in case 1 (respectively from $\QHom(\FF,\FF')$ to $\Hom_{\P,\Lambda}(\ulM,\ulM')\otimes_{\TA}Q$ in case 2) and another monomorphism from the target of the first to the $Q$-vector space 
\[
H\;:=\;\bigl\{\,f_{0,\eta}:\F_{0,\eta}\to\F'_{0,\eta}\enspace \text{with}\enspace f_{0,\eta}\circ\P_{0,\eta}^{-1}\circ\t_{0,\eta}\;=\;(\P'_{0,\eta})^{-1}\circ\t'_{0,\eta}\circ\sigma^\ast(f_{0,\eta})\,\bigr\}
\]
introduced in proposition~\ref{PropAltDescrQHom} such that the composition of the two monomorphisms is the isomorphism from \ref{PropAltDescrQHom}.

Let $f\in\QHom(\FF,\FF')$. By the Riemann-Roch Theorem we can find some $a\in Q$ such that $a\cdot f$ maps from $\FF$ into $\FF'(n\cdot D)$ for some $n>0$. Since $a$ and $f$ commute with the $\P$'s and $\t$'s, we get for the first monomorphism
\[
\begin{array}{crccll}
f&\mapsto &a\cdot f_0|_{C_L\setminus D}\otimes a^{-1}&\in&\Hom(\ulM,\ulM')\otimes_A Q&\quad\text{in case 1, and}\\[1mm]
f&\mapsto &a\cdot (f_0\oplus\ldots\oplus f_{l-1})|_{C_L\setminus D}\otimes a^{-1}&\in&\Hom_{\P,\Lambda}(\ulM,\ulM')\otimes_A Q&\quad\text{in case 2.}
\end{array}
\]
\forget{
\[
\begin{array}{rclcrcll}
a\cdot f_0|_{C_L\setminus D}&\in&\Hom(\ulM,\ulM')&\text{ and }& f\mapsto a\cdot f_0|_{C_L\setminus D}\otimes a^{-1}&\in&\Hom(\ulM,\ulM')\otimes_A Q&\quad\text{in case 1, and}\\[1mm]
a\cdot (f_0\oplus\ldots f_{l-1})|_{C_L\setminus D}&\in&\Hom_{\P,\Lambda}(\ulM,\ulM')&\text{ and }& f\mapsto a\cdot (f_0\oplus\ldots f_{l-1})|_{C_L\setminus D}\otimes a^{-1}&\in&\Hom_{\P,\Lambda}(\ulM,\ulM')\otimes_A Q&\quad\text{in case 2}
\end{array}
\]
}
To construct the second monomorphism we treat each case separately.

\smallskip
\noindent
1. The localization $\Hom(\ulM,\ulM')\otimes_{\TA}Q\to H$, $g\otimes a\mapsto ag_\eta$ at the generic point $\eta$ of $C_L$ gives the desired monomorphism.

\smallskip
\noindent
2. Let $(g:\,\bigoplus M_i \rightarrow \bigoplus M_i)\otimes a\in \Hom_{\P,\Lambda}(\ulM,\ulM')\otimes_{\TA}Q$. Then $g$ corresponds to a block matrix \mbox{$(g_{ij}:M_j\to M_i)_{0\le i,j<l}$} with $g_{ij}\cdot\lambda^{q^j} = \lambda^{q^i}\cdot g_{ij}$ for all $\lambda\in \Ff_{q^l}$. Therefore, we have $g_{ij}=0$ for $i\ne j$. We map $g\otimes a$ to the localization $a\cdot(g_{00})_\eta$ at $\eta$. Since $\P g=g\P$ this map is injective and our proof is complete.
\end{proof}

\noindent
{\it Remark.} Again note the importance of the assumption that $\FF$ and $\FF'$ must have the same weight, since otherwise $\QHom(\FF,\FF')=(0)$ by \ref{PROP.1} whereas $\Hom(\ulM^{(D)}(\FF),\ulM^{(D)}(\FF'))$ could be non-zero. Consider for example the abelian $\tau$-sheaves on $C=\PP^1_\Fq$ with $C\setminus\{\infty\}=\Spec \Fq[t]$ given by $\F_i=\O_{\PP^1_L}(i\cdot\infty)$, $\t_i=t$ and $\F'_i=\O_{\PP^1_L}(2i\cdot\infty)$, $\t'_i=t^2$, where $\P$ and $\P'$ are the natural inclusions. They have $\weight(\FF)=1$ and $\weight(\FF')=2$. If we choose $D={\rm V}(t)$ and $z=t^{-1}\in \Gamma(C\setminus D,\O_{\PP^1_\Fq})$ as uniformizing parameter at $\infty$ then $\ulM^{(D)}(\FF)=(L[z],1)=\ulM^{(D)}(\FF')$.

\begin{Definition}\label{Def1.38b}
Let $\ulM$ and $\ulM'$ be pure Anderson motives. Then the elements of $\Hom(\ulM,\ulM')\otimes_A Q$ which admit an inverse in $\Hom(\ulM',\ulM)\otimes_A Q$ are called \emph{quasi-isogenies}.
\end{Definition}

\begin{Corollary}\label{Cor2.9d}
Let the characteristic be different from $\infty$. Then the functor $\FF\mapsto\ulM(\FF)$ defines an equivalence of categories between 
\begin{suchthat}
\item 
the category with abelian $\tau$-sheaves as objects and with $\QHom(\FF,\FF')$ as the set of morphisms,
\item 
 and the category with pure Anderson motives as objects and with $\Hom(\ulM,\ulM')\otimes_A Q$ as the set of morphisms.
\end{suchthat}
We call these the \emph{quasi-isogeny categories} of abelian $\tau$-sheaves of characteristic different from $\infty$ and of pure Anderson motives respectively.
\end{Corollary}

\begin{proof}
This is just a reformulation of Theorem~\ref{PropX.1} and Proposition~\ref{CONNECTION}.
\end{proof}

% -----------------------------------------------------------------------------

\bigskip

\subsection{Simple and semisimple abelian $\tau$-sheaves and pure Anderson motives}

In the last section of this chapter we want to draw some first conclusions in our study of the structure of\/ $\QEnd(\FF)$. %For this purpose, we start considering the simple case.

\begin{Definition}
Let $\FF$ be an abelian $\tau$-sheaf.
\begin{suchthat}
\item $\FF$ is called \emph{simple}, if\/ $\FF\not=\ZZ$ and\/ $\FF$ has no abelian quotient $\tau$-sheaves other than $\ZZ$ and\/ $\FF$.
\item $\FF$ is called \emph{semisimple}, if\/ $\FF$ admits, up to quasi-isogeny, a decomposition into a direct sum $\FF\approx\FF_1\oplus\cdots\oplus\FF_n$ of simple abelian $\tau$-sheaves $\FF_j$ $(1\le j\le n)$.
\item $\FF$ is called \emph{primitive}, if its rank and its dimension are relatively prime.
\end{suchthat}
We make the same definition for a pure Anderson motive.
\end{Definition}

\begin{Remark} \label{RemDualBehaviour}
It is not sensible to try defining \emph{simple} abelian $\tau$-sheaves via abelian sub-$\tau$-sheaves, since for example the shifted abelian $\tau$-sheaf $(\F_{i-n},\P_{i-n},\t_{i-n})$ by $n\in\N$ is always a proper abelian sub-$\tau$-sheaf of $(\F_i,\P_i,\t_i)$. 
Furthermore we have for every positive divisor $D$ on $C$ a strict inclusion $\FF(-D)\subset\FF$. This shows that abelian $\tau$-sheaves behave dually to abelian varieties. Namely an abelian variety is called simple if it has no non-trivial abelian subvarieties.
\end{Remark}

\begin{Proposition}\label{Prop1.29a}
Let $\FF$ be an abelian $\tau$-sheaf with characteristic different from $\infty$. Then $\FF$ is \mbox{(semi-)}simple if and only if the pure Anderson motive $\ulM(\FF)$ is (semi-)simple.
\end{Proposition}

\begin{proof}
First let $\FF$ be simple and let $f:\ulM(\FF)\to \ulM'$ be a surjective morphism of pure Anderson motives. By Theorem~\ref{PropX.1} there is an abelian $\tau$-sheaf $\FF'$ with $\ulM'=\ulM(\FF')$. By \ref{CONNECTION} there is an integer $n$ such that $f\in\Hom\bigl(\FF,\FF'(n\cdot\infty)\bigr)$ and $\im f$ is an abelian quotient $\tau$-sheaf of $\FF$ by \ref{IMISABELIAN}. Hence $f$ is injective or $f=0$ proving that $f:\ulM(\FF)\to\ulM'$ is an isomorphism or $\ulM'=0$.

Conversely let $\ulM(\FF)$ be simple and let $f:\FF\to\FF'$ be an abelian quotient $\tau$-sheaf of $\FF$. Then $\ulM(f):\ulM(\FF)\to\ulM(\FF')$ is surjective. So $f=0$ or $f$ is injective proving that $\FF'=\ZZ$ or $f$ is an isomorphism.

Clearly if $\FF$ is semisimple then so is $\ulM(\FF)$. Conversely if $\ulM(\FF)$ is isogenous to a direct sum $\ulM_1\oplus\ldots\oplus\ulM_n$ with $\ulM_i$ simple, then we obtain from \ref{PropX.1} and \ref{Cor2.9b} simple abelian $\tau$-sheaves $\FF_i$ of the same weight with $\ulM(\FF_i)=\ulM_i$ and $\FF\approx\FF_1\oplus\ldots\oplus\FF_n$ by \ref{Cor2.9d}.
\end{proof}

\begin{Proposition}\label{PROP.5}%
Let $\FF$ be an abelian $\tau$-sheaf. If\/ $\FF$ is primitive, then $\FF$ is simple.
\end{Proposition}

\begin{proof}
Let $\TFF$ be an abelian quotient $\tau$-sheaf of $\FF$. Clearly, we have $\Tr\le r$. If $\Tr=0$, then $\TFF=\ZZ$. Otherwise, the surjection $f\in\Hom(\FF,\TFF)$ is non-zero, and by \ref{PROP.1} we get $\Td r = d\Tr$. Since $r$ and $d$ are relatively prime, it follows $\Tr=r$ and $\Td=d$. Therefore, considering the ranks in $\smallexact{0}{}{\ker f_i}{}{\F_i}{f_i}{\TF_i}{}{0}$, we conclude $\ker f_i=0$ and hence $f_i$ is an isomorphism.
\end{proof}

\begin{Corollary}\label{Cor1.50b}
If $\ulM$ is a pure Anderson motive of rank $r$ and dimension $d$ with $(r,d)=1$ then $\ulM$ is simple.\qed
\end{Corollary}

\begin{Proposition}\label{PROP.1.42B}%
Let $\FF$ and $\FF'$ be abelian $\tau$-sheaves of the same rank and dimension. If the characteristic is different from $\infty$ and if\/ $\FF$ is simple, then every non-zero morphism between $\FF$ and\/ $\FF'$ is an isogeny.
\end{Proposition}

\begin{proof}
Let $f\in\Hom(\FF,\FF')$ be a non-zero morphism. Since the characteristic is different from $\infty$, we know by \ref{IMISABELIAN} that $\im f$ is an abelian quotient $\tau$-sheaf. As $\FF$ is simple, we have $\FF\cong\im f$ and therefore $f$ is injective. Thus, by \ref{PROP.1.42A}, $f$ is an isogeny.
\end{proof}

\begin{Remark}\label{Rem1.51b}
Note that the assumption on the characteristic in the proposition and the  theorem below is essential. Namely, the abelian $\tau$-sheaf $\FF$ of Example~\ref{ExTRichter} is primitive, hence simple, but the endomorphism $f$ of $\FF$ constructed there violates the assertions of the proposition and the theorem below.
\end{Remark}

\begin{Theorem}\label{QEND-DIVISION-MATRIX}
Let $\FF$ be an abelian $\tau$-sheaf of characteristic different from $\infty$.
\begin{suchthat}
\item If\/ $\FF$ is simple, then $\QEnd(\FF)$ is a division algebra over $Q$.
\item If\/ $\FF$ is semisimple with decomposition $\FF\approx\FF_1\oplus\cdots\oplus\FF_n$ up to quasi-isogeny into simple abelian $\tau$-sheaves $\FF_j$, then $\QEnd(\FF)$ decomposes into a finite direct sum of full matrix algebras over the division algebras $\QEnd(\FF_j)$ over $Q$.
\end{suchthat}
\end{Theorem}

\noindent {\it Remark.}
We will show in theorem \ref{Thm3.8} below that over a finite field also the converses to these statements are true.

\begin{proof}
1. We saw in \ref{QEND-Q-ALGEBRA} that $\QEnd(\FF)$ is a $Q$-algebra. By \ref{QISOG-GROUP}, we can invert every quasi-isogeny in $\QIsog(\FF)$. Thus, by proposition \ref{PROP.1.42B}, $\QEnd(\FF)$ is a division algebra.

\smallskip
\noindent
2. Let $\FF\approx\FF_1\oplus\cdots\oplus\FF_n$ be a decomposition into simple abelian $\tau$-sheaves $\FF_j$. By \ref{QHOM-QEND-ISOMORPHIC}, we know that $\QEnd(\FF)\cong\QEnd(\FF_1\oplus\cdots\oplus\FF_n)$, so we just have to consider the decomposition. By proposition \ref{PROP.1.42B}, we only get non-zero morphisms between $\FF_j$ and $\FF_i$, if $\FF_j\approx\FF_i$. Hence we can group the quasi-isogenous $\FF_j$ and decompose $\QEnd(\FF)=E_1\oplus\cdots\oplus E_m$ into their collective endomorphism rings $E_\nu=\QEnd(\FF_{j_{\nu,1}}\oplus\cdots\oplus\FF_{j_{\nu,\mu_\nu}})$, $1\le\nu\le m$, $\sum_{\nu=1}^m \mu_\nu = n$. By \ref{QHOM-QEND-ISOMORPHIC}, every $\QHom(\FF_{j_{\nu,\alpha}},\FF_{j_{\nu,\beta}})$, $1\le\alpha,\beta\le\mu_\nu$, is isomorphic to $\QEnd(\FF_{j_{\nu,1}})$. Hence we conclude that each $E_\nu$ is isomorphic to a matrix algebra over $\QEnd(\FF_{j_{\nu,1}})$ which completes the proof.
\end{proof}

For example if $\FF$ is an abelian $\tau$-sheaf associated to a Drinfeld module, then  $d=1$ and $\FF$ is primitive, hence simple. Also $\chr\ne\infty$ and so $\QEnd(\FF)$ is a division algebra. Together with Corollary~\ref{Cor2.9d} this gives another proof that the endomorphism $Q$-algebra of a Drinfeld module is a division algebra over $Q$.

% =============================================================================

\section{The Associated Tate Modules and Local Shtukas} \label{Chapt2}

% =============================================================================

\subsection{Local Shtukas}\label{SectLS}

Before treating Tate modules in Section~\ref{SectTateModules} we want to attach another local structure to abelian $\tau$-sheaves or pure Anderson motives which is in a sense intermediate on the way to the $v$-adic Tate module, namely the local (iso-)shtuka at $v$. It is the analog of the Dieudonn\'e module of the $p$-divisible group attached to an abelian variety. Note however one fundamental difference. While the Dieudonn\'e module exists only if $p$ equals the characteristic of the base field, there is no such restriction in our theory here. And in fact this would even allow to dispense with Tate modules at all and only work with local (iso-)shtukas. Being not so radical here we shall nevertheless prove the standard facts about Tate modules through the use of local (iso-)shtukas.

To give the definition we introduce the following notation. Let $v\in C$ be a place of $Q$ and let $L\supset\Fq$ be a field. Recall that $A_{v,L}$ denotes the completion of $\O_{C_L}$ at the closed subscheme $v\times\Spec L$ and that $Q_{v,L}=A_{v,L}[\frac{1}{v}]$. Note that $v\times\Spec L$ may consist of more than one point if the intersection of $L$ with the residue field of $v$ is larger than $\Fq$. Then $A_{v,L}$ is not an integral domain and $Q_{v,L}$ is not a field. Local (iso-)shtukas were introduced in \cite{Hl} under the name \emph{Dieudonn\'e $\Fq\dbl z\dbr$-modules} (respectively \emph{Dieudonn\'e $\Fq\dpl z\dpr$-modules}). They are studied in detail in \cite{Anderson2,Kim}; see also \cite{HartlPSp}. Over a field their definition takes the following form.

\begin{Definition}\label{DefLS1}
An \emph{(effective) local $\sigma$-shtuka} at $v$ of rank $r$ over $L$ is a pair $\ulHM=(\hat M,\phi)$ consisting of a free $A_{v,L}$-module $\hat M$ of rank $r$ and an injective $A_{v,L}$-module homomorphism \mbox{$\s\hat M\to \hat M$}.

A \emph{local $\sigma$-isoshtuka} at $v$ of rank $r$ over $L$ is a pair $\ulHN=(\hat N,\phi)$ consisting of a free $Q_{v,L}$-module $\hat N$ of rank $r$ and an isomorphism $\phi:\s\hat N\to\hat N$ of $Q_{v,L}$-modules.

A \emph{morphism} between two local $\sigma$-shtukas $(\hat M,\phi)$ and $(\hat M',\phi')$ at $v$ is an $A_{v,L}$-homomorphism $f:\hat M\to\hat M'$ with $f\phi=\phi'\s(f)$. We denote the set of morphisms from $\ulHM\to\ulHM'$ by $\Hom_{A_{v,L}[\phi]}(\ulHM,\ulHM')$. The similar definition and notation applies to local isoshtukas.
\end{Definition}

\begin{Remark}
Note that so far in the literature \cite{Anderson2,Hl,HartlPSp, Kim, Laumon} it is always assumed that $A_v$ has residue field $\Fq$, the fixed field of $\sigma$ on $L$. So in particular $A_{v,L}$ is an integral domain and $Q_{v,L}$ is a field. For applications to pure Anderson motives this is not a problem since we may reduce to this case by Propositions~\ref{PropLS3} and \ref{PropLS4} below.
\end{Remark}

\begin{Definition}\label{DefLS1b}
A local shtuka $\ulHM=(\hat M,\phi)$ is called \emph{\'etale} if $\phi$ is an isomorphism. The \emph{Tate module} of an \'etale local $\sigma$-shtuka $\ulHM$ at $v$ is the $G:=\Gal(L^\sep/L)$-module of $\phi$-invariants
\[
T_v\ulHM\;:=\;\bigl(\ulHM\otimes_{A_{v,L}}A_{v,L^\sep}\bigr)^\phi\,.
\]
The \emph{rational Tate module} of $\ulHM$ is the $G$-module
\[
V_v\ulHM\;:=\;T_v\ulHM\otimes_{A_v}Q_v\,.
\]
\end{Definition}

It follows from \cite[Proposition~6.1]{TW} that $T_v\ulHM$ is a free $A_v$-module of the same rank than $\ulHM$ and that the natural morphism
\[
T_v\ulHM\otimes_{A_v}A_{v,L^\sep}\isoto \ulHM\otimes_{A_{v,L}}A_{v,L^\sep}
\]
is a $G$- and $\phi$-equivariant isomorphism of $A_{v,L^\sep}$-modules, where on the left module $G$-acts on both factors and $\phi$ is $\id\otimes\s$. Since $(L^\sep)^{G}=L$ we obtain

\begin{Proposition}\label{Prop2.13'}
Let $\ulHM$ and $\ulHM{}'$ be \emph{\'etale} local $\sigma$-shtukas at $v$ over $L$. Then
\begin{suchthat}
\item 
$\DS\ulHM\;=\;(T_v\ulHM\otimes_{A_v}A_{v,L^\sep})^{G}$, the Galois invariants,
\item 
$\DS\Hom_{A_{v,L}[\phi]}(\ulHM,\ulHM{}')\;\isoto\;\Hom_{A_v[G]}(T_v\ulHM,T_v\ulHM{}')\;,\es f\mapsto T_vf$ is an isomorphism.
\end{suchthat}
In particular the Tate module functor yields a fully faithful embedding of the category of \'etale local shtukas at $v$ over $L$ into the category of $A_v[G]$-modules, which are finite free over $A_v$. \qed
\end{Proposition}

\bigskip

If the residue field $\BF_v$ of $v$ is larger than $\Fq$ one has to be a bit careful with local (iso-)shtukas since the ring $Q_{v,L}$ is then in general not a field. Namely let $\#\BF_v=q^n$ and let $\BF_{q^f}:=\{\alpha\in L:\alpha^{q^n}=\alpha\}$ be the ``intersection'' of $\BF_v$ with $L$. Choose and fix an $\BF_q$-homomorphism $\BF_{q^f}\hookrightarrow\BF_v\subset A_v$. Then
\[
\BF_v\otimes_\Fq L = \prod_{\Gal(\BF_{q^f}/\Fq)}\BF_v\otimes_{\BF_{q^f}}L = 
\prod_{i\in\Z/f\Z}\BF_v\otimes_\Fq L\,/\,(b\otimes 1-1\otimes b^{q^i}:b\in \BF_{q^f})
\]
and $\s$ transports the $i$-th factor to the $(i+1)$-th factor. Denote by $\Fa_i$ the ideal of $A_{v,L}$ (or $Q_{v,L}$) generated by $\{b\otimes 1-1\otimes b^{q^i}:b\in \BF_{q^f}\}$. Then
\[
A_{v,L} = \prod_{\Gal(\BF_{q^f}/\Fq)}A_v\wh\otimes_{\BF_{q^f}}L = 
\prod_{i\in\Z/f\Z}A_{v,L}\,/\,\Fa_i
\]
and similarly for $Q_{v,L}$.
Note that the factors in this decomposition and the ideals $\Fa_i$ correspond precisely to the places $v_i$ of $C_{\BF_{q^f}}$ lying above $v$.

\begin{Proposition}\label{PropLS3}
Fix an $i$. The reduction modulo $\Fa_i$ induces equivalences of categories
\begin{suchthat}
\item
$\DS (\hat N,\phi)\longmapsto \bigl(\hat N/\Fa_i\hat N\,,\es\phi^f\mod\Fa_i:(\s)^f \hat N/\Fa_i\hat N \to \hat N/\Fa_i\hat N\bigr)$\\[2mm]
between local $\sigma$-isoshtukas at $v$ over $L$ and local $\sigma^f$-isoshtukas at $v_i$ over $L$ of the same rank.
\item 
\vspace{2mm}
$\DS (\hat M,\phi)\longmapsto \bigl(\hat M/\Fa_i\hat M\,,\es\phi^f\mod\Fa_i:(\s)^f \hat M/\Fa_i\hat M \to \hat M/\Fa_i\hat M\bigr)$\\[2mm]
between\/ \emph{\'etale} local $\sigma$-shtukas at $v$ over $L$ and \emph{\'etale} local $\sigma^f$-shtukas at $v_i$ over $L$ preserving Tate modules
\[
T_v(\hat M,\phi) \isoto T_{v_i}(\hat M/\Fa_i\hat M,\phi^f)\,.
\]
\end{suchthat}
\end{Proposition}

\begin{proof}
Since $\s\Fa_{i-1}=\Fa_{i}$ the isomorphism $\phi$ yields for every $i$ an isomorphism $\phi\mod\Fa_i:\s(\hat N/\Fa_{i-1}\hat N)\to \hat N/\Fa_{i}\hat N$ and similarly for $\hat M$. These allow to reconstruct the other factors from $(\hat N/\Fa_i\hat N,\,\phi^f\mod\Fa_i)$. More precisely we describe the quasi-inverse functor. Let $\ulHN'=(\hat N'\,,\,\psi\!:\!(\s)^f\hat N'\isoto\hat N')$ be a local $\sigma^f$-isoshtuka at $v_i$ over $L$. Define the $Q_{v,L}/\Fa_{i+j}$-module $\hat N_{i+j}:=(\s)^j\hat N'$ for $0\le j<f$ and the $Q_{v,L}/\Fa_{i+j}$-homomorphism 
\[
\begin{array}{rcc@{\!\,}cll}
\phi_{i+j}&:=&\id_{\hat N_{i+j}}&:&\s\hat N_{i+j-1}&\isoto\hat N_{i+j} \qquad\qquad\text{for $0<j< f$ and}\\
\phi_i&:=&\psi&:&\s\hat N_{i+f-1}=(\s)^f\hat N'&\isoto\hat N_i\,. 
\end{array}
\]
The quasi-inverse functor sends $\ulHN'$ to the local $\sigma$-isoshtuka $(\bigoplus_{0\le j<f}\hat N_{i+j}\,,\,\bigoplus_{0\le j<f}\phi_{i+j})$ at $v$ over $L$. Reducing the latter modulo $\Fa_i$ clearly gives back $\ulHN'$. Also note that this quasi-inverse functor sends a morphism $h'$ of local $\sigma^f$-isoshtukas at $v_i$ to the morphism $h:=\bigoplus_{0\le i<f}(\s)^j(h')$ of the corresponding $\sigma$-isoshtukas at $v$.

It remains to show that this functor is indeed quasi-inverse to the reduction modulo $\Fa_i$ functor. For this we need that $\phi\mod\Fa_{i+j}$ is an isomorphism for all $0<j<f$. Namely the required isomorphism
\[
\Bigl(\bigoplus_{0\le j<f}(\s)^j(\hat N/\Fa_i\hat N)\,,\,(\phi^f\mod\Fa_i)\oplus\bigoplus_{0<j<f}\id\Bigr)\;\isoto\;\Bigl(\bigoplus_{0\le j<f}\hat N/\Fa_{i+j}\hat N\,,\,\bigoplus_{0\le j<f}\phi\mod\Fa_{i+j}\Bigr)\;=\;(\hat N,\phi)
\]
is given by $\bigoplus_{0\le j<f}\phi^j\mod\Fa_{i+j}$ and one easily checks that this is a natural transformation. Note that the entry for $j=0$ is $\id_{\hat N/\Fa_i\hat N}$. So we \emph{do not need} that $\phi\mod\Fa_i$ is an isomorphism. Also if $\phi\mod\Fa_i$ is not an isomorphism then $\phi^j\mod\Fa_{i+l}$ is still an isomorphism for $l=j$, but not for $0\le l<j$ which is harmless. This will be crucial in the variant which we prove in Proposition~\ref{PropLS4} below.

For \'etale local shtukas we can use the same argument because there again all Frobenius maps are isomorphism. Finally, the isomorphism between the Tate modules follows from the observation that an element $(x_j)_{j\in\Z/f\Z}$ is $\phi$-invariant if and only if $x_{j+1}=\phi(\s x_j)$ for all $j$ and $x_i=\phi^f((\s)^f x_i)$.
\end{proof}

{\it Remark.} \es The advantage of the (\'etale) local $\sigma^f$-(iso-)shtuka at $v_i$ is that it is a free module over $A_{v,L}/\Fa_i=A_v\wh\otimes_{\BF_{q^f}}L$, and the later ring is an integral domain. So the results from \cite{Anderson2,Hl,HartlPSp, Kim, Laumon} apply.

\bigskip

Now let $\FF$ be an abelian $\tau$-sheaf and $v\in C$ an arbitrary place of $Q$. We define the \emph{local $\sigma$-isoshtuka of $\FF$ at $v$} as 
\[
\ulN_v(\FF)\;:=\;\Bigl(\F_0\otimes_{\O_{C_L}}Q_{v,L}\,,\,\P_0^{-1}\circ\t_0\Bigr)\,.
\]
If $v\ne\infty$ we define the \emph{local $\sigma$-shtuka of $\FF$ at $v$} as
\[
\ulM_v(\FF)\;:=\;\Bigl(\F_0\otimes_{\O_{C_L}}A_{v,L}\,,\,\P_0^{-1}\circ\t_0\Bigr)\,.
\]
Likewise if $\ulM$ is a pure Anderson motive over $L$ and $v\in \Spec A$ we define the \emph{local $\sigma$-(iso-)shtuka of $\ulM$ at $v$} as
\[
\ulM_v(\ulM)\;:=\;\ulM\otimes_{A_L}A_{v,L} \qquad\text{respectively}\qquad \ulN_v(\ulM)\;:=\;\ulM\otimes_{A_L}Q_{v,L}\,.
\]
These local (iso-)shtukas all have rank $r$. The local shtukas are \'etale if $v\ne\chr$. For $v=\infty$ we also define $\ulN_\infty(\ulM)$ in the same way. Note that $\ulN_\infty(\FF)$ and $\ulN_\infty(\ulM)$ do not contain a local $\sigma$-shtuka since they are isoclinic of slope $-\weight(\FF)<0$. 

However, if $v=\infty$ the periodicity condition allows to define a different local (iso-)shtuka at $\infty$ which is of slope $\ge0$. Namely, choose a finite closed subscheme $D\subset C$ as in Section~\ref{SectRelation} with $\infty\notin D$ and set $\TA=\Gamma(C\setminus D,\O_C)$. We define the \emph{big local $\sigma$-shtuka of $\FF$ at $\infty$} as 
\[
\ulTM_\infty(\FF)\;:=\;\ulM^{(D)}(\FF)\otimes_{\TA\otimes_\Fq L}A_{v,L}\;=\;\bigoplus_{i=0}^{l-1}\F_i\otimes_{\O_{C_L}}A_{v,L}
\]
with $\tau$ from (\ref{EQ.Phi}), and the \emph{big local $\sigma$-isoshtuka of $\FF$ at $\infty$} as
\[
\ulTN_\infty(\FF)\;:=\;\ulTM_\infty(\FF)\otimes_{A_{\infty,L}} Q_{\infty,L}\,.
\]
Both have rank $rl$ and depend on the choice of $k,l$ and $z$. If the characteristic is different from $\infty$ then $\ulTM_\infty(\FF)$ is \'etale. Note that $\ulTM_\infty(\FF)$ and $\ulTN_\infty(\FF)$ were used in \cite{Hl} to construct the uniformization at $\infty$ of the moduli spaces of abelian $\tau$-sheaves.

The big local (iso-)shtukas at $\infty$, $\ulTM_\infty(\FF)$ and $\ulTN_\infty(\FF)$ are always equipped with the automorphisms $\P$ and $\Lambda(\lambda)$ for $\lambda\in \Ff_{q^l}\cap L$ from (\ref{EQ.Pi+Lambda}). We let $\Delta_\infty$ be ``the'' central division algebra over $Q_\infty$ of rank $l$ with Hasse invariant $-\frac{k}{l}$, or explicitly
\begin{equation}\label{EqDelta}
\Delta_\infty \,:=\, \Ff_{q^l}\dpl z\dpr\II{\P}\,/\,(\P^l-z^k,\, \lambda z-z\lambda,\, \P\lambda^q-\lambda\P\text{ for all }\lambda\in \Ff_{q^l})\,.
\end{equation}
If $\Ff_{q^l}\subset L$ we identify $\Delta_\infty$ with a subalgebra of $\End_{Q_{\infty,L}[\phi]}\bigl(\ulTN_\infty(\FF)\bigr)$ by mapping $\lambda\in\BF_{q^l}\subset\Delta_\infty$ to $\Lambda(\lambda)$.

\medskip

\forget{
Now let $\FF$ be an abelian $\tau$-sheaf and $v\in C$ an arbitrary place of $Q$. Choose a finite closed subscheme $D\subset C$ as in section~\ref{SectRelation} and set $\TA=\Gamma(C\setminus D,\O_C)$. If $v\ne\infty$ we choose $D=\{\infty\}$, then $\ulM^{(D)}(\FF)$ is the pure Anderson motive associated with $\FF$. We define the \emph{local $\sigma$-shtuka of $\FF$ at $v$} as 
\[
\ulM_v(\FF)\;:=\;\ulM^{(D)}(\FF)\otimes_{\TA\otimes_\Fq L}A_{v,L}
\]
and the \emph{local $\sigma$-isoshtuka of $\FF$ at $v$} as
\[
\ulN_v(\FF)\;:=\;\ulM^{(D)}(\FF)\otimes_{\TA\otimes_\Fq L} Q_{v,L}\,.
\]
Likewise if $\ulM$ is a pure Anderson motive over $L$ and $v\in \Spec A$ we define the \emph{local $\sigma$-(iso-)shtuka of $\ulM$ at $v$} as
\[
\ulM_v(\ulM)\;:=\;\ulM\otimes_{A_L}A_{v,L} \qquad\text{respectively}\qquad \ulN_v(\ulM)\;:=\;\ulM\otimes_{A_L}Q_{v,L}\,.
\]
These local shtukas are \'etale if $v\ne\chr$, even if $v=\infty$. In the later case we will also need another local isoshtuka, namely
\[
\ulTN_\infty(\FF)\;:=\;\bigl(\F_0\otimes_{\O_{C_L}}Q_{\infty,L}\,,\,\P_0^{-1}\circ\t_0\bigr)
\]
for an abelian $\tau$-sheaf $\FF$ and
\[
\ulTN_\infty(\ulM)\;:=\;\ulM\otimes_{A_L}Q_{\infty,L}
\]
for a pure Anderson motive $\ulM$. In contrast to $\ulN_\infty(\FF)$ and $\ulN_\infty(\ulM)$ the later do not contain \'etale local shtukas (except for $\ulTN_\infty(\FF)$ when $\chr=\infty$).

\medskip
}

\begin{Theorem}\label{ThmLS1}
Let $\FF$ and $\FF'$ be abelian $\tau$-sheaves of the same weight over a \emph{finite} field $L$ and let $v$ be an arbitrary place of $Q$. 
\begin{suchthat}
\item 
Then there is a canonical isomorphism of $Q_v$-vector spaces
\[
\QHom(\FF,\FF')\otimes_Q Q_v \isoto \Hom_{Q_{v,L}[\phi]}\bigl(\ulN_v(\FF),\ulN_v(\FF')\bigr)\,.
\]
\item 
If $v=\infty$ choose an $l$ which satisfies \ref{Def1.1}/2 for both $\FF$ and $\FF'$ and assume $\BF_{q^l}\subset L$. Then there is a canonical isomorphism of $Q_\infty$-vector spaces
\[
\QHom(\FF,\FF')\otimes_Q Q_\infty \isoto \Hom_{\Delta_\infty\wh\otimes_\Fq L[\phi]}\bigl(\ulTN_\infty(\FF),\ulTN_\infty(\FF')\bigr)\,.
\]
\end{suchthat}
\end{Theorem}

\begin{proof}
1. Since in the notation of proposition~\ref{PropAltDescrQHom} the condition
\begin{equation}
\label{EqLS1}
f_{0,\eta}\circ\P_{0,\eta}^{-1}\circ\t_{0,\eta}-(\P'_{0,\eta})^{-1}\circ\t'_{0,\eta}\circ\s(f_{0,\eta})\;=\;0  
\end{equation}
is $Q$-linear in $f_{0,\eta}$ we see that the left hand side of the asserted isomorphy is
\[
\{\,f_{0,\eta}:\F_{0,\eta}\otimes_Q Q_v \to \F'_{0,\eta}\otimes_Q Q_v \,|\es f_{0,\eta}\text{ satisfies (\ref{EqLS1})}\,\}\,.
\]
Since $L/\Fq$ is finite, $Q_v\otimes_\Fq L=Q_{v,L}$ and $\F_{0,\eta}\otimes_Q Q_v$ equals $\ulN_v(\FF)$, and 1 is proved.

\smallskip
\noindent
2. Consider the isomorphism 
\[
\Hom\bigl(\ulM^{(D)}(\FF),\ulM^{(D)}(\FF')\bigr)\otimes_{\TA} Q_\infty\es\cong\es\Hom_{Q_{\infty,L}[\phi]}\bigl(\ulTN_\infty(\FF),\ulTN_\infty(\FF')\bigr)
\]
whose existence is proved as in 1. Now 2 follows by applying \ref{CONNECTION} and noting that the commutation with $\P$ and $\Lambda(\lambda)$ cuts out linear subspaces on both sides which become isomorphic.
\end{proof}

\begin{Theorem}\label{ThmLS2}
Let $\ulM$ and $\ulM'$ be pure Anderson motives over a finite field $L$ and let $v\in\Spec A$ be an arbitrary maximal ideal. Then
\[
\Hom(\ulM,\ulM')\otimes_A A_v \isoto \Hom_{A_{v,L}[\phi]}(\ulM_v(\ulM),\ulM_v(\ulM'))\,.
\]
\end{Theorem}

\begin{proof}
The argument of the previous theorem also works here since $A_v$ is flat over $A$.
\end{proof}

\noindent
{\it Remark.}
If one restricts to places $v\ne\chr$, where the local $\sigma$-shtuka is \'etale, Theorems~\ref{ThmLS1} and \ref{ThmLS2} even hold for finitely generated fields. This was shown by Tamagawa~\cite{Tam}; see also Corollary~\ref{Cor2.17'} below.

\medskip

Let now the characteristic be finite and $v=\chr$ be the characteristic point. Consider an abelian $\tau$-sheaf $\FF$ of characteristic $c$, its local $\sigma$-shtuka $\ulM_\chr(\FF)=(\hat M,\phi)$ at $\chr$ and the decomposition of the later described before proposition~\ref{PropLS3}
\[
\ulM_\chr(\FF)=\prod_{i\in\Z/f\Z}\ulM_\chr(\FF)/\Fa_i\ulM_\chr(\FF)\,.
\]
From the morphism $c:\Spec L\to\Spec\BF_\chr\subset C$ we obtain a canonical $\BF_q$-homomorphism $c^\ast:\BF_\chr\hookrightarrow L$, $f=[\BF_\chr:\Fq]$ and a distinguished place $v_0$ of $C_{\BF_\chr}$ above $v=\chr$, namely the image of $c\times c:\Spec L\to C\times\BF_\chr$. Then $\coker\phi$ on $\ulM_\chr(\FF)$ is annihilated by a power of $\Fa_0$ and therefore $\phi$ has no cokernel on $\ulM_\chr(\FF)/\Fa_i\ulM_\chr(\FF)$ for $i\ne0$ and the proof of proposition~\ref{PropLS3} yields

\begin{Proposition}\label{PropLS4}
The reduction modulo $\Fa_0$ 
\[
\ulM_\chr(\FF)\longmapsto \bigl(\ulM_\chr(\FF)/\Fa_0\ulM_\chr(\FF)\,,\,\phi^f\bigr)
\]
induces an equivalence of categories between the category of local
$\sigma$-shtukas at $\chr$ associated with abelian $\tau$-sheaves of
characteristic $c$ and the category of local $\sigma^f$-shtukas at $v_0$ associated with abelian $\tau$-sheaves of characteristic $c$. 
The same is true for pure Anderson motives. \qed
\end{Proposition}

{\it Remark.} \es Now the fixed field of $\sigma^f$ on $L$ equals $\BF_\chr$, the residue field of $A_\chr$. Also $\ulM_\chr(\FF)/\Fa_0\ulM_\chr(\FF)$ is a module over the integral domain $A_\chr\wh\otimes_{\BF_\chr}L$. So again \cite{Anderson2,Hl,HartlPSp, Kim, Laumon} apply to $\bigl(\ulM_\chr(\FF)/\Fa_0\ulM_\chr(\FF),\phi^f\bigr)$.

\begin{Proposition}\label{Prop2.18b}
Let $\ulM$ be a pure Anderson motive over $L$ and let $\ulHM{}'_\chr$ be a local $\sigma^f$-subshtuka of $\ulM_\chr(\ulM)/\Fa_0\ulM_\chr(\ulM)$ of the same rank. Then there is a pure Anderson motive $\ulM'$ and an isogeny $f:\ulM'\to\ulM$ with $\ulM_\chr(f)\bigl(\ulM_\chr(\ulM')/\Fa_0\ulM_\chr(\ulM')\bigr)=\ulHM{}'_\chr$. The same is true for abelian $\tau$-sheaves.
\end{Proposition}

\begin{proof}
Extend $\ulHM{}'_\chr$ to the local $\sigma$-subshtuka $\bigoplus_{i\in\Z/f\Z}\phi^i\bigl((\s)^i\ulHM{}'_\chr\bigr)$ of $\ulM_\chr(\ulM)$ and consider 
\[
\ulK\;:=\;\ulM_\chr(\ulM)\,/\,\bigoplus_{i\in\Z/f\Z}\phi^i\bigl((\s)^i\ulHM{}'_\chr\bigr)\,.
\]
The induced morphism $\phi_K:\s K\to K$ has its kernel and cokernel supported on the graph of $c$. Set $\ulM'=(M',\tau'):=\bigl(\ker(\ulM\to\ulK),\tau|_{M'}\bigr)$. Then $\ulM'$ is a pure Anderson motive with the required properties by Proposition~\ref{Prop1.5b}.
\end{proof}

There is a corresponding result at the places $v\ne\chr$ which is stated in Proposition~\ref{Prop2.7b}.

\begin{Theorem}\label{ThmW5.2}
For pure Anderson motives over a finite field, being isogenous via a separable isogeny is an equivalence relation.
\end{Theorem}

\begin{proof} 
(cf.~\cite[Theorem 5.2]{Wat})
Since the composition of separable isogenies is again separable we only need to prove symmetry. So let $f:\ulM'\to \ulM$ be a separable isogeny. If the support of $\coker f$ does not meet $\chr$ we can find a dual isogeny which is separable by Remark~\ref{Rem1.26'}. In general we write $\coker f=\ulK^\chr\oplus\ulK_\chr$ with $\Spec(A_L/\chr A_L)\cap\supp\ulK^\chr=\emptyset$ and $\supp\ulK_\chr\subset\Spec(A_L/\chr A_L)$. We factor $f$ as $\ulM'\to\ulM''\to\ulM$ with $\ulM'':=\ker(\ulM\shortonto\coker f\shortonto\ulK_\chr)$ according to Proposition~\ref{Prop1.5b}. By the above we may replace $\ulM'$ by $\ulM''$ and are reduced to the case where $\supp(\coker f)\subset\Spec(A_L/\chr A_L)$. There is a power of $\chr$ which is principal $\chr^n=aA$ for $a\in A$ such that $a$ annihilates $\coker f$. Since our base field is perfect, Lemma~\ref{Lemma1.5d} yields a decomposition
\[
\ulM/a\ulM=(\ulM/a\ulM)^\et\oplus(\ulM/a\ulM)^\nil\,.
\]
We let $\ulM'':=\ker\bigl(\ulM\shortonto(\ulM/a\ulM)^\et\bigr)$ and consider the factorization of $a\cdot\id_\ulM$
\[
\ulM\longto\ulM''\xrightarrow{\es h\;}\ulM'\xrightarrow{\es f\;}\ulM
\]
obtained from the natural surjection $(\ulM/a\ulM)^\et\shortonto \coker f$. Clearly $\coker h$ equals the kernel $\ker\bigl((\ulM/a\ulM)^\et\shortonto\coker f\bigr)$ and $h$ is separable. 

Consider the local $\sigma$-shtuka $\ulM_\chr(\ulM)$ at $\chr$ and the associated local $\sigma^f$-shtuka $\ulM_{v_0}(\ulM):=\ulM_\chr(\ulM)/\Fa_0\ulM_\chr(\ulM)$ from Proposition~\ref{PropLS4}. By \cite[Proposition 2.4.6]{Laumon} the later decomposes $\ulM_{v_0}(\ulM)=\ulM_{v_0}(\ulM)^\et\oplus\ulM_{v_0}(\ulM)^\nil$ into an \'etale part $\ulM_{v_0}(\ulM)^\et$ on which $\sigma^f$ is an isomorphism and a nilpotent part $\ulM_{v_0}(\ulM)^\nil$ on which $\sigma^f$ is topologically nilpotent. Via \ref{PropLS4} we obtain the induced decomposition $\ulM_\chr(\ulM)=\ulM_\chr(\ulM)^\et\oplus\ulM_\chr(\ulM)^\nil$ in which again $\phi$ is an isomorphism on $\ulM_\chr(\ulM)^\et$ and topologically nilpotent on $\ulM_\chr(\ulM)^\nil$. By construction $(\ulM/a\ulM)^\et=\ulM_\chr(\ulM)/\ulM_\chr(\ulM'')$ and
$\ulM_\chr(\ulM'')=a\cdot\ulM_\chr(\ulM)^\et\oplus\ulM_\chr(\ulM)^\nil$. The later is isomorphic to $\ulM_\chr(\ulM)$ as $A_{\chr,L}[\phi]$-module, so $\Hom_{A_{\chr,L}[\phi]}\bigl(\ulM_\chr(\ulM)\,,\,\ulM_\chr(\ulM'')\bigr)$ contains an isomorphism. Since the set of isomorphisms is open we find by Theorem~\ref{ThmLS2} an isogeny $g:\ulM\to\ulM''$ with $\ulM_\chr(g)$ an isomorphism (here we use the assumption that the base field is finite). In particular $g$ is separable and $h\circ g:\ulM\to\ulM'$ is the desired separable isogeny.
\end{proof}

\begin{Example} \label{Ex8.10}
We give an example showing that the preceding theorem is false over infinite fields. This parallels the situation for abelian varieties. Let $C=\PP^1_\Fq$, $A=\Fq[t]$, and $L=\Fq(\gamma)$ where $\gamma$ is transcendental over $\Fq$. Set
\[
T\;:=\;\matr{t+1}{\gamma^{-q}}{\es -\gamma^q}{\es t-1} \qquad\text{and}\qquad T'\;:=\;\matr{t+1-\gamma^{1-q}}{\gamma^{-q}t}{\quad\gamma^q-\gamma}{\quad t}
\]
and consider the pure Anderson motives $\ulM=(L[t]^2,\t=T)$ and $\ulM'=(L[t]^2,\t'=T')$ of characteristic $c^\ast:A\to L, t\mapsto 0$. There is a separable isogeny $f:\ulM'\to\ulM$ given by
\[
\xymatrix @R=0pc
{0 \ar[r] & \ulM' \ar[r]^{\matr{t}{0}{\gamma}{1}} & \ulM \ar[r] &\bigl(L, \t=(1-\gamma^{1-q})\bigr) =\coker f \\
& & {x\choose y} \ar@{|->}[r] & **{ !R(0.45) =<8.5pc,0pc>} \objectbox{(x-\gamma y) \mod t\quad.}
}
\]
We claim that $\End(\ulM')=A=\Fq[t]$. From this it will follow that there is no separable isogeny $g:\ulM\to\ulM'$. Indeed, assume there exists a separable $g$. Then $gf\in \End(\ulM')=A$ is also separable. But $gf$ is not an isomorphism on $\ulM'/t\ulM'$ since already $f$ is not injective modulo $t$. Therefore $gf$ is divisible by $t$, which is not separable. This contradicts the separability of $gf$.

It remains to prove the claim $\End(\ulM')=\Fq[t]$. Let $\sum_{i\ge0}\matr{a_i}{c_i}{b_i}{d_i}t^i\in\End(\ulM')$. The commutation with $\t'$ yields the equations
\begin{eqnarray*}
a_{i-1}+(1-\gamma^{1-q})a_i + \gamma^{-q}b_{i-1} & = & a_{i-1}^q+(1-\gamma^{1-q})a_i^q+(\gamma^q-\gamma)c_i^q\,, \\
(\gamma^q-\gamma)a_i+b_{i-1} & = & b_{i-1}^q+(1-\gamma^{1-q})b_i^q+(\gamma^q-\gamma)d_i^q\,, \\
c_{i-1}+(1-\gamma^{1-q})c_i+\gamma^{-q}d_{i-1} & = & \gamma^{-q}a_{i-1}^q+c_{i-1}^q \,,\\
(\gamma^q-\gamma)c_i+d_{i-1} & = & \gamma^{-q} b_{i-1}^q+d_{i-1}^q\,.
\end{eqnarray*}
For $i=0$ one obtains $c_0=0$ and $a_0\in\Fq$. By subtracting the endomorphism $\matr{a_0}{0}{0}{a_0}$ we may assume that $a_0=0$ and hence $b_0=-\gamma d_0$. When $i=1$ we multiply the first equation by $\gamma^q$ and subtract the second to obtain
\[
b_0^q\;=\;(\gamma^q-\gamma)\bigl(a_1+\gamma c_1-\gamma^{-1}b_1-d_1\bigr)^q\,.
\]
Since $\gamma^q-\gamma$ is not a $q$-th power in $L$ we must have $b_0=d_0=0$ and iterating this argument proves the claim.\qed
\end{Example}

% =============================================================================

\bigskip

\subsection{The Tate Conjecture for pure Anderson motives and Abelian $\tau$-Sheaves}\label{SectTateModules}

In this section we define Tate modules for pure Anderson motives and abelian $\tau$-sheaves and we prove the standard facts on the finiteness of the $A$-module $\Hom(\ulM,\ulM')$ and its relation with Tate modules by using local (iso-)shtukas. We also state the analog of the Tate conjecture for abelian varieties, which was proved by Taguchi~\cite{Taguchi95b} and Tamagawa~\cite{Tam}.

\begin{Definition}
Let $\ulM$ be a\/ $\t$-module on $\TA$ over $L$ (Definition~\ref{Def1.16}) and let\/ $v\in\Spec\TA$ such that the support of\/ $\coker\t$ does not meet\/ $v$. We set
\[
T_v\ulM \;:=\; \liminv{n\in\N}((M/v^nM)\otimes_L\Lsep)^{\textstyle\t} \qquad\text{and}\qquad V_v\ulM \;:=\; T_v\ulM\otimes_\Av\Qv\,,
\]
where the superscript $(\mbox{...})^{\textstyle\t}$ denotes the $\t$-invariants. We call\/ $T_v\ulM$ (respectively $V_v\ulM$) the \emph{(rational) $v$-adic Tate module of $\ulM$}. 
This definition applies in particular if $\ulM$ is a pure Anderson motive.
\end{Definition}

\noindent {\it Remark.} 
Our functor $T_v$ is covariant. In the literature usually the $A_v$-dual of our $T_v \ulM$ is called the $v$-adic Tate module of $\ulM$. With that convention the Tate module functor is contravariant on $\tau$-modules but covariant on Drinfeld modules and Anderson's abelian $t$-modules~\cite{Anderson} (which both give rise to $\tau$-modules). Similarly the classical Tate module functor on abelian varieties is covariant. We chose our non standard convention here solely to avoid perpetual dualizations. This agrees also with the remark that abelian $\tau$-sheaves behave dually to abelian varieties; see~\ref{RemDualBehaviour}.

\smallskip

Next we make similar definitions for abelian $\tau$-sheaves.

\begin{Definition} \label{DefTateMod}
Let $\FF$ be an abelian $\tau$-sheaf and let $v\in C$ be a place different from the characteristic point $\chr$. We choose a finite closed subset $D\subset C$ as in section~\ref{SectRelation} with $v\notin D$ and $\infty\in D$ if $v\ne\infty$ and set
\[
T_v\FF \;:=\; T_v\ulM^{(D)}(\FF)\qquad\text{and}\qquad \VvFF \;:=\; V_v\ulM^{(D)}(\FF)\,.
\]
We call\/ $T_v\FF$ (respectively $\VvFF$) the \emph{(rational) $v$-adic Tate module associated to $\FF$}. It is independent of the particular choice of $D$, but if $v=\infty$ it depends on $k,l$ and $z$; see page~\pageref{EQ.Pi+Lambda}.
\end{Definition}

By \cite[Proposition~6.1]{TW}, $T_v\FF$ (and $V_v\FF$) are free $A_v$-modules (respectively $Q_v$-vector spaces) of rank $r$ for $v\ne\infty$ and $rl$ for $v=\infty$, which carry a continuous $G=\Gal(L^\sep/L)$-action.

\smallskip

Also the Tate modules $T_\infty\FF$ and $V_\infty\FF$ are always equipped with the automorphisms $\P$ and $\Lambda(\lambda)$ for $\lambda\in \Ff_{q^l}\cap L$ from (\ref{EQ.Pi+Lambda}). 
And if $\Ff_{q^l}\subset L$ we identify the algebra $\Delta_\infty$ from (\ref{EqDelta}) with a subalgebra of $\End_{Q_\infty}(V_\infty\FF)$ by mapping $\lambda\in\BF_{q^l}\subset\Delta_\infty$ to $\Lambda(\lambda)$.

\medskip

The following is evident from the definitions.

\begin{Proposition}\label{PropLS2b}
If $\FF$ is an abelian $\tau$-sheaf over $L$, respectively $\ulM$ a pure Anderson motive over $L$ and $v\in C$ (respectively $v\in \Spec A$) is a place of $Q$ different from the characteristic point $\chr$, then
\begin{eqnarray*}
T_v\FF=T_v(\ulM_v(\FF)) \quad&\text{and}&\quad V_v\FF=V_v(\ulM_v(\FF)) \quad\text{for }v\ne\infty,\\[2mm]
T_\infty\FF=T_\infty(\ulTM_\infty(\FF)) \quad&\text{and}&\quad V_\infty\FF=V_\infty(\ulTM_\infty(\FF)),\\[2mm]
\text{respectively}\quad T_v\ulM=T_v(\ulM_v(\ulM)) \quad&\text{and}&\quad V_v\ulM=V_v(\ulM_v(\ulM))\,.\qed
\end{eqnarray*}
\end{Proposition}

\bigskip

In order to prove the finiteness of $\Hom(\ulM,\ulM')$ we first need the following facts.

\begin{Proposition}\label{PropT.2}
Let $\FF$ and $\FF'$ be abelian $\tau$-sheaves over an arbitrary field $L$ and let $v$ be a place of $Q$ different from $\chr$.
\begin{suchthat}
\item
If $v\ne\infty$ then the natural map is injective
\[
\QHom(\FF,\FF')\otimes_Q Q_v\longto\Hom_{Q_v[G]}(V_v\FF,V_v\FF')\,.
\]
\item 
If $v=\infty$ then the natural map is injective
\[
\QHom(\FF,\FF')\otimes_Q Q_\infty \longto \Hom_{Q_\infty[\P,\Lambda,G]}(V_\infty\FF,V_\infty\FF')\,.
\]
\end{suchthat}
In particular $\QHom(\FF,\FF')$ is a $Q$-vector space of dimension $\le rr'$.
\end{Proposition}

\begin{proof}
1. Consider the morphisms
\[
\QHom(\FF,\FF')\otimes_Q Q_v\into\Hom_{Q_L}(\F_{0,\eta},\F'_{0,\eta})\otimes_{Q_L} (Q_v\otimes_\Fq L) \into\Hom_{Q_{v,L}}\bigl(\ulN_v(\FF),\ulN_v(\FF')\bigr)\,.
\]
obtained from \ref{PropAltDescrQHom} and the definition of $\ulN_v(\FF)$. Clearly the composition factors through $\Hom_{Q_{v,L}[\phi]}\bigl(\ulN_v(\FF),\ulN_v(\FF')\bigr)$. Since in both cases $\ulM_v(\FF)$ and $\ulM_v(\FF')$ are \'etale local shtukas, the isomorphy of the later $Q_v$-vector space with the one stated in the proposition follows by tensoring \ref{Prop2.13'} with $\otimes_{A_v}Q_v$.

\smallskip
\noindent 
2. We adapt the argument from 1 by replacing $\ulN_v$ and $Q_{v,L}$ by $\ulTN_\infty$ and $Q_{\infty,L}[\P,\Lambda]$ and the assertion follows as above.
\end{proof}

The following fact is well known and proved for instance in \cite[Proposition 1.2.4]{Taelman} even without the purity assumption. Nevertheless, we include a proof for the sake of completeness and to illustrate the use of the $\infty$-adic Tate module $V_\infty\FF$ arising from the big local shtuka $\ulTM_\infty(\FF)$.

\begin{Theorem}\label{ThmT.3}
Let $\ulM$ and $\ulM'$ be pure Anderson motives over an arbitrary field $L$. Then $\Hom(\ulM,\ulM')$ is a projective $A$-module of rank $\le rr'$.
\end{Theorem}

\begin{proof}
Since $M'$ is a locally free $A_L$-module, $H:=\Hom(\ulM,\ulM')$ is a torsion free, hence flat $A$-module, because all local rings of $A$ are principal ideal domains. We prove that $H$ is finitely generated by showing that $H$ is a discrete submodule of a finite dimensional $Q_\infty$-vector space. Let $\FF$ and $\FF'$ be abelian $\tau$-sheaves with $\ulM=\ulM(\FF)$ and $\ulM'=\ulM(\FF')$. Then Corollary~\ref{Cor2.9d} and Proposition~\ref{PropT.2} yield inclusions
\[
H\into H\otimes_A Q=\QHom(\FF,\FF')\into\Hom_{Q_\infty}(V_\infty\FF,V_\infty\FF')
\]
The later $Q_\infty$-vector space is finite dimensional and we claim that $H$ is a discrete $A$-lattice. Indeed choose $Q\otimes_\Fq L$-bases $\ul m=(m_1,\ldots,m_{rl})$ of $\bigoplus_{i=0}^{l-1}\F_{i,\eta}$ with $m_j\in M$ and $\ul m'=(m'_1,\ldots,m'_{r'l'})$ of $\bigoplus_{i=0}^{l'-1}\F'_{i,\eta}$ such that $\bigoplus_{i=0}^{l'-1}M'\subset \bigoplus_{j=1}^{r'l'}A_L m'_j$. With respect to these bases every element of $H$ corresponds to a matrix in $M_{r'l'\times rl}(A_L)$. Now choose $Q_\infty$-bases $\ul n$ of $V_\infty\FF$ and $\ul n'$ of $V_\infty\FF'$ and denote the base change matrix from $\ul m$ to $\ul n$ by $B\in\GL_{rl}(Q_{\infty,L^\sep})$, and the base change matrix from $\ul m'$ to $\ul n'$ by $B'\in\GL_{r'l'}(Q_{\infty,L^\sep})$. Then $H$ is contained in
\[
M_{r'l'\times rl}(Q_\infty) \cap B'\cdot M_{r'l'\times rl}(A_L)\cdot B^{-1}
\]
which is discrete in $M_{r'l'\times rl}(Q_\infty)$. This proves that $H$ is a projective $A$-module. The estimate on the rank of $H$ follows from \ref{Cor2.9d} and \ref{PropT.2}.
\end{proof}

\begin{Corollary}\label{CorT.4}
The minimal polynomial of every endomorphism of a pure Anderson motive $\ulM$ lies in $A[x]$. \qed
\end{Corollary}

\begin{Proposition}\label{PropT.1}
Let $\ulM$ and $\ulM'$ be pure Anderson motives over an arbitrary field $L$ and let $v\in\Spec A$ be a maximal ideal different from $\chr$. Then the natural map
\[
\Hom(\ulM,\ulM')\otimes_A A_v \longto \Hom_{A_v[G]}(T_v\ulM,T_v\ulM')
\]
is injective with torsion free cokernel.
\end{Proposition}

\begin{proof}
Consider the morphisms
\[
\Hom(\ulM,\ulM')\otimes_A A_v \into \Hom_{A_L}(M,M')\otimes_{A_L}(A_v\otimes_\Fq L)\into \Hom_{A_L}(M,M')\otimes_{A_L}A_{v,L}
\]
which are injective because $A_v$ is flat over $A$, respectively because $\Hom_{A_L}(M,M')$ is flat over $A_L$. Again the composite morphism factors through
\[
\Hom_{A_{v,L}[\phi]}\bigl(\ulM_v(\ulM),\ulM_v(\ulM')\bigr)\;=\;\Hom_{A_v[G]}(T_v\ulM,T_v\ulM')
\]
(use \ref{Prop2.13'}). To prove that the cokernel is torsion free let $f_v$ be an element of the cokernel which is torsion. Since a power of $v$ is principal we may assume that $g_v=a f_v\in\Hom(\ulM,\ulM')\otimes_A A_v$ for an $a\in A$ with $(a)=v^m$ for some $m\in\N$. Fix a positive integer $n$. There exists a $g\in\Hom(\ulM,\ulM')$ with $g\equiv g_v\mod v^{n+m}$. In particular $a$ divides $g$ in $\Hom_{A_v[G]}(T_v\ulM,T_v\ulM')$. Since
\[
\Bigl((T_v\ulM'/a\cdot T_v\ulM')\otimes_{A/(a)} A_{L^\sep}/(a)\Bigr)^G\;\cong\;\ulM'/a\ulM'
\]
(compare \ref{Prop2.13'}) we see that $g$ maps $\ulM$ into $a\ulM'$. Thus $g$ factors, $g=af$ with $f\in\Hom(\ulM,\ulM')$ and $f\equiv f_v\mod v^n$. As $n$ was arbitrary and $\Hom(\ulM,\ulM')$ is a finitely generated $A$-module the proposition follows.
\end{proof}

If $L$ is finitely generated, Proposition~\ref{PropT.1} was strengthened by Taguchi~\cite{Taguchi95b} and Tamagawa~\cite[\S 2]{Tam} to the following analog of the Tate conjecture for abelian varieties.

\begin{Theorem}[Tate conjecture for $\t$-modules]\label{TATE-CONJECTURE-MODULES}
Let\/ $\ulM$ and $\ulM'$ be two $\t$-modules on $\TA$ over a finitely generated field $L$ and let $G:=\Gal(\Lsep/L)$. Let\/ $v\in\Spec\TA$ such that the support of\/ $\coker\t'$ does not meet\/ $v$. For instance $\ulM$ and $\ulM'$ could be pure Anderson motives, $\TA=A$, and $v\ne\chr$. Then the Tate conjecture holds:
\[
\Hom(\ulM,\ulM')\otimes_{\TA}\Av \;\cong\; \Hom_\AvG(T_v\ulM,T_v\ulM')\,.\qed
\]
\end{Theorem}

\noindent {\it Remark.}
Note one interesting consequence of this result. If $\ulM$ and $\ulM'$ are pure Anderson motives of different weights over a finitely generated field then $\Hom_{A_v[G]}(T_v\ulM,T_v\ulM')=(0)$.

\begin{Theorem}[Tate conjecture for abelian $\tau$-sheaves]\label{TATE-CONJECTURE}
Let $\FF$ and $\FF'$ be abelian $\tau$-sheaves over a finitely generated field $L$ and let $G:=\Gal(\Lsep/L)$. Let\/ $v\in C$ be a place different from the characteristic point $\chr$. 
\begin{suchthat}
\item
If $v\ne\infty$ assume the characteristic $\chr$ is different from $\infty$ or $\weight(\FF)=\weight(\FF')$. Then 
\[
\QHom(\FF,\FF')\otimes_Q\Qv\;\cong\;\Hom_\QvG(\VvFF,\VvFF')\,.
\]
\item
If $v=\infty$ choose an integer $l$ which satisfies condition 2 of \ref{Def1.1} for both $\FF$ and $\FF'$ and assume $\Ff_{q^l}\subset L$. Then 
\[
\QHom(\FF,\FF')\otimes_Q Q_\infty\;\cong\;\Hom_{\Delta_\infty[G]}(V_\infty\FF,V_\infty\FF')\,.
\]
\end{suchthat}
\end{Theorem}

\begin{proof}
1. Set $\ulM:=\ulM(\FF)$ and $\ulM':=\ulM(\FF')$. By \ref{CONNECTION} and \ref{TATE-CONJECTURE-MODULES}, we have
\[
\QHom(\FF,\FF')\otimes_Q Q_v\;\cong\;\Hom(\ulM,\ulM')\otimes_AQ_v \;\cong\;
\Hom_\QvG(V_v\ulM,V_v\ulM')\,.
\]

\smallskip\noindent
2. Let $D\subset C$ be a finite closed subscheme as in Section~\ref{SectRelation} with $\chr,\infty\notin D$ and set $\ulM:=\ulM^{(D)}(\FF)$ and $\ulM':=\ulM^{(D)}(\FF')$. By \ref{CONNECTION} and \ref{TATE-CONJECTURE-MODULES}, we have
\[
\QHom(\FF,\FF')\otimes_Q Q_\infty\;\cong\;\Hom_{\P,\Lambda}(\ulM,\ulM')\otimes_{\TA}Q_\infty \;\cong\;
\Hom_{\Delta_\infty[G]}(V_\infty \ulM,V_\infty \ulM')\,.
\]
Here the last isomorphism comes from the fact that the commutation with $\P$ and $\Lambda(\lambda)$ are linear conditions on $\Hom(\ulM,\ulM')$ and $\Hom(\ulM,\ulM')\otimes_{\TA}Q_\infty\cong\Hom_{Q_\infty[G]}(V_\infty \ulM,V_\infty \ulM')$ thus cutting out isomorphic subspaces.
\end{proof}

As a direct consequence of the theorem together with Proposition~\ref{Prop2.13'} we obtain:

\begin{Corollary}\label{Cor2.17'}
\begin{suchthat}
\item 
Let $\ulM$ and $\ul M'$ be pure Anderson motives over a finitely generated field and let $v\in\Spec A$ be a maximal ideal different from the characteristic point $\chr$, then 
\[
\Hom(\ulM,\ulM')\otimes_A A_v\;\cong\;\Hom_{A_{v,L}[\phi]}\bigl(\ulM_v(\ulM),\ulM_v(\ulM')\bigr)\,.
\]
\item
Let $\FF$ and $\FF'$ be abelian $\tau$-sheaves over a finitely generated field $L$ and let $v$ be a place of $Q$ different from $\chr$ and $\infty$. If $\chr=\infty$ assume $\weight(\FF)=\weight(\FF')$. Then
\[
\QHom(\FF,\FF')\otimes_Q Q_v\;\cong\;\Hom_{Q_{v,L}[\phi]}\bigl(\ulN_v(\FF),\ulN_v(\FF')\bigr)\,.\qed
\]
\end{suchthat}
\end{Corollary}

\bigskip

Finally we establish the relation between Tate modules and isogenies.

\begin{Proposition}\label{Prop2.7b}
\begin{suchthat}
\item 
Let $f:\ulM'\to\ulM$ be an isogeny between pure Anderson motives then $T_vf(T_v\ulM')$ is a $G$-stable lattice in $V_v\ulM$ contained in $T_v\ulM$.
\item 
Conversely if $\ulM$ is a pure Anderson motive and $\Lambda_v$ is a $G$-stable lattice in $V_v\ulM$ contained in $T_v\ulM$, then there exists a pure Anderson motive $\ulM'$ and a separable isogeny $f:\ulM'\to \ulM$ with $T_vf(T_v\ulM')=\Lambda_v$.
\end{suchthat}
\end{Proposition}

\begin{proof}
1 follows from the $G$-invariance of $f$.

\noindent
2. Consider the $A_{v,L^\sep}[G,\phi]$-module $\Lambda_v\otimes_{A_v}A_{v,L^\sep}$. The action of $\phi$ through $\s$ on the right factor and of $G$ through both factors commute. This module is a submodule of 
\[
T_v\ulM\otimes_{A_v}A_{v,L^\sep}\;=\; \ulM_v(\ulM)\otimes_{A_{v,L}}A_{v,L^\sep}
\]
(see Proposition~\ref{Prop2.13'}), and contains $a\cdot \ulM_v(\ulM)\otimes_{A_{v,L}}A_{v,L^\sep}$ for a suitable $a\in A$. Therefore the $G$-invariants $(\Lambda_v\otimes_{A_v}A_{v,L^\sep})^G$ form an \'etale local $\sigma$-subshtuka $\ulHM{}'$ of $\ulM_v(\ulM)$ of the same rank. Now by Proposition~\ref{Prop1.5b} the kernel of the surjection $\ulM\shortonto\ulM_v(\ulM)/\ulHM{}'$ is a pure Anderson motive $\ulM'$ together with a separable isogeny $f:\ulM'\to\ulM$. Clearly $\ulM_v f(\ulM_v\ulM')=\ulHM{}'$ and $T_vf(T_v\ulM')=\Lambda_v$.
\end{proof}

\begin{Proposition}\label{FACTORSHEAF-FACTORSPACE}
Let\/ $\FF'$ be an abelian quotient $\tau$-sheaf of\/ $\FF$. Then $\VvFF'$ is a $G$-quotient space of\/ $\VvFF$. The same holds if $\ulM'$ is a  quotient motive of a pure Anderson motive $\ulM$.
\end{Proposition}

\begin{proof}
Let $f\in\Hom(\FF,\FF')$ be surjective and let $\ulHM$ and $\ulHM'$ be the (big, if $v=\infty$) local $\sigma$-shtuka of $\FF$, respectively $\FF'$, at $v$. Then the induced morphism $\ulM_v(f)\in\Hom(\ulHM,\ulHM')$ is surjective and $\ulHM'':=\ker \ulM_v(f)$ is also a local $\sigma$-shtuka at $v$. We get an exact sequence of local $\sigma$-shtukas which we tensor with $A_{v,L^\sep}$ yielding
\[
\bigexact{0}{}{\ulHM''\otimes_{A_{v,L}}A_{v,L^\sep}}{}{\ulHM\otimes_{A_{v,L}}A_{v,L^\sep}}{\ulM_v(f)}{\ulHM'\otimes_{A_{v,L}}A_{v,L^\sep}}{}{0\,.}
\]
The Tate module functor is left exact, because considering the morphism of $A_{v,L^\sep}$-modules
\[
1-\t:\es\ulHM\otimes_{A_{v,L}}A_{v,L^\sep} \longto \ulHM\otimes_{A_{v,L}}A_{v,L^\sep}
\]
we have by definition $T_v\ulHM = \ker(1-\t)$, and the desired left exactness follows from the snake lemma. After tensoring with $\otimes_\Av\Qv$ we get
\[
\bigexact{0}{}{V_v\ulHM''}{}{V_v\ulHM}{V_v f}{V_v\ulHM'\;.}{}{}
\]
Counting the dimensions of these $\Qv$-vector spaces, we finally also get right exactness, as desired.
\end{proof}

\forget{
\begin{proof}
Let $f\in\Hom(\FF,\FF')$ be surjective, let $\ulM:=\ulM^{(D)}(\FF)$ \? and let $\ulM':=\ulM^{(D)}(\FF')$. Then the induced morphism $\ulM^{(D)}(f)\in\Hom(\ulM,\ulM')$ is surjective and $\ulM'':=\ker \ulM^{(D)}(f)$ is a $\t$-module on $A$. Thus we get the exact sequence
\[
\bigexact{0}{}{M''}{}{M}{\ulM^{(D)}(f)}{M'}{}{0\,.}
\]
The exactness being preserved since $M'$ is locally free, we consider the following diagram
\[
\xymatrix{
0 \ar[r] &
M''/v^nM'' \otimes_L\Lsep \ar[r] &
M/v^nM \otimes_L\Lsep \ar[r] &
M'/v^nM' \otimes_L\Lsep \ar[r] &
0 \\
0 \ar@{-->}[r] &
(M''/v^nM''\otimes_L\Lsep)^{\textstyle\t} \ar@{-->}[r]\ar[u] &
(M/v^nM\otimes_L\Lsep)^{\textstyle\t} \ar@{-->}[r]\ar[u] &
(M'/v^nM'\otimes_L\Lsep)^{\textstyle\t} \ar[u] \;.&
\\
}
\]
The $\t$-invariant functor is left exact, because considering the morphism of $A/{v^n}$-modules
\[
1-\t:\,(M/v^nM)\otimes_L\Lsep\rightarrow(M/v^nM)\otimes_L\Lsep
\]
we have by definition $((M/v^nM)\otimes_L\Lsep)^{\textstyle\t} = \ker\, 1-\t$, and the desired left exactness follows from the snake lemma. Since the projective limit preserves left exactness as well, we get after tensoring with $\otimes_\Av\Qv$
\[
\bigexact{0}{}{V_v\ulM''}{}{V_v\ulM}{V_v\ulM^{(D)}(f)}{V_v\ulM'\,.}{}{}
\]
Counting the dimensions of these $\Qv$-vector spaces, we finally also get right exactness, as desired.
\end{proof}

}

% -----------------------------------------------------------------------------

\bigskip

\subsection{The Frobenius Endomorphism}

Suppose that the characteristic is finite, that is, the characteristic point $\chr$ is a closed point of $C$ with finite residue field $\Ff_\chr$, and the map $c:\Spec L\to C$ factors through the finite field $\chr=\Spec\Ff_\chr$.

\begin{Definition}[$s$-Frobenius on abelian $\tau$-sheaves] \label{DefFrob}
Let $\FF$ be an abelian $\tau$-sheaf with finite characteristic point $\chr=\Spec \Ff_\chr$ and let $s=q^e$ be a power of the cardinality of $\Ff_\chr$. We define the \emph{$s$-Frobenius on $\FF$} by
\[
\pi \;:=\; (\pi_i):\, (\sigma^\ast)^e\FF\rightarrow\FF\II{e}, \quad 
\pi_i \;:=\; \t_{i+e-1}\circ\cdots\circ(\s)^{e-1}\t_{i}:\, (\sigma^\ast)^e\F_{i}\rightarrow\F_{i+e}\;.
\]
Clearly $\pi$ is an isogeny. Observe that $\Ff_\chr\subset\Fs$ implies that $(\sigma^\ast)^e\FF$ has the same characteristic as $\FF$. 
\end{Definition}

Similarly if $\chr\in\Spec A$ is a closed point we define

\begin{Definition}[$s$-Frobenius on pure Anderson motives]\label{Def2.19b}
Let $\ulM$ be a pure Anderson motive with finite characteristic point $\chr=\Spec \Ff_\chr$ and let $s=q^e$ be a power of the cardinality of $\Ff_\chr$. We define the \emph{$s$-Frobenius isogeny on $\ulM$} by
\[
\pi\;:=\;\t\circ\ldots\circ (\s)^{e-1}\t:\,(\sigma^\ast)^e\ulM\to\ulM\,.
\]
\end{Definition}

\begin{Remark}
Classically for (abelian) varieties $X$ over a field $K$ of characteristic $p$ one defines the Frobenius morphism $X\to\phi^\ast X$ where $\phi$ is the $p$-Frobenius on $K$. There $p$ equals the cardinality of the ``characteristic field'' $\im(\Z\to K)=\Ff_p$. 
In view of the dual behavior of abelian $\tau$-sheaves and pure Anderson motives our definition is a perfect analog since here we consider the $s$-Frobenius for $s$ being the cardinality of (a power of) the ``characteristic field'' $\im(c^\ast:A\to L)=\BF_\chr$.
\end{Remark}

Now we suppose $L=\Fs$ to be a finite field with $s=q^e$ $(e\in\N$). Let $\FS$ denote a fixed algebraic closure of $\Fs$ and set $G=\GalFSFs$. It is topologically generated by ${\rm Frob}_s:x\mapsto x^s$. The following results for the Frobenius endomorphism of $\t$-modules can be found in Taguchi and Wan~\cite[\S 6]{TW}.

\begin{Definition} \label{Def2.5}
Let $\Spec\TA\subset C$ be an affine open subscheme and let $\ulM$ be a\/ $\t$-module on $\TA$ over $\Fs$. Since $\sigma^e=\id_\CFs$, we define the \emph{$s$-Frobenius on $\ulM$} by
\[
\pi \;:=\; \t^e \;:=\; \t\circ\s\t\circ\cdots\circ(\s)^{e-1}\t :\; (\s)^e M = M\rightarrow M\,.
\]
\end{Definition}

\begin{Proposition}\label{ABSOLUTE-GALOIS}
Let $\ulM$ be a $\t$-module on $\TA$ over $\Fs$ of rank $r$ and let\/ $v\in\Spec\TA$ such that the support of\/ $\coker\t$ does not meet\/ $v$. For example $\ulM$ can be a pure Anderson motive, $\TA=A$, and $v\ne\chr$.
\begin{suchthat}
\item The generator ${\rm Frob}_s$ of $G$ acts on $T_v\ulM$ like\/ $(T_v\pi)^{-1}$.
\item Let\/ $\Psi:\,\AvG\rightarrow\End_\Av(T_v\ulM)$ denote the continuous morphism of $\Av$-modules which is induced by the action of\/ $G$ on $T_v\ulM$. Then\/ $\im\Psi = \Av\II{T_v\pi}$.
\end{suchthat}
\end{Proposition}

\begin{proof}
1 was proved in \cite[Ch. 6]{TW} and 2 follows from the continuity of $\Psi$.
\end{proof}

\noindent {\it Remark.}
The inversion of $T_v\pi$ in the first statement results from the dual definition of our Tate module.

\begin{Proposition}\label{PI-IS-QISOG}
Let\/ $\FF$ be an abelian $\tau$-sheaf over $L=\Fs$ with $s=q^e$ and let\/ $\pi$ be its $s$-Frobenius. Then $(\sigma^\ast)^e\FF=\FF$. Let\/ $v\in C$ be a place different from $\infty$ and from the characteristic point $\chr$. 
\begin{suchthat}
\item The $s$-Frobenius $\pi$ can be considered as a quasi-isogeny of\/ $\FF$.
\item The generator ${\rm Frob}_s$ of $G$ acts on $T_v\FF$ like\/ $(T_v\pi)^{-1}$.
\item The image of the continuous morphism of $\Qv$-vector spaces $\QvG\rightarrow\End_\Qv(V_v\FF)$ is $\Qv\II{V_v\pi}$.
\item $\ulM^{(D)}(\pi)$ coincides with the $s$-Frobenius on $\ulM^{(D)}(\FF)$ from definition~\ref{Def2.5}.
\end{suchthat}
\end{Proposition}

\begin{proof}
1. Due to the periodicity condition, we have $\FF\II{e}\subset\FF(nk\cdot\infty)$ for a sufficiently large $n\in\N$, since $\F_{i+e}\subset\F_{i+nl}=\F_i(nk\cdot\infty)$ for $e\le nl$. Thus $\pi\in\Hom(\FF,\FF(nk\cdot\infty))$, and therefore $\pi\in\QEnd(\FF)$. By \ref{PROP.1.42A}, we have $\pi\in\QIsog(\FF)$.

\noindent 
4. This follows from the definition of $\pi$ and the commutation of the $\P$'s and the $\t$'s and then 2 and 3 follow from \ref{ABSOLUTE-GALOIS}.
\end{proof}

\vspace{1cm}

% =============================================================================

\section{Applications to Pure Anderson Motives over Finite Fields}

\bigskip

\subsection{The Poincar\'e-Weil Theorem}

In this section we study the analog for pure Anderson motives and abelian $\tau$-sheaves of the Poincar\'e-Weil theorem. Originally, this theorem states that every abelian variety is semisimple, that is, isogenous to a product of simple abelian varieties, see \cite[Corollary of Theorem~II.1/6]{La}. Unfortunately, we cannot expect a full analog of this statement for abelian $\tau$-sheaves or pure Anderson motives as our next example illustrates. On the positive side we show that every abelian $\tau$-sheaf or pure Anderson motive over a finite field becomes semisimple after a finite base field extension.

\begin{Example} \label{Ex3.1}%
Let $C=\PP^1_\Fq$, $C\setminus\{\infty\}=\Spec\Fq\II{t}$ and $\zeta:=c^\ast(1/t)\in\Fq^\times$. We construct an abelian $\tau$-sheaf $\FF$ over $L=\Fq$ with $r=d=2$ which is not semisimple. Let 
\[\textstyle
\Delta=\tmatr{1}{0}{0}{1} + \tmatr{\alpha}{\gamma}{\beta}{\delta}\cdot t
\]
with $\alpha,\beta,\gamma,\delta\in\Fq$.
To obtain characteristic $c$ we need $\det\Delta=(1-\zeta t)^2$, and thus we require the conditions $\alpha+\delta=-2\zeta$ and $\alpha\delta-\beta\gamma=\zeta^2$. We set $\F_i:=\O_{C_L}(i\cdot\infty)^{\oplus 2}$, we let $\P_i$ be the natural inclusion, and  we let $\t_i:=\Delta$. Then $\FF$ is an abelian $\tau$-sheaf with $r=d=2$ and $k=l=1$. The associated pure Anderson motive is $\ulM=(L[t]^{\oplus 2},\Delta)$.

We see that $\FF$ is not simple. If $\Delta=\smatr{1-\zeta t}{0}{0}{1-\zeta t}$ then $\FF$ is semisimple as a direct sum of two simple abelian $\tau$-sheaves. Otherwise, if $\Delta\ne\smatr{1-\zeta t}{0}{0}{1-\zeta t}$ which is the case for example if $\beta\ne 0$, consider
\[\textstyle 
\tilde{\Delta} \,:=\, \tmatr{\beta}{\delta+\zeta}{0}{1}^{-1} \!\!\cdot\Delta\cdot \s\tmatr{\beta}{\delta+\zeta}{0}{1} \,=\, \tmatr{1-\zeta t}{0}{t}{1-\zeta t}
\]
and the abelian $\tau$-sheaf with $\wt\F_i=\O_{C_L}(i\cdot\infty)^{\oplus 2}$ and $\tilde\t_i=\tilde\Delta$ which is isomorphic to $\FF$.
There is an exact sequence 
\[
\begin{array}{c}
\bigexact{0}{}{\FF'}{\varphi}{\TFF}{\psi}{\FF''}{}{0} \\[1ex]
\scriptstyle \t'\,=\,1-\zeta t \qquad\quad\Tt\quad\qquad \t''\,=\,1-\zeta t
\end{array}
\]
with $\varphi:\, 1\mapsto\tvect{1}{0}$ and $\psi:\,\tvect{x}{y}\mapsto y$ where $\FF'=\FF''$ is the abelian $\tau$-sheaf with $\F'_i=\O_{C_L}(i\cdot\infty)$ and $\t'_i=1-\zeta t$. If $\TFF$ were semisimple, then there would be a quasi-morphism $\omega:\,\FF''\rightarrow\TFF$ with $\psi\circ\omega=\id_{\FF''}$, hence $\omega:\, y\mapsto\tvect{e}{1}\cdot y$ for some $e\in \Fq(t)$. Thus, a necessary condition for the semisimplicity of $\FF$ is
\[\textstyle
(1-\zeta t)\cdot\s(y)\cdot \tvect{e}{1}
\;=\;
\tmatr{1-\zeta t}{0}{t}{1-\zeta t}\cdot\tvect{\s(e)}{1}\cdot\s(y)
\]
which is equivalent to the condition
\[
e-\s(e)\;=\;\frac{t}{1-\zeta t}\enspace.
\]
But this cannot be true since $e-\s(e)=0$, thus $\FF$ is \emph{not} semisimple. However, this last formula is satisfied if $e=\lambda\cdot\frac{t}{1-\zeta t}$ for $\lambda\in\Ff_{q^q}$ with $\lambda^q-\lambda=-1$. That means that after field extension $\Fq(\lambda)\,/\,\Fq$ we get $\FF\cong{\FF'}\mbox{\,}^{\oplus 2}$  and we have $\QEnd(\FF)=M_2(\QEnd(\FF'))=M_2(Q)$. Note that this phenomenon generally appears, and we will state and prove it in Theorem \ref{Thm3.8b}.
\end{Example}

From now on we fix a place $v\in\Spec A$ which is different from the characteristic point $\chr$ of $c$. For a morphism $f\in\QHom(\FF,\FF')$ between two abelian $\tau$-sheaves $\FF$ and $\FF'$ we denote its image $V_vf\in\Hom_\QvG(\VvFF,\VvFF')$ just by $f_v$. If $\FF$ is defined over $\Fs$ this applies in particular to the $s$-Frobenius endomorphism $\pi$ of $\FF$ (Definition~\ref{DefFrob}).

Let $\FF$ be an abelian $\tau$-sheaf over the finite field $L=\Fs$. We set
\begin{equation}\label{Eq3.1}
\begin{array}{l@{\;}c@{\;}l@{\qquad\qquad}l@{\;}c@{\;}l}
E &:=& \QEnd(\FF) \ni\pi & E_v &:=& \End_\QvG(\VvFF) \ni\pi_v \\[0.5ex]
F &:=& Q\II{\pi} \subset E & F_v &:=& \im(\QvG\rightarrow \End_\Qv(\VvFF)) \subset E_v \\
\end{array}
\end{equation}
with $\QvG\rightarrow \End_\Qv(\VvFF)$ induced by the action of $G$ on $\VvFF$. Clearly, we have $F\subset E$ and $F_v\subset E_v$ by Proposition~\ref{PI-IS-QISOG}/3. By \ref{PropT.2}, we know that $\dim_Q E<\infty$. Thus $\pi$ is algebraic over $Q$. We denote its minimal polynomial by $\mu_\pi\in Q[x]$, and the characteristic polynomial of the endomorphism $\pi_v$ of $V_v\FF$ by $\chi_v\in Q_v[x]$. If $\chr\ne\infty$ Corollary~\ref{CorT.4} shows that $\pi$ is integral over $A$, $\mu_\pi\in A[x]$. The zeroes of $\pi$ in $\Spec A[\pi]$ all lie above $\chr$ because $\pi$ is an isomorphism locally at all $v\ne\chr$; compare with Remark~\ref{Rem1.26'}.

Due to the Tate conjecture, our situation can be represented by the following diagram where we want to fit the missing bottom right arrow with an isomorphism.
\[
\xymatrix{
\,E^{\,\!}_{\,\!} \ar[r] & 
\,E\xQv\, \ar[r]^<<{\sim} &
\,E_v \\
\,F^{\,\!}_{\,\!} \ar[r]\ar[u] & 
\,F\xQv\, \ar@{-->}[r]^<<{\sim}\ar[u] &
\,F_v \,.\!\!\ar[u] \\
}
\]

\begin{Lemma}\label{Lemma3.4}
The natural morphism between $F\xQv$ and\/ $F_v$ is an isomorphism.
\end{Lemma}

\begin{proof}
Consider the isomorphism $\psi:\,E\xQv\cong E_v\subset\End_\Qv(\VvFF)$ and set $\varphi:=\psi|_{F\xQv}$. Then $\varphi$ is injective and maps into $F_v$. Since $\im\varphi=\Qv\II{\pi_v}$, the surjectivity follows from Proposition~\ref{PI-IS-QISOG}.
\end{proof}

To evaluate the dimension of $E$ we need the following notation.

\begin{Definition}\label{Def3.3}
Let\/ $K$ be a field. Let\/ $f,g\in K\II{x}$ be two polynomials and let
\[
f\, = \!\!\!\prod\limits_{\genfrac{}{}{0pt}{}{\scriptstyle\mu\in K\II{x}}{\mbox{\rm\scriptsize irred.}}}\! \mu^{m(\mu)}
,\qquad
g\, = \!\!\!\prod\limits_{\genfrac{}{}{0pt}{}{\scriptstyle\mu\in K\II{x}}{\mbox{\rm\scriptsize irred.}}}\! \mu^{n(\mu)}
\]
be their respective factorizations in powers of irreducible polynomials. Then we define the integer
\[
r_K(f,g)\; := \!\!\!\prod\limits_{\genfrac{}{}{0pt}{}{\scriptstyle\mu\in K\II{x}}{\mbox{\rm\scriptsize irred.}}}\! m(\mu)\cdot n(\mu)\cdot\deg\mu \;.
\]
\end{Definition}

\noindent {\it Remark.}
In contrast to characteristic zero, we have for $\charakt(K)\ne0$ in general different values of the integer $r_K$ for different ground fields $K$. Namely, if $K\subset L$ then $r_K(f,g)\le r_L(f,g)$ with equality if and only if all irreducible $\mu\in K\II{x}$ which are contained both in $f$ and in $g$ have no multiple factors in $L\II{x}$. This is satisfied for example if the greatest common divisor of $f$ and $g$ has only separable irreducible factors, or if $L$ is separable over $K$. See \ref{LAST-EXAMPLE} below for an example where $r_K(f,g)<r_L(f,g)$.

\medskip

Before we discuss semisimplicity criteria in \ref{PROP.EQ} -- \ref{Cor3.8c}, let us compute the dimension of $\QHom(\FF,\FF')$.

\begin{Lemma}\label{PROP.15}\label{PI-END-SEMISIMPLE}
Let $v$ be a place of $Q$ different from $\chr$ and $\infty$. Let $\FF$ and $\FF'$ be abelian $\tau$-sheaves over $\Fs$ and assume that $\pi_v$ and $\pi'_v$ are semisimple. Factor their characteristic polynomials $\chi_v \,=\, \mu_1^{m_1}\cdot\ldots\cdot\mu_n^{m_n}$ and $\chi'_v\,=\, \mu_1^{m'_1}\cdot\ldots\cdot\mu_n^{m'_n}$ with distinct monic irreducible polynomials $\mu_1,\dots,\mu_n\in\Qv\II{x}$ and $m_i,m'_i\in\N_0$. Then
\begin{suchthat}
\item
$\DS\Hom_{Q_v[G]}(V_v\FF,V_v\FF')\,\cong\,\bigoplus_{i=1}^nM_{m_i'\times m_i}\bigl(Q_v[x]/(\mu_i)\bigr)$ as $Q_v$-vector spaces,
\item 
$\DS\End_{Q_v[G]}(V_v\FF)\,\cong\,\bigoplus_{i=1}^nM_{m_i}\bigl(Q_v[x]/(\mu_i)\bigr)$ as $Q_v$-algebras, and
\item
$\dim_\Qv\Hom_\QvG(\VvFF,\VvFF') \,=\, r_\Qv(\chi_v,\chi'_v)$\,.
\end{suchthat}
\end{Lemma}

\begin{proof}
Clearly 2 and 3 are consequences of 1 which we now prove. Since $\pi_v$ and $\pi'_v$ are semisimple, we have the following decomposition of $Q_v[G]$-modules
\[
\VvFF \,\cong\, \bigoplus_{i=1}^n\,(\Qv\II{x}/(\mu_i))^{\oplus m_i}, \qquad
\VvFF'\,\cong\, \bigoplus_{i=1}^n\,(\Qv\II{x}/(\mu_i))^{\oplus m'_i}
\]
where $\Qv\II{x}/(\mu_i)=:K_i$ are fields. Obviously, we only have non-zero $Q_v[G]$-morphisms $K_i\rightarrow K_j$ if $i=j$, since otherwise $\mu_i(\pi)\ne0$ in $K_j$. Since $\pi_v$ operates on $K_i^{\oplus m_i}$ as multiplication by the scalar $x$, the lemma follows.
\end{proof}

\begin{Theorem}\label{Thm3.5a}
Let $\FF$ and $\FF'$ be abelian $\tau$-sheaves of the same weight over $\BF_s$ and assume that the endomorphisms $\pi_v$ and $\pi'_v$ of $V_v\FF$ and $V_v\FF'$ are semisimple at a place $v\ne\chr,\infty$ of $Q$. Let $\chi_v$ and $\chi'_v$ be their characteristic polynomials. Then
\[
\dim_Q\QHom(\FF,\FF')\;=\;r_{Q_v}(\chi_v,\chi'_v).
\]
\end{Theorem}

\begin{proof}
This follows from the lemma and the Tate conjecture, Theorem~\ref{TATE-CONJECTURE}.
\end{proof}

\begin{Corollary}\label{Cor3.11b}
Let $\FF$ be an abelian $\tau$-sheaf of rank $r$ over $\Fs$ with Frobenius endomorphism $\pi$ and let $\mu_\pi$ be the minimal polynomial of $\pi$. Assume that $F=Q[x]/(\mu_\pi)$ is a field and set $h:=[F:Q]=\deg\mu_\pi$. Then
\begin{suchthat}
\item
$h|r$ and $\dim_Q \QEnd(\FF) =\frac{r^2}{h}$ and $\dim_F\QEnd(\FF)=\frac{r^2}{h^2}$.
\item 
For any place $v$ of $Q$ different from $\chr$ and $\infty$ we have $\QEnd(\FF)\otimes_Q Q_v\cong M_{r/h}(F\otimes_Q Q_v)$ and $\chi_v=(\mu_\pi)^{r/h}$ independent of $v$.
\end{suchthat}
\end{Corollary}

\begin{proof}
Since $F$ is a field, $\pi_v$ is semisimple by \ref{PROP.EQ} below. So general facts of linear algebra imply that $\mu_\pi=\mu_1\cdot\ldots\cdot\mu_n$ with pairwise different irreducible monic polynomials $\mu_i\in Q_v[x]$ and $\chi_v=\mu_1^{m_1}\cdot\ldots\cdot\mu_n^{m_n}$ with $m_i\ge1$. We set $K_i=Q_v[x]/(\mu_i)$ and use the notation from (\ref{Eq3.1}). By Lemma~\ref{PROP.15} the semisimple $Q_v$-algebra $E_v$ decomposes $E_v\cong \bigoplus_{i=1}^n E_i$ into the simple constituents $E_i=M_{m_i}(K_i)$. By \cite[Th\'eor\`eme~5.3/1 and Proposition~5.4/12]{Bou}, $E_i=E_v\cdot e_i$ where $e_i$ are the central idempotents with $K_i=F_v\cdot e_i$. Thus there are epimorphisms of $K_i$-vector spaces
\[
\QEnd(\FF)\otimes_F K_i \,=\,E_v\otimes_{F_v}K_i \,\onto\, E_i\,.
\]
This shows that $m_i^2\le \dim_F E$. So by Lemma~\ref{PROP.15}
\begin{eqnarray*}
[F:Q]\cdot\dim_F E\es=\es\dim_{Q_v}E_v&=&\sum_{i=1}^n m_i^2\deg\mu_i\es\le \\[-2mm]
\le\es\dim_F E\cdot\sum_{i=1}^n\deg\mu_i &=& \dim_F E\cdot\deg\mu_\pi\es =\es [F:Q]\cdot\dim_F E\,.
\end{eqnarray*}
Therefore $m_i^2=\dim_F E$ for all $i$. Since $r=\deg\chi_v=\sum_i m_i\deg\mu_i=\sqrt{\dim_F E}\cdot[F:Q]$. We find $r=m_ih$ and $\dim_FE=\frac{r^2}{h^2}$, proving 1. For 2 we use that 
\[
E_v\es\cong\es\bigoplus_i\, M_{r/h}\bigl(Q_v[x]/(\mu_i)\bigr)\es=\es M_{r/h}\bigl(\bigoplus_i Q_v[x]/(\mu_i)\bigr)\es=\es M_{r/h}\bigl(Q_v[x]/(\mu_\pi)\bigr)\,.
\]

\end{proof}

\medskip

Next we investigate when $\pi_v$ is semisimple.

\begin{Remark}\label{SEPARABLE}
Notice that the completion $\Qv$ is separable over $Q$. Namely, in terms of \cite[IV.7.8.1--3]{EGA}, we can state that $\O_{C,v}$ is an excellent ring. Thus the formal fibers of $\widehat{\O}_{C,v} \longrightarrow {\O_{C,v}}$ and therefore $Q_v = \widehat{\O}_{C,v}\!\otimes_{\O_{C,v}}\!Q \longrightarrow Q$ are geometrically regular. This means that $Q_v\otimes_Q K$ is regular for every finite field extension $K$ over $Q$. Since \quotes{regular} implies \quotes{reduced}, we conclude that $Q_v$ is separable over $Q$.
\end{Remark}

\begin{Proposition}\label{PROP.EQ}
In the notation of (\ref{Eq3.1}) the following statements are equivalent:
\begin{suchthat}
\item $\pi$ is semisimple.
\item $F$ is semisimple.
\item $F\xQv\cong F_v$ is semisimple.
\item $\pi_v$ is semisimple.
\item $E\xQv\cong E_v$ is semisimple.
\item $E$ is semisimple.
\end{suchthat}
\end{Proposition}

\begin{proof}
1. and 2. are equivalent by definition. So we show the equivalences from 2. to 6. 

Let $F$ be semisimple. Since $\Qv$ is separable over $Q$, we conclude that $F\xQv\cong \Qv\II{\pi_v}$ is semisimple by \cite[Corollaire~7.6/4]{Bou}.
Hence $\pi_v$ is semisimple by definition, and we showed in Lemma~\ref{PI-END-SEMISIMPLE}/2 that then $E_v\cong E\xQv$ is semisimple. Again by \cite[Corollaire~7.6/4]{Bou} this implies that $E$ is semisimple. Since $F\subset Z(E)$ is a finite dimensional $Q$-subalgebra of the center of $E$, we conclude by \cite[Corollaire de Proposition~6.4/9]{Bou} that $F$ is semisimple, and our proof is complete.
\end{proof}

\noindent {\it Remark.}
If more generally $\FF$ is defined over a finitely generated field, then one cannot consider $\pi$, $\pi_v$, nor $F$. Nevertheless 5 and 6 remain equivalent and are still implied by 3 due to the following well-known lemma. Namely $E_v$ is the commutant of $F_v$ in $\End_{Q_v}(V_v\FF)$. 
We thank O.\ Gabber for mentioning this fact to us and we include its proof for lack of reference.

\begin{Lemma}\label{Lemma3.5b}
Let $B$ be a central simple algebra of finite dimension over a field $K$ and let $F$ be a semisimple $K$-subalgebra of $B$. Then the commutant of $F$ in $B$ is semisimple.
\end{Lemma}

\begin{proof}
Let $F=\bigoplus_i F_i$ be the decomposition into simple constituents and let $e_i$ be the corresponding central idempotents, that is, $F_i=Fe_i$. Consider $B_i=e_i Be_i$ which is again central simple over $K$ by \cite[Corollaire~6.4/4]{Bou}, since if $I\subset B_i$ is a non-zero two sided ideal then $BIB$ contains $1$ and so $I$ contains the unit $e_i$ of $B_i$. By \cite[Th\'eor\`eme~10.2/2]{Bou} the commutant $E_i$ of $F_i$ in $B_i$ is simple. Clearly the commutant $E$ of $F$ in $B$ satisfies $E_i=e_iEe_i=Ee_i$ and $E=\bigoplus_i E_i$ proving the lemma.
\end{proof}

\begin{Corollary} \label{CorFisCenter}
Let\/ $\FF$ be an abelian $\tau$-sheaf over $\Fs$ of rank $r$ with semisimple Frobenius endomorphism\/ $\pi$. 
Then the algebra $F=Q(\pi)$ is the center of the semisimple algebra $E=\QEnd(\FF)$.
\end{Corollary}

\begin{proof}
Since $F_v$ is semisimple, we know by \cite[Proposition~5.1/1]{Bou} that the $F_v$-module $\VvFF$ is semisimple. The commutant of $F_v$ in $\End_{Q_v}(\VvFF)$ is $E_v$ by definition. Trivially $\VvFF$ is of finite type over $E_v$. Thus, by the theorem of bicommutation \cite[Corollaire~4.2/1]{Bou}, the commutant of $E_v$ in $\End(\VvFF)$ is again $F_v$. We conclude $Z(E_v) = E_v \cap F_v = F_v$ and we have $F\xQv = F_v = Z(E_v) = Z(E)\xQv$ by \cite[Corollaire de Proposition~1.2/3]{Bou}. Considering the dimensions, we obtain $\dim_Q F= \dim_Q Z(E)$.
Since $F\subset Z(E)$ and the dimensions are finite, we finish by $F=Z(E)$.
\end{proof}

\begin{Theorem}\label{Thm3.8} %\label{POINCARE-WEIL}
Let $\FF$ be an abelian $\tau$-sheaf over a finite field $L$.
\begin{suchthat}
\item If $\QEnd(\FF)$ is a division algebra over $Q$ then $\FF$ is simple. If in addition $\chr\ne\infty$ then both statements are equivalent.
\item If the characteristic point $\chr$ is different from $\infty$ then
$\FF$ is semisimple if and only if $\QEnd(\FF)$ is semisimple.
\end{suchthat}
\end{Theorem}

\begin{proof}
1. Let $\QEnd(\FF)=E$ be a division algebra and let $f\in\Hom(\FF,\FF')$ be the morphism onto a non-zero quotient sheaf $\FF'$ of $\FF$. We show that $f$ is an isomorphism. We know by \ref{FACTORSHEAF-FACTORSPACE} that \mbox{$f_v\in\Hom_\QvG(\VvFF,\VvFF')$} is surjective. By the semisimplicity of $E$ and Proposition~\ref{PROP.EQ}, $F_v$ is semisimple, and therefore $\VvFF$ is a finitely generated semisimple $F_v$-module. Thus we get a morphism $g_v\in\Hom_{Q_v[G]}(\VvFF',\VvFF)$ with \mbox{$f_v\circ g_v=\id_{\VvFF'}$}. Consider the integral Tate modules $T_v\FF$ and $T_v\FF'$. We can find some $n\in\N$ such that 
\[
v^n g_v \in \Hom_{\AvG}(T_v\FF',T_v\FF) \;\cong\; \Hom(\ulM(\FF'),\ulM(\FF))\otimes_A\Av
\]
and we choose $g\in\Hom(\ulM(\FF'),\ulM(\FF))\subset\QHom(\FF',\FF)$ with $g\equiv v^n g_v$ modulo $v^m$ for a sufficiently large $m>n$. If $g\circ f=0$ in $E$, then $f\circ g\circ f=0$, and therefore $f\circ g=0$ in $\QEnd(\FF')$ due to the surjectivity of $f$. This would imply
\[
v^n\cdot\id_{\VvFF'} \;=\; v^n(f_v\circ g_v) \;=\; f_v\circ(v^n g_v) \;\equiv\; f\circ g \;=\; 0 
\qquad\mbox{(modulo $v^m$)}
\]
which is a contradiction. Thus $g\circ f\not=0$ is invertible in $E$, and therefore $f$ is injective. By that, $f$ gives the desired isomorphism between $\FF'$ and $\FF$. The second assertion follows from Theorem~\ref{QEND-DIVISION-MATRIX}.

\noindent 
2. We already saw one direction in Theorem~\ref{QEND-DIVISION-MATRIX}/2. So let now $\QEnd(\FF)$ be semisimple and let 
\[
\QEnd(\FF)=\bigoplus_{j=1}^m\, M_{\lambda_j}(E_j)
\]
be the decomposition into full matrix algebras $M_{\lambda_j}(E_j)$ over division algebras $E_j$ over $Q$ $(1\le j\le m)$. For each $j$ we find $\lambda_j$ distinct idempotents $e_{j,1},\dots,e_{j,\lambda_j}\in M_{\lambda_j}(E_j)$ such that $e_{j,\alpha}\cdot\QEnd(\FF)\cdot e_{j,\alpha} = E_j$ for all $1\le\alpha\le\lambda_j$ with $\sum_{\alpha=1}^{\lambda_j} e_{j,\alpha}=1$ in $M_{\lambda_j}(E_j)$. Let $e_1,\dots,e_n$ denote all these idempotents, $n=\sum_{j=1}^m \lambda_j$, and choose a divisor $D$ on $C$ such that $e_i\in\Hom(\FF,\FF(D))$ for all $1\le i\le n$. Then $\sum_{i=1}^n e_i=\id_\FF$ in $\QEnd(\FF)$ and therefore 
\[
\begin{CD}
\FF\; @>{\sum_i e_i}>> \;{\displaystyle\bigoplus_{i=1}^n\, \im e_i}\; \subset \,\FF(D)\;.
\end{CD}
\]
The image $\FF_i:=\im e_i$ is an abelian $\tau$-sheaf by Proposition~\ref{IMISABELIAN} because $\chr\ne\infty$. Since $\sum_i e_i$ is injective it is an isogeny by \ref{PROP.1.42A}. Since $\QEnd(\FF_i) = e_i\cdot\QEnd(\FF)\cdot e_i$ is a division algebra, $\FF_i$ is a simple abelian $\tau$-sheaf by 1. Thus $\FF\approx\FF_1\oplus\cdots\oplus\FF_n$ gives the decomposition into a direct sum of simple abelian $\tau$-sheaves $\FF_i$ as desired. 
\end{proof}

\begin{Remark}\label{Rem3.10b}
Unfortunately the theorem fails if $L$ is not finite, as Example~\ref{Ex3.10c} below shows. The reason is, that then $E_v$ may still be semisimple while the image $F_v$ of $Q_v[G]$ in $\End_{Q_v}(V_v\FF)$ is not. Nevertheless, if one adds the assumption that $F_v$ is semisimple, the assertions of Theorem~\ref{Thm3.8} remain valid over an arbitrary field $L$. (See also the remark after Proposition~\ref{PROP.EQ}.)
\end{Remark}

\begin{Example}\label{Ex3.10c}
We construct a pure Anderson motive $\ulM$ over a non-finite field $L$ which is not semisimple, but has $\End(\ulM)=A$. Any associated abelian $\tau$-sheaf $\FF$ has $\QEnd(\FF)=Q$. Let $C=\PP^1_\Fq$, $A=\Fq[t]$ with $q>2$, and $L=\Fq(\alpha)$ where $\alpha$ is transcendental over $\Fq$. Let $M=A_L^{\oplus 2}$ and $\t=\matr{\alpha t}{0}{t}{t}$. Then $\ulM=(M,\t)$ is a pure Anderson motive of rank and dimension $2$. Clearly $\ulM$ is not simple, since $\ulM'=(A_L,\t'=t)$ is a quotient motive by projecting onto the second coordinate. We will see below that $\ulM$ is not even semisimple.

Let $\matr{e}{g}{f}{h}\in M_2(A_L)$ be an endomorphism of $\ulM$, that is,
\[
\matr{\alpha\s e+\s g}{\s g}{\;\;\alpha\s f+\s h}{\;\s h}\es=\es\matr{\alpha e}{\alpha g}{\;\;e+f}{\;\;g+h}\,.
\]
Choose $\beta\in\BF_q(\alpha)^\alg\setminus\BF_q(\alpha)$ satisfying $\beta^{q-1}=\alpha$ (for $\beta\notin\BF_q(\alpha)$ we use $q>2$). Then $\s g=\alpha g$ implies $g\in\beta\cdot\BF_q[t]$. Since also $g\in\BF_q(\alpha)[t]$ we must have $g=0$. Now $\s e=e$ and $\s h=h$ yielding $e,h\in\BF_q[t]$.

Let $\gamma\in\BF_q(\alpha)^\alg\setminus\BF_q(\beta)$ with $\gamma^q-\gamma=\beta$ and set $\tilde f:=\beta f-\gamma\cdot(e-h)$. Then $\alpha \s f-f=e-\s h=e-h$ implies $\s\tilde f-\tilde f=\beta^q\s f-\beta f-(\gamma^q-\gamma)(e-h)=\beta(\alpha\s f-f-(e-h))=0$. Thus $\tilde f\in\BF_q[t]$ and $\gamma\cdot(e-h)\in\BF_q(\beta)[t]$. So we must have $e=h$ and then $\beta f=\tilde f\in\BF_q[t]$ implies $f=0$. 
This shows that $\End(\ulM)=\Fq[t]=A$.

\forget{
Using the inclusion $A_L=\Fq(\alpha)[t]\subset\Fq(t)\dpl\alpha\dpr$ we write $g=\sum_{i=N}^\infty g_i\alpha^i$ with $N\in\Z$, $g_i\in\Fq(t)$, and $g_N\ne0$. Then $\s g=\sum_{i=N}^\infty g_i\alpha^{qi}$. Comparing the coefficients of $\alpha^i$ in the equation 
\[
\sum_{i=N}^\infty g_i\alpha^{qi}\;=\;\s g\;=\;\alpha g\;=\;\sum_{i=N}^\infty g_i\alpha^{i+1}
\]
and observing $q>2$ one easily sees that all $g_i$ must be zero. Therefore $g=0$ and $\s e=e$, $\s h=h$, hence $e,h\in\Fq[t]$. Now we write $f=\sum_{i=N}^\infty f_i\alpha^i\in\Fq(t)\dpl\alpha\dpr$ with $f_N\ne0$ to solve the remaining equation $\alpha\s f+\s h=e+f$, that is,
\[
\sum_{i=N}^\infty f_i\alpha^i-\sum_{i=N}^\infty f_i\alpha^{qi+1}\;=\;h-e\,.
\]
It follows that $N=0$, $f_0=h-e=f_i$ whenever $i=1+q+\ldots+q^j$ for some $j\in\N_0$ and $f_i=0$ otherwise. But for $e\ne h$ this solution $f$ does not belong to $\Fq(\alpha)[t]$ due to the well known characterization
\[
\Fq(t)(\alpha)\;=\;\bigl\{\,\sum_{i=N}^\infty f_i\alpha^i\in\Fq(t)\dpl\alpha\dpr:\;\exists\,m,n\in\Z_{>0}\text{ with }f_{i+n}=f_i\;\forall\,i\ge m\,\bigr\}\,.
\]
We conclude $e=h$ and $f=0$, that is, $\End(\ulM)=\Fq[t]=A$.
}

The same argument shows that $\ulM$ is not even semisimple. Namely, the projection $\ulM\to\ulM'$ has no section $\ulM'\to\ulM, 1\mapsto{f\choose 1}$, since there is no solution $f$ for the equation $\alpha t\s f+t=tf$.

It is also not hard to compute $F_v$ for instance at the place $v=(t-1)$. Let $z=t-1$ and $\beta\in L^\sep$ with $\beta^{q-1}=\alpha$, and consider the basis ${y/\beta\choose 0},{x\choose y}$ of the Tate module $T_v(\ulM)$ with 
\[
{x\choose y}\;=\;\sum_{i=0}^\infty{x_i\choose y_i}z^i\quad\text{and}\quad x_i,y_i\in L^\sep\,,\;y_0\ne0\,.
\]
They are subject to the equations $y=t\s y=(1+z)\s y$ and $x=\alpha t\s x+t\s y=\alpha(1+z)\s x+y$, that is,
\begin{eqnarray*}
y_i-y_i^q&=&y_{i-1}^q\,,\text{ and}\\[1mm]
x_i-\alpha x_i^q&=&\alpha x_{i-1}^q+y_i\,.
\end{eqnarray*}
There are elements $\gamma$ and $\delta$ of $G=\Gal(L^\sep/L)$ operating as $\gamma(y_i)=y_i$, $\gamma(x_i)=x_i$, $\gamma(\beta)=\beta/\eta$ for an $\eta\in\Fq^\times\setminus\{1\}$, respectively as $\delta(y_i)=y_i$, $\delta(\beta)=\beta$, $\delta(x_i)=x_i+y_i/\beta$. With respect to our basis of $T_v(\ulM)$ they correspond to matrices $\gamma_v=\matr{\eta}{0}{0}{1}$ and $\delta_v=\matr{1}{0}{1}{1}$. We conclude that $F_v$ is the $Q_v$-algebra of upper triangular matrices. Its commutant in $M_2(Q_v)$ equals $Q_v\cdot\Id_2\cong\End(\ulM)\otimes_A Q_v$.
\end{Example}

\noindent
{\it Remark.} If $q=2$ any pure Anderson motive of rank $\rk\ulM=2$ on $A=\BF_q[t]$, which is not semisimple has  $\End(\ulM)\supsetneq A$. One easily sees this by choosing a basis of $\ulM$ for which $\t$ has the form $\matr{\alpha (t-\theta)^d}{0}{\ast}{\beta(t-\theta)^d}$ with $\alpha,\beta,\theta\in L$. Then $\matr{0}{0}{\beta/\alpha}{0}$ is an endomorphism.

However, we expect that also for $q=2$ there are examples similar to \ref{Ex3.10c} (of $\rk\ulM\ge3$), although we have not tried to find one.

\bigskip

Let\/ $\FF$ be an abelian $\tau$-sheaf over $\Fs$ and let\/ $\Fss/\Fs$ be a finite field extension. The base extension
\[
\FF\xFss := (\F_i\otimes_{\O_{\CFs}}\!\!\!\O_{C_{\Fs\mbox{\raisebox{-0.5ex}{$\scriptscriptstyle'$}}}}, \P_i\otimes1,\t_i\otimes1)
\]
is an abelian $\tau$-sheaf over $\Fss$ with $\pi'=(\pi\otimes1)^t$ for $s'=s^t$, and we have a canonical isomorphism between $\VvFF$ and $\VvFF'$.

\bigskip

For the next result recall that an endomorphism $\varphi$ of a finite dimensional vector space $V$ over a field $K$ is called \emph{absolutely semisimple} if for every field extension $K'/K$ the endomorphism $\varphi\otimes1\in\End_{K'}(V\otimes_K K')$ is semisimple. The following characterization is taken from \cite[Proposition~9.2/4 and Proposition~9.2/5]{Bou}.

\begin{Lemma}\label{XY99}
Let\/ $K$ be a field and let\/ $V$ be a finite dimensional $K$-vector space. Let\/ $\varphi\in\End_K(V)$ be an endomorphism. 
\begin{suchthat}
\item $\varphi$ is absolutely semisimple, if and only if there exists a perfect field extension\/ $K'/K$ such that $\varphi\otimes1\in\End_{K'}(V\otimes_K K')$ is semisimple.
\item $\varphi$ is absolutely semisimple, if and only if its minimal polynomial is separable.
\end{suchthat}
\end{Lemma}

\medskip

\begin{Theorem}\label{Thm3.8b}
Let $\FF$ be an abelian $\tau$-sheaf over the finite field $\Fs$. Then there exists a finite field extension\/ $\Fss/\Fs$ whose degree is a power of $\charakt\Fs$ such that\/ $\FF\xFss$ has an absolutely semisimple Frobenius endomorphism. Thus if moreover\/ $\chr\ne\infty$ then $\FF\xFss$ is semisimple.
\end{Theorem}

\noindent
{\it Remark.} It suffices to take $[\BF_{s'}:\BF_s]$ as the smallest power of $\charakt\BF_s$ which is $\ge\rk\FF$.

\begin{proof}
Let $s'=s^t$ for some arbitrary $t\in\N$. Let $\FF':=\FF\xFss$ be the abelian $\tau$-sheaf over $\Fss$ induced by $\FF$. Let $v\in\Spec A$ be a place different from $\chr$. Over $\Qv^\alg$ we can write $\pi_v\in\End_\Qv(\VvFF)$ in Jordan normal form 
\[
\makebox[5em][r]{$B^{-1}\,(\pi_v\otimes1)\,B$} 
\;=\; 
\left(\begin{array}{c@{\;\;\;}c@{\;\;\;}c@{\;\;\;}c}
\lambda_1 & \ast & & 0 \\
& \lambda_2 & \ddots & \\
& & \ddots & \ast \\
0 & & & \lambda_r \\
\end{array}\right)
\]
for $B\in GL_r(\Qv^\alg)$ and for some $\lambda_j\in\Qv^\alg$, $1\le j\le r$. Thus, by a suitable choice of $t\in\N$ as a power of $\charakt\Fq$ (as in the remark), we can achieve that $\pi'_v=(\pi_v\otimes1)^t$ is of the form
\[
\makebox[5em][r]{$B^{-1}\,\pi'_v\,B$} 
\;=\; 
\left(\begin{array}{c@{\;\;\;}c@{\;\;\;}c@{\;\;\;}c}
\lambda_1^t & & & 0 \\
& \lambda_2^t & & \\
& & \ddots & \\
0 & & & \lambda_r^t \\
\end{array}\right)\;.
\]
Since $\Qv^\alg$ is perfect, we conclude by \ref{XY99}/1 that $\pi'_v$ and thus $\pi'$ is absolutely semisimple.
\end{proof}

The following corollary illustrates that, in contrast to endomorphisms of vector spaces, there is no need of the term \quotes{absolutely semisimple} for abelian $\tau$-sheaves or pure Anderson motives over finite fields.

\begin{Corollary} \label{Cor3.8c}
Let\/ $\FF$ be an abelian $\tau$-sheaf over $\Fs$ of characteristic different from $\infty$. If\/ $\FF$ is semisimple, then\/ $\FF\xFss$ is semisimple for every finite field extension $\Fss/\Fs$. The same is true for pure Anderson motives.
\end{Corollary}

\begin{proof}
Let $\FF$ be semisimple and let $\Fss/\Fs$ be a finite field extension with $s'=s^t$. We set $\FF':=\FF\xFss$. By \ref{Thm3.8} and \ref{PROP.EQ}, we know that $\QEnd(\FF)\xQv \cong \End_{\Qv\II{\pi_v}}(\VvFF)$ is semisimple. Since $\Qv\II{\pi_v^t}\subset\Qv\II{\pi_v}$ we conclude by \cite[Corollaire de Proposition~6.4/9]{Bou} that $\Qv\II{\pi_v^t}$ is semisimple, as well. As $\VvFF'=\VvFF$, we have $\pi'_v=\pi_v^t$, and therefore $\pi'_v$ is semisimple. Thus, by \ref{PROP.EQ}, $\QEnd(\FF')$ is semisimple and $\FF'$ is semisimple by \ref{Thm3.8}/2.
\end{proof}

% -----------------------------------------------------------------------------

\bigskip

\subsection{Zeta Functions and Reduced Norms} \label{SectZeta}

In this section we generalize Gekeler's results \cite{Gekeler} on Zeta functions for Drinfeld modules to pure Anderson motives.
But let us begin by recalling a few facts about reduced norms; see for instance \cite[\S 9]{Re}. Let $\ulM$ be a semisimple pure Anderson motive over a finite field and let $\pi$ be its Frobenius endomorphism. Then $F=Q(\pi)$ is the center of the semisimple algebra $E$ by Corollary~\ref{CorFisCenter}. Write $F=\bigoplus_i F_i$ and $E=\bigoplus_i E_i$ where the $F_i$ are fields and $E_i$ is central simple over $F_i$. Note that by~\ref{Thm3.8} the pure Anderson motive $\ulM$ decomposes correspondingly up to isogeny $\ulM\approx\bigoplus_i\ulM_i$ with $E_i=\End(\ulM_i)\otimes_AQ$. We apply \ref{Cor3.11b} to $\ulM_i$ and obtain $\sum_i[E_i:F_i]^{1/2}\cdot[F_i:Q]=r$. Let $f\in E$ and write it as $f=\sum_i f_i$ with $f_i\in E_i$. Choose for each $i$ a splitting field $K_i$ of $E_i$ with $\alpha_i:E_i\otimes_{F_i}K_i\isoto M_{n_i}(K_i)$ where $n_i^2=[E_i:F_i]$. The \emph{reduced norm of $f$} is then defined by
\[
N(f)\;:=\;nr_{E/Q}(f)\;:=\;\prod_i N_{F_i/Q}\bigl(\det\alpha_i(f_i\otimes 1)\bigr),
\]
where $N_{F_i/Q}$ is the usual field norm. The reduced norm is an element of $Q$ which is independent of the choices of $K_i$ and $\alpha_i$. It satisfies $N(a)=a^r$ for all $a\in Q$, and $N(f)\ne0$ if and only if $f\in E^{\SSC\times}$, that is, $f$ is a quasi-isogeny. If $f\in \End(\ulM)$ or more generally $f$ is contained in a finite $A$-algebra then $N(f)\in A$ since $A$ is normal.

\begin{Theorem}\label{Thm3.4.1a}
Let $\FF$ be a semisimple abelian $\tau$-sheaf over a finite field $L$ and let $f\in \QEnd(\FF)$ be a quasi-isogeny. Then
for any place $v\ne\chr,\infty$ of $Q$ we have $N(f)=\det V_vf$, the determinant of the endomorphism $V_vf\in\End_{Q_v}(V_v\FF)$. For $v=\infty\ne\chr$ we have $N(f)^l=\det V_\infty f$, where $l$ comes from Definition~\ref{DefTateMod} and satisfies $\dim_{Q_\infty}V_\infty\FF=l\cdot\rk\FF$.
\end{Theorem}

\begin{proof}
Clearly, if $t$ is a power of $q$ then $N(f^t)=\det V_vf^t$ implies $N(f)=\det V_vf$ since $1$ is the only $t$-th root of unity in $Q_v$ for $v\ne\infty$, and likewise for $v=\infty$.
Writing $V_vf$ in Jordan canonical form over $Q_v^\alg$ we find as in the proof of Theorem~\ref{Thm3.8b} a power $t$ of $q$ such that $V_vf^t$ is absolutely semisimple over $Q_v$ and hence its minimal polynomial is separable by \ref{XY99}.
Then $F_v(f^t)$ and $F(f^t)$ are semisimple by \cite[Proposition 9.1/1 and Corollaire 7.7/4]{Bou}. We now replace $f$ by $f^t$ and thus assume that $F(f)$ is semisimple. 

As is well known there is a semisimple commutative subalgebra $H=\bigoplus_i H_i$ of $E$ containing $F(f)$ with $\dim_{F_i}H_i=n_i$ and hence $\dim_QH=r$. Then $nr_{E/Q}(f)$ equals the determinant of the $Q$-endomorphism $\tilde f:x\mapsto fx$ of $H$. The reason for this is that $H_i\otimes_{F_i}K_i$ is still semisimple and commutative if we choose a splitting field $K_i$ which is separable over $F_i$. By Lemma~\ref{Lemma3.4.1b} below $H_i\otimes_{F_i}K_i$ is isomorphic to $K_i^{n_i}$ as left $H_i\otimes_{F_i}K_i$-modules, and this implies that $nr_{E_i/F_i}(f_i)=\det\alpha_i(f_i)=\det\tilde f_i$, the determinant of the $F_i$-endomorphism $\tilde f_i:x\mapsto f_i x$ of $H_i$, and $N(f)=\det\tilde f$ the determinant of the $Q$-endomorphism $\tilde f$ of $H$.

If $v\ne\infty$ then again by Lemma~\ref{Lemma3.4.1b}, $H_v$ is $H_v$-isomorphic to $V_v\FF$ and $N(f)=\det\tilde f=\det V_vf$.

If $v=\infty$ we embed $E_\infty^{\oplus l}$ into $\End_{Q_{\infty,L}[\phi]}\bigl(\ulTN_\infty(\FF)\bigr)$. Namely, if $(f^{(0)},\ldots,f^{(l-1)})\in E_\infty^{\oplus l}$, where $f^{(m)}=\bigl(f^{(m)}_i:\F_i\otimes_{\O_{C_L}}Q_{\infty,L}\to\F_i\otimes_{\O_{C_L}}Q_{\infty,L}\bigr)$, we set
\[
g_{ij}\;:=\;\left\{
\begin{array}{ll}
\P_{i-1}\circ\ldots\circ\P_j\circ f_j^{(i-j)} & \text{if }0\le j\le i\le l-1 \\[2mm]
z^k\P_i^{-1}\circ\ldots\circ\P_{j-1}^{-1}\circ f_j^{(l+i-j)} & \text{if }0\le i<j\le l-1\;.
\end{array}
\right.
\]
Then $g_{ij}:\F_j\otimes_{\O_{C_L}}Q_{\infty,L}\to\F_i\otimes_{\O_{C_L}}Q_{\infty,L}$ and a straightforward computation shows that the homomorphism $g=(g_{ij})_{i,j=0\ldots l-1}$ commutes with $\phi$ from (\ref{EQ.Phi}) on page~\pageref{EQ.Phi}, that is, $g$ is an element of $\End_{Q_{\infty,L}[\phi]}\bigl(\ulTN_\infty(\FF)\bigr)=\End_{Q_\infty[G]}(V_\infty\FF)$; use Proposition~\ref{Prop2.13'}. Now we apply Lemma~\ref{Lemma3.4.1b} to $H_\infty^{\oplus l}\subset E_\infty^{\oplus l}\subset\End_{Q_\infty}(V_\infty\FF)$, and we compute $N(f)^l=(\det\tilde f)^l=\det_{Q_\infty}\bigl(H_\infty^{\oplus l}\to H_\infty^{\oplus l}, h\mapsto fh\bigr)=\det V_\infty f$ as desired.
\end{proof}

\begin{Lemma}\label{Lemma3.4.1b}
Let $K$ be a field and let $H\subset M_n(K)$ be a semisimple commutative $K$-algebra with $\dim_K H=n$. Then as a (left) module over itself $H$ is isomorphic to $K^n$.
\end{Lemma}

\begin{proof}
Decomposing $H$ into a direct sum of fields $\bigoplus_\kappa L_\kappa$ and $K^n$ into a direct sum $\bigoplus_\lambda V_\lambda$ of simple $H$-modules, each $V_\lambda$ is isomorphic to an $L_{\kappa(\lambda)}$. The injectivity of $H\to M_n(K)$ and $\dim_K H=n$ imply that $H$ is isomorphic to $\bigoplus_\lambda \End_{L_{\kappa(\lambda)}}(V_\lambda)$ and a fortiori isomorphic as left module over itself to $K^n$.
\end{proof}

\begin{Theorem}\label{Prop3.4.1}
Let $\ulM$ be a semisimple pure Anderson motive of rank $r$ over a finite field $L$ and let $f\in \End(\ulM)$ be an isogeny. Then
\begin{suchthat}
\item 
$\dim_L\coker N(f)=r\cdot\dim_L\coker f$.
\item 
The ideal $\deg(f)=N(f)\cdot A$ is principal and has a canonical generator.
\item 
There exists a canonical dual isogeny $\dual{f}\in\End(\ulM)$ satisfying $f\circ\dual{f}=N(f)=\dual{f}\circ f$.
\end{suchthat}
\end{Theorem}

\noindent
{\it Remark.} 1. This shows that $N(1-\pi^n)\in A$ is the analog for pure Anderson motives of the number of rational points $X(\BF_{q^n})=\deg(1-{\rm Frob}_q^n)\in\Z$ on an abelian variety $X$ over the finite field $\Fq$; see also Theorem~\ref{ThmZeta} below.

2. The dual isogeny satisfies $\dual{(fg)}=\dual{g}\dual{f}$, because $N(fg)=N(f)N(g)$. Note however, that we cannot expect that $\dual{(f+g)} =\dual{f}+\dual{g}$ unless $r=2$ because for $f=a\in A$ we have $N(a)=a^r$ and $\dual{a}=a^{r-1}$.

\begin{proof}
1. Clearly for any $a\in A$ we have $\dim_L\ulM/a\ulM=r\cdot\dim_\Fq A/(a)=-r\cdot\infty(a)$ where $\infty(a)$ denotes the $\infty$-adic valuation of $a$. Now let $\FF$ be an abelian $\tau$-sheaf with $\ulM=\ulM(\FF)$, and let $f:\FF\to\FF(n\cdot\infty)$ for some $n$ be the isogeny induced by $f$. Using Theorem~\ref{Thm3.4.1a} we compute the dimension
\begin{eqnarray*}
l\cdot\dim_L\coker f & =& nrl-\dim_L{\TS\bigoplus_{j=0}^{l-1}\bigl(\F_j(n\cdot\infty)/f_j(\F_j)}\bigr)_\infty\es\\[1mm]
& = & nrl-\dim_L \ulTM_\infty\bigl(\FF(n\!\cdot\!\infty)\bigr)/\ulTM_\infty(f)\bigl(\ulTM_\infty(\FF)\bigr)\\[1mm]
& = & nrl-\dim_\Fq\bigl(T_\infty\FF(n\cdot\infty)/T_\infty f(T_\infty \FF)\bigr)\\[1mm]
& = & -\infty(\det V_\infty f)\es =\es -l\cdot\infty\bigl(N(f)\bigr)\,.
\end{eqnarray*}
Here the first equality follows from the identities $\F_j(n\cdot\infty)/f_j(\F_j)=\bigl(\F_j(n\cdot\infty)/f_j(\F_j)\bigr)_\infty\oplus\coker f$ and $\dim_L\bigl(\F_j(n\cdot\infty)/f_j(\F_j)\bigr)=\deg\F_j(n\cdot\infty)-\deg f_j(\F_j)=nr$. The second equality is the definition of $\ulTM_\infty$, and the third follows from the isomorphism $\ulTM_\infty(\FF)\otimes_{A_{\infty,L}}A_{\infty,L^\sep}\cong T_\infty\FF\otimes_{A_\infty}A_{\infty,L^\sep}$. The fourth equality follows from the elementary divisor theorem. From this we obtain 1.

\smallskip
\noindent
2. Let $v\ne\chr$ be a maximal ideal of $A$. Using Theorem~\ref{Thm3.4.1a} we compute the $v$-adic valuation of $N(f)$
\[
v(N(f))\;=\;v(\det T_vf)\;=\;\dim_{\BF_v}\bigl(T_v\ulM/T_vf(T_v\ulM)\bigr)\;=\;\dim_{\BF_v}\bigl((\coker f)_v\otimes_L L^\sep\bigr)^\t\;=\;v(\deg f)\,.
\]
Again the second equality follows from the elementary divisor theorem, the third equality comes from the fact that the $\t$-invariants of the $v$-primary  part $(\coker f)_v\otimes_L L^\sep$ are isomorphic to $T_v\ulM/T_vf(T_v\ulM)$, and the last equality is the definition of $\deg f$. From 1 and Lemma~\ref{Lemma1.7.7} we obtain
\begin{eqnarray*}
r\cdot\dim_\Fq A/\deg(f)&=& r\cdot\dim_L\coker f \es=\es \dim_L\coker N(f)\es\\
&=&\dim_L\bigl((A/N(f))^r\otimes_\Fq L\bigr)\es=\es r\cdot\dim_\Fq A/N(f)\,.
\end{eqnarray*}
From the identity $\dim_\Fq A/\Fa=\sum_v[\BF_v:\Fq]\cdot v(\Fa)$ for any ideal $\Fa\subset A$ we conclude $\chr(\deg f)=\chr(N(f))$ and therefore $\deg (f)=N(f)\cdot A$.

\smallskip
\noindent
Finally 3 is immediate since $N(f)$ annihilates $\coker f$ by Proposition~\ref{Prop3.28a}.
\end{proof}

\begin{Remark}\label{Rem3.32b}
We do not know of a proof of 1 and 2 for arbitrary pure Anderson motives which does not make use of the associated abelian $\tau$-sheaf $\FF$. In the special case when $\ulM$ comes from a Drinfeld module, Gekeler~\cite[Lemma 3.1]{Gekeler} argued that both sides of the equation in 2 are extensions to $E$ of the $\infty$-adic valuation on $Q$. But this argument fails in general, since there may be more than one such extension as one sees from Example~\ref{Ex3.15} below.
\end{Remark}

\begin{Corollary}\label{Cor3.4.1c}
Let $\ulM$ be a semisimple pure Anderson motive of dimension $d$ over a finite field $L$ and let $\pi$ be its Frobenius endomorphism. Let $v\ne\chr$ be a maximal ideal of $A$ and let $\chi_v$ be the characteristic polynomial of $\pi_v$. Then
\begin{suchthat}
\item
$\chi_v\in A[x]$ is independent of $v$ and $\chi_v(a)\cdot A=\det V_v(a-\pi)\cdot A=\deg(a-\pi)$ for every $a\in A$,
\item 
$\chr^{d\cdot[L:\BF_\chr]}=\deg(\pi)=\chi_v(0)\cdot A=N(\pi)\cdot A$ is principal.
\end{suchthat}
\end{Corollary}

\begin{proof}
1 is a direct consequence of Theorems~\ref{Thm3.4.1a} and \ref{Prop3.4.1} and the Lagrange interpolation theorem applied to the fact that $\chi_v(a)=N(a-\pi)=\chi_w(a)\in A$ for all $a\in A$.\\
2 follows from the fact that $\coker\pi$ is supported on $\chr$ and from the equation $\dim_L\coker\pi=[L:\Fq]\cdot\dim_L\coker\t=d\cdot[L:\Fq]$.
\end{proof}

\begin{Definition}\label{DefZetaFkt}
We define the \emph{Zeta function} of a pure Anderson motive $\ulM$ over a finite field $\Fs$ as
\[
Z_\ulM(t)\;:=\;\prod_{0\le i\le r} \det(1-t\wedge^i\pi_v)^{(-1)^{i+1}}
\]
where $\chr\ne v\in \Spec A$ is a maximal ideal and $\wedge^i\pi_v\in\End_{Q_v}(\wedge^iV_v\ulM)$. 
\end{Definition}

By \ref{Cor3.4.1c}/1 the Zeta function $Z_\ulM(t)$ is independent of the place $v$ and lies in $Q(t)$. This also follows from work of B\"ockle~\cite{Boeckle02} and Gardeyn~\cite[\S7]{Gardeyn2}. The name ``Zeta function'' is justified by the following theorem (see also the remark after Theorem~\ref{Prop3.4.1}).

\begin{Theorem} \label{ThmZeta}
If $\ulM$ is semisimple and $\sum_i a_it^i$ is the power series expansion of $t\frac{d}{dt}\log Z_\ulM(t)$, then $a_i=N(1-\pi^i)\in A$.
\end{Theorem}

\begin{proof}
By standard arguments $a_i=\det(1-\pi_v^i)$; see \cite[Lemma 5.6]{Gekeler}. Now our assertion follows from Theorem~\ref{Thm3.4.1a} 
\end{proof}

This Zeta function satisfies the Riemann hypothesis:

\begin{Theorem}\label{ThmRH}
In an algebraic closure of $Q_\infty$ all eigenvalues of $\wedge^i\pi_v\in\End_{Q_v}(\wedge^iV_v\ulM)$ have the same absolute value $(\#\Fs)^{i\weight(\ulM)}$.
\end{Theorem}

\begin{proof}
This was proved by Goss~\cite[Theorem~5.6.10]{Go} for $i=1$ and follows for the remaining $i$ by general arguments of linear algebra.
\end{proof}

% -----------------------------------------------------------------------------

\bigskip

\subsection{A Quasi-Isogeny Criterion}

Similarly to the theory for abelian varieties, the characteristic polynomials of the Frobenius endomorphisms on the associated Tate modules play an important role for the study of abelian $\tau$-sheaves. For example, we can decide on quasi-isogeny of two abelian $\tau$-sheaves $\FF$ and $\FF'$ just by considering these characteristic polynomials.

\begin{Theorem}\label{THEOREM-1}
Let\/ $\FF$ and\/ $\FF'$ be abelian $\tau$-sheaves over $\Fs$ with respective Frobenius endomorphisms\/ $\pi$ and\/ $\pi'$, and let $\mu_\pi$ and $\mu_{\pi'}$ be their minimal polynomials over $Q$. Let\/ $v\in C$ be a place different from $\infty$ and $\chr$. Let\/ $\chi_v$ and\/ $\chi'_v$ be the characteristic polynomials of\/ $\pi_v$ and\/ $\pi'_v$, respectively, and let $G:=\Gal(L^\sep/L)$. Assume in addition that $\chr\ne\infty$, or that $\FF$ and $\FF'$ have the same weight. 
\begin{suchthat}
\item Consider the following statements:
\begin{suchthat}
\item[\labelenumi 1. ] $\FF'$ is quasi-isogenous to an abelian quotient $\tau$-sheaf of\/ $\FF$.
\item[\labelenumi 2. ] $\VvFF'$ is $G$-isomorphic to a $G$-quotient space of\/ $\VvFF$.
\item[\labelenumi 3. ] $\chi'_v$ divides $\chi_v$ in $Q_v[x]$.
\item[\labelenumi 4. ] $\mu_{\pi'}$ divides $\mu_\pi$ in $Q[x]$ and $\rk\FF'\le\rk\FF$
\end{suchthat}%\vspace{3\parsep}
\noindent
$\begin{array}{rc@{\:\:\:}c@{\:\:\:}c@{\:\:}c@{\:\:}c@{\:\:}c@{\:\:}c@{\:\:\quad}l}
\mbox{We have}&\mbox{1.1} &\Rightarrow& \mbox{1.2} &\Rightarrow& \mbox{1.3} & \mbox{and} & \mbox{1.4} & \makebox[0.5\textwidth][l]{always,} \\
& & & \mbox{1.2} &\Leftarrow & \mbox{1.3} & & & \mbox{if\/ $\pi_v$ and\/ $\pi'_v$ are semisimple,} \\
& & & \mbox{1.2} & \Leftarrow & \mbox{1.3} &\Leftarrow& \mbox{1.4} & \mbox{if\/ $\mu_\pi$ is irreducible in\/ $Q[x]$,}\\
&\mbox{1.1} &\Leftarrow & \mbox{1.2} & & & & & \mbox{if the characteristic is different from $\infty$.}
\end{array}$
\smallskip
\item Consider the following statements:
\begin{suchthat}
\item[\labelenumi 1. ] $\FF$ and $\FF'$ are quasi-isogenous.
\item[\labelenumi 2. ] $\VvFF$ and $\VvFF'$ are $G$-isomorphic.
\item[\labelenumi 3. ] $\chi_v=\chi'_v$.
\item[\labelenumi 4. ] $\mu_\pi=\mu_{\pi'}$ and $\rk\FF=\rk\FF'$.
\item[\labelenumi 5. ] There is an isomorphism of $Q$-algebras $\QEnd(\FF)\cong\QEnd(\FF')$ mapping $\pi$ to $\pi'$.
\item[\labelenumi 6. ] There is a $Q_v$-isomorphism $\QEnd(\FF)\otimes_Q Q_v\cong\QEnd(\FF')\otimes_Q Q_v$ mapping $\pi_v$ to $\pi'_v$.
\item[\labelenumi 7. ] If $\chr\ne\infty$ also consider the statement $Z_{\ulM(\FF)}=Z_{\ulM(\FF')}$.
\end{suchthat}%\vspace{3\parsep}
\noindent
$\begin{array}{rc@{\:\:}c@{\:\:}c@{\:\:}c@{\:\:}c@{\:\:}c@{\:\:}c@{\:\:}l@{\:\:}c@{\:\:}c@{\:\:\quad}l}
\mbox{\hspace{-1ex}We have}&\mbox{2.1} &\Leftrightarrow& \mbox{2.2} &\Rightarrow& \multicolumn{4}{l}{\mbox{\hspace{-1.07ex}2.3\,,\,2.4\,,\,2.5}} & & & \makebox[0.5\textwidth][l]{always,} \\
& & & & & & & \multicolumn{2}{r}{\mbox{2.5}} &\Rightarrow &\mbox{2.6} & \makebox[0.5\textwidth][l]{always,} \\
& & & & & \multicolumn{5}{l}{\hspace{-1.07ex}\mbox{2.3} \quad\Leftrightarrow \quad \mbox{2.7}} & & \mbox{if the characteristic is different from $\infty$,} \\
& & & \mbox{2.2} &\Leftarrow & \mbox{2.3}& &\Leftarrow & & &\mbox{2.6} & \mbox{if\/ $\pi_v$ and\/ $\pi'_v$ are semisimple,} \\
& & & \mbox{2.2} & \Leftarrow & \mbox{2.3}&\Leftarrow& \mbox{2.4} & \multicolumn{2}{l}{\hspace{-1.1ex}\Leftarrow} & \mbox{2.6} & \mbox{if $\mu_\pi$ and $\mu_{\pi'}$ are irreducible in $Q[x]$.}
\end{array}
$
\end{suchthat}
\end{Theorem}

\begin{proof}
1. For the implication $1.1\Rightarrow 1.2$ without loss of generality, $\FF'$ can itself be considered as abelian quotient $\tau$-sheaf of $\FF$ and the implication follows from Proposition~\ref{FACTORSHEAF-FACTORSPACE}. The implication $1.2\Rightarrow 1.3$ is obvious.

For $1.2\Rightarrow 1.4$ note that $\mu_\pi$ is also the minimal polynomial of $\pi_v$ over $Q_v$ by Lemma~\ref{Lemma3.4}. By Proposition~\ref{PI-IS-QISOG} statement 1.2 implies $\mu_{\pi}(\pi'_v)=0$, whence 1.4.

For $1.3\Rightarrow 1.2$ let $\pi_v$ and $\pi'_v$ be semisimple. Let $\chi_v=\mu_1\cdot\ldots\cdot\mu_n$ and $\chi'_v=\mu'_1\cdot\ldots\cdot\mu'_{n'}$ be the factorization in $Q_v[x]$ into irreducible factors and set $V_i:=Q_v[x]/(\mu_i)$ and $V'_i:=Q_v[x]/(\mu'_i)$.
Then we can decompose $\VvFF=V_1\oplus\cdots\oplus V_n$ and $\VvFF'=V'_1\oplus\cdots\oplus V'_{n'}$. Since $\chi'_v$ divides $\chi_v$, we can now easily construct a surjective $G$-morphism from $\VvFF$ onto $\VvFF'$ which gives the desired result.

Next if $\mu_\pi$ is irreducible, 1.4 implies $\mu_{\pi'}=\mu_\pi$ and 1.3 follows from Corollary~\ref{Cor3.11b}. It further follows from Proposition~\ref{PROP.EQ} that $\pi_v$ and $\pi'_v$ are semisimple and this implies 1.2 by the above.

For $1.2\Rightarrow 1.1$ we first do not assume that $\chr\ne\infty$. Let $f_v:\, \VvFF\rightarrow\VvFF'$ be a surjective morphism of $\QvG$-modules. We may multiply $f_v$ by a suitable power of $v$ to get a morphism $f_v:\, T_v\FF\rightarrow T_v\FF'$ of the integral Tate modules which is not necessarily surjective, but satisfies $v^n T_v\FF'\subset f_v(T_v\FF)$ for a sufficiently large $n$. Let $\ulM:=\bigl(\Gamma(C_L\setminus\{\infty\},\F_0),\,\P_0^{-1}\circ\t\bigr)$. This is a ``$\t$-module on $A$'' in the sense of Definition~\ref{Def1.16}. If $\chr\ne\infty$ then $\ulM$ is the pure Anderson motive $\ulM(\FF)$ associated with $\FF$ in (\ref{Eq1.1}). Also let $\ulM':=\bigl(\Gamma(C_L\setminus\{\infty\},\F'_0),\,\P_0'{}^{-1}\circ\t'\bigr)$. By Theorem~\ref{TATE-CONJECTURE-MODULES} (or Theorem~\ref{TATE-CONJECTURE-MODULES} if $\chr\ne\infty$), $f_v$ lies inside $\Hom(\ulM,\ulM')\otimes_A A_v$, so we can approximate $f_v$ by some $f\in\Hom(\ulM,\ulM')$ with $T_v(f)\equiv f_v$ modulo $v^{n+1} T_v\ulM'$. Since $v^n T_v\ulM'\subset f_v(T_v \ulM)$ we find inside $\im T_v(f)$ generators of $v^n T_v\ulM'/v^{n+1} T_v\ulM'$. They generate an $\Av$-submodule of $v^nT_v\ulM'$ whose rank must at least be $r'$ since $v^n T_v\ulM'/v^{n+1} T_v\ulM'\cong (\Av/v\Av)^{r'}$. Thus $\im T_v(f)$ has rank $r'$. Either by assumption or by Corollary~\ref{Cor2.9b} if $\chr\ne\infty$, both $\FF$ and $\FF'$ have the same weight. So by Proposition~\ref{CONNECTION}/1, $f$ comes from a quasi-morphism $f\in\QHom(\FF,\FF')$, that is, a morphism $f:\FF\to\FF'(D)$ for a suitable divisor $D$. Now we finally assume that the characteristic is different from $\infty$. By Proposition~\ref{IMISABELIAN}, the image $\im \bigl(f:\FF\to\FF'(D)\bigr)$ is an abelian factor $\tau$-sheaf of $\FF$ and $\im f\rightarrow\FF'(D)$ is an injective morphism between abelian $\tau$-sheaves of the same rank and weight, hence an isogeny by Proposition~\ref{PROP.1.42A}.

\medskip

\noindent
2. A large part of 2 follows from 1. We prove the rest. To show $2.2\Rightarrow 2.1$ without the hypothesis on the characteristic, we just replace the last argument of the proof of $1.2\Rightarrow 1.1$ by the following: Since $r=\dim_{Q_v}V_v\FF=\dim_{Q_v}V_v\FF'=r'$, the morphism $f:\FF\to\FF'(D)$ is an injective morphism between abelian $\tau$-sheaves of the same rank and weight, hence an isogeny by Proposition~\ref{PROP.1.42A}.

For the implication $2.1\Rightarrow 2.5$ let $g\in\QIsog(\FF,\FF')$. Then the map $\QEnd(\FF)\to\QEnd(\FF')$ sending $f\mapsto gfg^{-1}$ is an isomorphism with $\pi'=g\pi g^{-1}$. The implication $2.5\Rightarrow 2.6$ is obvious.

For the implication $2.3\Rightarrow 2.7$ note that knowledge of $\chi_v$ yields the knowledge of $\det(1-t\wedge^i\pi_v)$ and thus of $Z_{\ulM(\FF)}$ by linear algebra. Conversely we know from Theorem~\ref{ThmRH} that all zeroes of $\det(1-t\wedge^i\pi_v)$ have absolute value $s^{-i\weight(\FF)}$ in an algebraic closure of $Q_\infty$. So we can recover $\chi_v$ from $Z_{\ulM(\FF)}$ by simply looking at this absolute value. This proves $2.3\Leftarrow 2.7$.

Next if $\pi_v$ and $\pi'_v$ are semisimple $2.6\Rightarrow 2.3$ follows from Lemma~\ref{PROP.15}/2, and $2.3\Rightarrow 2.2$ was already established in 1.

Finally if $\mu_\pi$ and $\mu_{\pi'}$ are irreducible, 2.4 follows from 2.6 by Corollary~\ref{Cor3.11b} since $\mu_\pi$ is also the minimal polynomial of $\pi_v$ over $Q_v$ by Lemma~\ref{Lemma3.4}. Also 2.3 follows from 2.4 by Corollary~\ref{Cor3.11b} and $\pi_v$ and $\pi'_v$ are semisimple, so $2.3\Rightarrow 2.2$ by the above.
\end{proof}

% -----------------------------------------------------------------------------

\bigskip

\subsection{The quasi-endomorphism ring}

In this section we study the structure of $\QEnd(\FF)$ for a semisimple abelian $\tau$-sheaf $\FF$ over a finite field and calculate the local Hasse invariants of $\QEnd(\FF)$ as a central simple algebra over $Q(\pi)$. For a detailed introduction to central simple algebras, Hasse invariants and the Brauer group, we refer to \cite[Ch. 7, \S\S 28--31]{Re}.

\begin{Theorem}\label{THEOREM-2}
Let\/ $\FF$ be an abelian $\tau$-sheaf over the finite field $\Fs$ of rank $r$ with semisimple Frobenius endomorphism\/ $\pi$, that is, $Q(\pi)$ is semisimple. Let\/ $v\in C$ be a place different from $\infty$ and from the characteristic point $\chr$. Let\/ $\chi_v$ be the characteristic polynomial of\/ $\pi_v$. 
\begin{suchthat}
\item \label{THEOREM-2.1}
The algebra $F=Q(\pi)$ is the center of the semisimple algebra $E=\QEnd(\FF)$.\vspace{\parsep}
\item \label{THEOREM-2.3}
We have $\quad r \,\le\, \II{E:Q} \,=\, r_\Qv(\chi_v,\chi_v) \,\le\, r^2\;.$
\item \label{THEOREM-2.4}
Consider the following statements:
\begin{suchthat}
\item[\ref{THEOREM-2.4}.1. ] $E=F$.
\item[\ref{THEOREM-2.4}.2. ] $E$ is commutative.
\item[\ref{THEOREM-2.4}.3. ] $\II{F:Q} = r$.
\item[\ref{THEOREM-2.4}.4. ] $\II{E:Q} = r$.
\item[\ref{THEOREM-2.4}.5. ] $\chi_v$ has no multiple factor in $\Qv\II{x}$.
\item[\ref{THEOREM-2.4}.6. ] $\chi_v$ is separable.
\end{suchthat}%\vspace{3\parsep}
\noindent
$\begin{array}{rc@{\:\:}c@{\:\:}c@{\:\:}c@{\:\:}c@{\:\:}c@{\:\:}c@{\:\:}c@{\:\:}c@{\:\:}c@{\:\:}c@{\:\:\quad}l}
\text{We have} &\mbox{\ref{THEOREM-2.4}.1} &\Leftrightarrow& \mbox{\ref{THEOREM-2.4}.2} &\Leftrightarrow& \mbox{\ref{THEOREM-2.4}.3} &\Leftrightarrow& \mbox{\ref{THEOREM-2.4}.4} &\Leftrightarrow& 
\mbox{\ref{THEOREM-2.4}.5} &\Leftarrow & \mbox{\ref{THEOREM-2.4}.6} & \makebox[0.55\textwidth][l]{always,} \\
& & & & & & & & & \mbox{\ref{THEOREM-2.4}.5} &\Rightarrow& \mbox{\ref{THEOREM-2.4}.6} & \mbox{if\/ $\pi_v$ is absolutely semisimple.}
\end{array}$
\smallskip

\item \label{THEOREM-2.5}
Consider the following statements:
\begin{suchthat}
\item[\ref{THEOREM-2.5}.1. ] $F=Q$.
\item[\ref{THEOREM-2.5}.2. ] $E$ is a central simple algebra over $Q$.
\item[\ref{THEOREM-2.5}.3. ] $\II{E:Q} = r^2$.
\item[\ref{THEOREM-2.5}.4. ] $\chi_v$ is the $r$-th power of a linear polynomial in $\Qv\II{x}$.
\item[\ref{THEOREM-2.5}.5. ] $\chi_v$ is purely inseparable.
\end{suchthat}%\vspace{3\parsep}
$\begin{array}{rc@{\:\:}c@{\:\:}c@{\:\:}c@{\:\:}c@{\:\:}c@{\:\:}c@{\:\:}c@{\:\:}c@{\:\:\quad}l}
\text{We have}& \mbox{\ref{THEOREM-2.5}.1} &\Leftrightarrow& \mbox{\ref{THEOREM-2.5}.2} &\Leftrightarrow& \mbox{\ref{THEOREM-2.5}.3} &\Leftrightarrow& 
 \mbox{\ref{THEOREM-2.5}.4} &\Rightarrow& \mbox{\ref{THEOREM-2.5}.5} & \makebox[0.55\textwidth][l]{always,} \\
& & & & & & & \mbox{\ref{THEOREM-2.5}.4} &\Leftarrow & \mbox{\ref{THEOREM-2.5}.5} & \mbox{if\/ $\pi_v$ is absolutely semisimple.}
\end{array}$
\smallskip
If\/ \ref{THEOREM-2.5}.2 holds and moreover the characteristic point $\chr:=c(\Spec\Fs)\in\CFs$ is different from $\infty$, $E$ is characterized by\/ $\inv_\infty E=\weight(\FF)$, $\inv_\chr E=-\weight(\FF)$ and\/ $\inv_v E=0$ for any other place $v\in C$.
\item \label{THEOREM-2.6}
In general the local Hasse invariants of $E$ at the places $v$ of $F$ equal $\inv_v E=-\frac{[\BF_v:\BF_q]}{[\BF_s:\BF_q]}\cdot v(\pi)$. In particular
\[
\inv_v E\,=\,\left\{\begin{array}{ll} 0 & \text{ if }v\nmid\chr\infty\,,\\
\weight(\FF)\cdot[F_v:Q_\infty] & \text{ if }v|\infty \text{ and }\chr\ne\infty\,.
\end{array}\right.
\]
(Here $F_v$ denotes the completion of $F$ at the place $v$ and $\BF_v$ is the residue field of the place $v$.)
\end{suchthat}
\end{Theorem}

\begin{Remark}\label{Rem3.14b}
If $\chr\ne\infty$ and $\FF$ is an elliptic sheaf, that is, $d=1$ and $\ulM(\FF)$ is the Anderson motive of a Drinfeld module, Gekeler~\cite[Theorem~2.9]{Gekeler} has shown that there is exactly one place $v$ of $F$ above $\chr$, and exactly one place $w$ of $F$ above $\infty$, and that $\inv_w E=[F:Q]\cdot\weight(\FF)$ and $\inv_v E=-[F:Q]\cdot\weight(\FF)$.
Note that Gekeler actually computes the Hasse invariants of the endomorphism algebra of the Drinfeld module. So his invariants differ from ours by a minus sign, since passing from Drinfeld modules to abelian $\tau$-sheaves is a contravariant functor, see \cite[Theorem~3.2.1]{BS}.
\end{Remark}

\begin{Corollary}\label{Cor3.13b}
Let $\FF$ be an abelian $\tau$-sheaf over the smallest possible field $L=\Fq$ such that $\QEnd(\FF)$ is a division algebra. Then $\QEnd(\FF)$ is commutative and equals $Q(\pi)$.
\end{Corollary}

\begin{proof}
$\QEnd(\FF)$ is a central division algebra over $F$ by Theorem~\ref{THEOREM-2}, which splits at all places of $F$ by \ref{THEOREM-2}/\ref{THEOREM-2.6}, hence equals $F$.
\end{proof}

\begin{proof}[\Proofof{of Theorem \ref{THEOREM-2}}]
\ref{THEOREM-2.1} was already proved in Corollary~\ref{CorFisCenter}. 

\smallskip\noindent
\ref{THEOREM-2.3}. Let
\[
\chi_v \,=\, \prod_{i=1}^n \,\mu_i^{m_i} \quad\in\Qv\II{x}
\]
with distinct irreducible $\mu_i\in\Qv\II{x}$ and $m_i>0$ for $1\le i\le n$. Then $\sum_{i=1}^n m_i\cdot\deg\mu_i=\deg\chi_v=r$, and by Theorem \ref{Thm3.5a} we have $\II{E:Q}=r_\Qv(\chi_v,\chi_v)=\sum_{i=1}^n m_i^2\cdot\deg\mu_i$. The result now follows from the obvious inequalities
\begin{equation}\label{Eq9.1}
r\;=\;\sum_{i=1}^n m_i\cdot\deg\mu_i \;\stackrel{(1)}{\le}\; 
\sum_{i=1}^n m_i^2\cdot\deg\mu_i \;\stackrel{(2)}{\le}\; 
\left(\sum_{i=1}^n m_i\cdot\deg\mu_i\right)^2 \;=\;r^2\, .
\end{equation}

\smallskip\noindent
\ref{THEOREM-2.4}. Since $F=Z(E)$, the equivalence $\ref{THEOREM-2.4}.1\Leftrightarrow \ref{THEOREM-2.4}.2$ is evident. We have equality in (1) of Equation (\ref{Eq9.1}) if and only if $m_i=1$ for all $1\le i\le s$ which establishes the equivalence $\ref{THEOREM-2.4}.4\Leftrightarrow \ref{THEOREM-2.4}.5$. In order to prove $\ref{THEOREM-2.4}.5\Rightarrow \ref{THEOREM-2.4}.3$ we consider the minimal polynomial $\mu_v$ of $\pi_v$ over $\Qv$. If $\chi_v$ has no multiple factor, then $\mu_v=\chi_v$ and therefore $\II{F:Q}=\II{\Qv(\pi_v):\Qv}=r$. Next $\ref{THEOREM-2.4}.3\Rightarrow \ref{THEOREM-2.4}.1$ because $F\subset E$ and $(\dim_{Q_v}F_v)(\dim_{Q_v}E_v)=\dim_{Q_v}\End_{Q_v}(V_v\FF)=r^2$ by \cite[Th\'eor\`eme 10.2/2]{Bou}, since $E_v$ is the commutant of $F_v$ in $\End_{Q_v}(V_v\FF)$. Note that $\ref{THEOREM-2.4}.3\Rightarrow \ref{THEOREM-2.4}.1$ also follows from Lemma~\ref{Lemma3.4.1b}. Conversely $\ref{THEOREM-2.4}.1\Rightarrow \ref{THEOREM-2.4}.4$ because $E=F$ implies $r\ge\II{\Qv(\pi_v):\Qv}=\II{F:Q}=\II{E:Q}\ge r$. For $\ref{THEOREM-2.4}.5\Rightarrow \ref{THEOREM-2.4}.6$ we use Lemma~\ref{XY99}/2 as we know that $\chi_v=\mu_v$. $\ref{THEOREM-2.4}.6\Rightarrow \ref{THEOREM-2.4}.5$ is clear.

\smallskip\noindent
\ref{THEOREM-2.5}. If $F=Q$, then $E$ is simple with center $Q$, so $E$ is a central simple algebra over $Q$. Since $F=Z(E)$, the converse is obvious. This shows $\ref{THEOREM-2.5}.1\Leftrightarrow \ref{THEOREM-2.5}.2.$ We have equality in (2) of (\ref{Eq9.1}) if and only if $n=1$, $\deg\mu_1=1$ and $m_1=r$ which establishes $\ref{THEOREM-2.5}.3\Leftrightarrow \ref{THEOREM-2.5}.4.$ In order to connect $\ref{THEOREM-2.5}.1\Leftrightarrow \ref{THEOREM-2.5}.2$ with $\ref{THEOREM-2.5}.3\Leftrightarrow \ref{THEOREM-2.5}.4$ let $\chi_v$ be a power of a linear polynomial. By \cite[Proposition~9.1/1]{Bou} the minimal polynomial of $\pi_v$ over $\Qv$ is linear and thus $F=Q$. The converse is trivial. For $\ref{THEOREM-2.5}.5\Rightarrow \ref{THEOREM-2.5}.4$ we use again \ref{XY99}/2 to see that $\mu_v$ is linear. $\ref{THEOREM-2.5}.4\Rightarrow \ref{THEOREM-2.5}.5$ is clear.

The statement about the Hasse invariants follows from \ref{THEOREM-2.6}. Nevertheless, we give a separate proof in case $(k,l)=1$ using Tate modules, since this is much shorter here and exhibits a different technique than \ref{THEOREM-2.6}.
By the Tate conjecture \ref{TATE-CONJECTURE}, $E\otimes_Q Q_v$ is isomorphic to $\End_{Q_v}(V_v\FF)\cong M_r(Q_v)$ for all places $v\in C$ which are different from $\chr$ and $\infty$, so the Hasse invariants of $E$ at these places are $0$. Since the sum of all Hasse invariants is $0$ (modulo $1$), we only need to calculate $\inv_\infty E$.

As a first step, we show that $\Ff_{q^l}$ is contained in $\Fs$. In our situation, $\pi$ lies inside $Q$. Thus, by \ref{ThmRH} we get $s^{k/l}=|\,\pi\,|_\infty=q^m$ for some $m\in\Z$ as $|\,Q_\infty^{\times}\,|_\infty=q^\Z$. Since $q^e=s$, we conclude that $e\cdot k/l=m\in\Z$ and hence $l \mid e$, since $k$ and $l$ are assumed to be relatively prime. Therefore $\Ff_{q^l}\subset\Ff_{q^e}=\Fs$.

Consider the rational Tate module $V_\infty(\FF)$ at $\infty$ and the isomorphism of $Q_\infty$-algebras 
\[
E\otimes_Q Q_\infty\,\cong\,\End_{\Delta_\infty[G]}(V_\infty\FF)\,=\,\End_{\Delta_\infty}(V_\infty\FF)
\]
from Theorem~\ref{TATE-CONJECTURE}.
Since $\dim_{Q_\infty}\Delta_\infty=l^2$ and $\dim_{Q_\infty}V_\infty\FF=rl$, we conclude that $V_\infty\FF$ is a left $r/l$-dimensional $\Delta_\infty$-vector space and hence isomorphic to $\Delta_\infty^{r/l}$. Thus we have 
\[
E\otimes_Q Q_\infty \,\cong\, \End_{\Delta_\infty}(\Delta_\infty^{r/l}) \,=\, M_{r/l}(\End_{\Delta_\infty}(\Delta_\infty))\,=\,M_{r/l}(\Delta_\infty^{\op})\,.
\]
\forget{
We claim that $\End_{\Delta_\infty}(\Delta_\infty) = \Delta_\infty^{\op}$ denoting the opposite division algebra of $\Delta_\infty$. Let $m\in\End_{\Delta_\infty}(\Delta_\infty)$ and consider $f:=m(1)\in\Delta_\infty$. For $g\in\Delta_\infty$ we have
\[
m(g) \,=\, m(g\cdot 1) \,=\, g\cdot m(1) \,=\, g\cdot f
\]
since $m$ is a $\Delta_\infty$-endomorphism. Thus, $m$ is uniquely determined by $f$, and the map $\alpha:\, \End_{\Delta_\infty}(\Delta_\infty)\rightarrow\Delta_\infty$, $m\mapsto f$ is an isomorphism of $Q_\infty$-vector spaces. Let now $n\in\End_{\Delta_\infty}(\Delta_\infty)$ and let $g:=n(1)$. Then
\[
(m\cdot n)(1) \,=\, m(n(1)) \,=\, m(g) \,=\, g\cdot f\;.
\]
This means that, under $\alpha$, the multiplication in $\End_{\Delta_\infty}(\Delta_\infty)$ passes to the opposite multiplication which shows our claim. 
}
Our proof now completes by $\inv_\infty E \,=\, \inv\Delta_\infty^{\op}\,=\, -\inv\Delta_\infty \,=\, \frac{k}{l} \,=\, \weight(\FF)$\,.

\smallskip\noindent
\ref{THEOREM-2.6}. We prove the general case using local (iso-)shtukas rather than Tate modules which were used in \ref{THEOREM-2.5}. Our method is inspired by Milne's and Waterhouse' computation for abelian varieties~\cite[Theorem~8]{WM}. However in the function field case this method can be used to calculate the Hasse invariant at \emph{all} places, whereas in the number field case it applies only to the place which equals the characteristic of the ground field.
Let $w$ be a place of $Q$ and let $\ulN_w:=\ulN_w(\FF)$ be the local $\sigma$-isoshtuka of $\FF$ at $w$. Let $\BF_w$ be the residue field of $w$ and $\BF_{q^f}=\BF_w\cap\Fs$ the intersection inside an algebraic closure of $\Fq$. Let $\Fa_0$ be the ideal $(b\otimes 1-1\otimes b:b\in\BF_{q^f})$ of $Q_w\otimes_\Fq \Fs$ and let $R:=(Q_w\otimes_{\BF_q}\BF_s/\Fa_0)[T]=Q_w\otimes_{\BF_{q^f}}\Fs[T]$ be the non-commutative polynomial ring with $T\cdot(a\otimes b)=(a\otimes b^{q^f})\cdot T$ for $a\in Q_w$ and $b\in\BF_s$. Since $Q_w\otimes_{\BF_{q^f}}\Fs$ is a field, $R$ is a non-commutative principal ideal domain as studied by Jacobson~\cite[Chapter 3]{Jacobson}. Its center is the commutative polynomial ring $Q_w[T^g]$ where $g=[\Fs:\BF_{q^f}]=\frac{e}{f}$. From Theorem~\ref{ThmLS1} and Proposition~\ref{PropLS3} we get isomorphisms
\[
\QEnd(\FF)\otimes_Q Q_w\,\cong\,\End_{Q_w\otimes_\Fq\Fs[\phi]}(\ulN_w)\,\cong\,\End_R(\ulN_w/\Fa_0\ulN_w)
\]
where $T$ operates on $\ulN_w/\Fa_0\ulN_w$ as $\phi^f$.

By \cite[Theorem~3.19]{Jacobson} the $R$-module $\ulN_w/\Fa_0\ulN_w$ decomposes into a finite direct sum indexed by some set $I$
\begin{equation}\label{Eq3.13.6}
\ulN_w/\Fa_0\ulN_w\,\cong\,\bigoplus_{v\in I}\ulN_v^{\oplus n_v}
\end{equation}
of indecomposable $R$-modules $\ulN_v$ with $\ulN_v\not\cong \ulN_{v'}$ for $v\ne v'$. The annihilator of $\ulN_v$ is a two sided ideal of $R$ generated by a central element $\mu_v\in Q_w[T^g]$ by \cite[\S 3.6]{Jacobson}, which can be chosen to be monic. In particular (\ref{Eq3.13.6}) is an isomorphism of $Q_w[T^g]$-modules and $\mu_v$ is the minimal polynomial of $T^g$ on $\ulN_v$ by \cite[Lemma 3.1]{Jacobson}. Therefore the least common multiple $\mu$ of the $\mu_v$ is the minimal polynomial of $T^g$ on $\ulN_w/\Fa_0\ulN_w$. Note that $T^g$ operates on $\ulN_w/\Fa_0\ulN_w$ as the Frobenius $\pi$, hence $\mu=mipo_{\pi|\FF}$ and $F=Q(\pi)=Q[T^g]/(\mu)$, where we write $mipo$ for the minimal polynomial. By the semisimplicity of $\pi$ (and Proposition~\ref{PROP.EQ}) $\mu$ has no multiple factors in $Q_w[T^g]$. Since the $\mu_v$ are powers of irreducible polynomials by \cite[Theorem~3.20]{Jacobson} we conclude that all $\mu_v$ are themselves irreducible in $Q_w[T^g]$. Again \cite[Theorem~3.20]{Jacobson} implies that $\mu_v\ne \mu_{v'}$ since $\ulN_v\not\cong \ulN_{v'}$ and
\[
\mu\,=\,mipo_{\pi|\FF}\,=\,\prod_{v\in I}\mu_v\quad\text{inside}\quad Q_w[T^g]\,.
\]
Thus $F\otimes_Q Q_w=Q_w[T^g]/(\mu)=\prod_{v\in I}Q_w[T^g]/(\mu_v)=\prod_{v|w}F_v$. So $I$ is the set of places of $F$ dividing $w$ and $F_v=Q_w[T^g]/(\mu_v)$ is the completion of $F$ at $v$, justifying our notation. Let $\pi_v$ be the image of $\pi$ in $F_v$. Its minimal polynomial over $Q_w$ is $\mu_v$. This implies that $E\otimes_Q Q_w$ decomposes further
\[
E\otimes_Q Q_w\,=\,\bigoplus_{v\in I}\End_R(\ulN_v^{\oplus n_v})\,=\,\bigoplus_{v\in I}E\otimes_F F_v
\]
and $E\otimes_F F_v\cong\End_R(\ulN_v^{\oplus n_v})$. 

Now fix a place $v$ above $w$ and consider the diagram of field extensions
\[
\xymatrix @R-1pc @C-1pc {& & \BF_v\BF_s\ar@{-}[d] \\
\BF_v\ar@{-}[ddd]\ar@{-}[urr]^{g/h}\ar@{-}[dr] & & \BF_w\BF_s\ar@{-}[dd]^i\\
& \BF_w\BF_s\cap\BF_v \ar@{-}[ur]^{g/h}\ar@{-}[d]\\
& \BF_w(\BF_v\cap\BF_s)\ar@{-}[d]^i & \BF_s\\
\BF_w \ar@{-}[d]^i\ar@{-}[ur]^{h\es} & \BF_v\cap\BF_s\ar@{-}[ur]_{g/h}\\
**{!L !<0.8pc,0pc> =<1.5pc,1.5pc>}\objectbox{\BF_{q^f}=\BF_w\cap\BF_s=\BF_w\cap(\BF_v\cap\BF_s)} \ar@{-}[ru]^h \ar@{-}[d]^f\\
\BF_q
}
\]
Let $h:=[\BF_v\cap\Fs:\BF_{q^f}]=\gcd([\BF_v:\BF_{q^f}],g)$. Let $i:=[\BF_w:\BF_{q^f}]$. From the formulas
\begin{eqnarray*}
[\BF_w\BF_s:\BF_w]&=&[\BF_s:\BF_{q^f}]\es=\es g,\\[1mm]
[\BF_w(\BF_v\cap\BF_s):\BF_w]&=&[\BF_v\cap\BF_s:\BF_{q^f}]\es=\es h,\\[1mm]
[\BF_w\BF_s:(\BF_w\BF_s\cap\BF_v)]&=&[\BF_v\BF_s:\BF_v]\es=\es[\BF_s:\BF_v\cap\BF_s]\es=\es{\TS\frac{g}{h}},\quad\text{and}\\[1mm]
\BF_w(\BF_v\cap\BF_s)&\subset&\BF_w\BF_s\cap\BF_v,
\end{eqnarray*}
we obtain $\BF_w\BF_s\cap\BF_v=\BF_w(\BF_v\cap\BF_s)=\BF_{q^{fhi}}$.
Let $F_{v,L}$ be the compositum of $Q_w\otimes_{\BF_{q^f}}\Fs$ and $F_v$ in an algebraic closure of $Q_w$. Note that $F_{v,L}$ is well defined since $\Fs/\BF_{q^f}$ is Galois. Let $F_{v,L}[T']$ be the non-commutative polynomial ring with
\[
T'\cdot(a\otimes b)=(a\otimes b^{q^{fhi}})\cdot T'\quad\text{and}\quad T'\cdot x= x\cdot T'
\]
for $a\in Q_w$, $b\in\Fs$, and $x\in F_v$ and set $\Delta_v=F_{v,L}[T']/\bigl((T')^{g/h}-\pi_v^i\bigr)$.
Observe that the commutation rules of $T'$ are well defined since $(Q_w\otimes_{\BF_{q^f}}\Fs)\cap F_v$ has residue field $\BF_w\Fs\cap\BF_v=\BF_{q^{fhi}}$ and is unramified over $Q_w$, because $Q_w\otimes_{\BF_{q^f}}\Fs$ is. Moreover, the extension $F_{v,L}/F_v$ is unramified of degree $[\BF_v\BF_s:\BF_v]=\frac{g}{h}$ and $\wt T:=(T')^{[\BF_v:\Fq]/fhi}$ is its Frobenius automorphism. Since $\wt T^{g/h}=\pi_v^{[\BF_v:\Fq]/fh}$ in $\Delta_v$, our $\Delta_v$ is just the cyclic algebra $\bigl(F_{v,L}/F_v,\wt T,\pi_v^{[\BF_v:\Fq]/fh}\bigr)$ and has Hasse invariant $\frac{[\BF_v:\Fq]}{[\Fs:\Fq]}\cdot v(\pi_v)$; compare \cite[p.\ 266]{Re}. We relate $\Delta_v$ to $E\otimes_F F_v$. Firstly by \cite[Theorem~3.20]{Jacobson} there exists a positive integer $u$ such that $\ulN_v^{\oplus u}\cong R/R \mu_v(T^g)$. Therefore
\[
M_u(E\otimes_F F_v)\,\cong\,M_u\bigl(\End_R(\ulN_v^{\oplus n_v})\bigr)\,=\,\End_R(\ulN_v^{\oplus un_v})\,=\,M_{n_v}\bigl((R/R\mu_v(T^g))^\op\bigr)\,.
\]
Secondly we choose integers $m$ and $n$ with $m>0$ and $mi+ng=1$. We claim that the morphism $R/R\mu_v(T^g)\to M_h(\Delta_v)$, which maps
\[
a\otimes b\longmapsto\left(\begin{array}{cccc}a\otimes b\\ & a\otimes b^{q^f}\\ & & \ddots \\ & & & a\otimes b^{q^{f(h-1)}} \end{array}\right)\qquad\text{and}\qquad T\longmapsto 
\pi_v^n\cdot\left( \raisebox{4ex}{$
\xymatrix @R=0.3pc @C=0.7pc {0 \ar@{.}[drdr] & 1\ar@{.}[dr]\\
& & 1 \\
(T')^m & & 0}$} \right)
\]
for $a\in Q_w$ and $b\in\Fs$, is an isomorphism of $F_v$-algebras. It is well defined since it maps $T\cdot(a\otimes b)$ and $(a\otimes b^{q^f})\cdot T$ to the same element because $(T')^m=(T')^{1/i}$ in $\Gal(F_{v,L}/F_v)$, and it maps $T^g=(T^h)^{g/h}$ to $\pi_v^{ng}(T')^{mg/h}\cdot\Id_h=\pi_v\cdot\Id_h$. Since $R\mu_v(T^g)\subset R$ is a maximal two sided ideal the morphism is injective. To prove surjectivity we compare the dimensions as $Q_w$-vector spaces. We compute
\begin{eqnarray*}
\dim_{F_v} M_h(\Delta_v)&=&\TS h^2\cdot(\frac{g}{h})^2\es=\es g^2\,,\\[2mm]
\dim_{Q_w\otimes_{{\SC\BF}_{q^f}}\Fs}\bigl(R/R\mu_v(T^g)\bigr) &=&g\cdot\deg \mu_v\es=\es g\cdot[F_v:Q_w]\,, \quad\text{and}\\[2mm]
\dim_{Q_w}\bigl(R/R\mu_v(T^g)\bigr) &=&g^2\cdot[F_v:Q_w]\es=\es\dim_{Q_w}M_h(\Delta_v)\,.
\end{eqnarray*}
Altogether $M_u(E\otimes_F F_v)\cong M_{hn_v}(\Delta_v^\op)$ and $\inv_v E=-\inv_v\Delta_v=-\frac{[\BF_v:\Fq]}{[\Fs:\Fq]}\cdot v(\pi_v)$ as claimed.

It remains to convert this formula into the special form asserted for $v\nmid\chr\infty$ or $v|\infty$. If $v|\infty$ and $\chr\ne\infty$, let $e_v$ be the ramification index of $F_v/Q_\infty$. Then we get from Theorem~\ref{ThmRH} the formula $q^{e\weight(\FF)}=|\pi|_\infty=q^{-v(\pi_v)/e_v}$, since the residue field of $Q_\infty$ is $\Fq$. This implies as desired
\[
-\frac{[\BF_v:\Fq]}{[\Fs:\Fq]}\cdot v(\pi_v)\,=\,-\frac{[\BF_v:\Fq]\cdot  (-e_ve\cdot\weight(\FF))}{e}\,=\,\weight(\FF)\cdot[F_v:Q_\infty]
\]

Finally if $w\ne\chr,\infty$ is a place of $Q$, the local $\sigma$-shtuka $\ulM_w(\FF)$ at $w$ is \'etale. So $\mu=mipo_{\pi|\FF}$ has coefficients in $A_w$ with constant term in $A_w^\times$. Therefore $v(\pi_v)=0$ for all places $v$ of $F$ dividing $w$.
\end{proof}

\begin{Example}\label{LAST-EXAMPLE}
Let $C=\PP^1_\Fq$, $C\setminus\{\infty\}=\Spec\Fq\II{t}$ and $L=\Fq$. Let $d$ be a positive integer. Let $\F_i:=\O(d\lceil\frac{i}{2}\rceil\cdot\infty)\oplus\O(d\lceil\frac{i-1}{2}\rceil\cdot\infty)$ for $i\in\Z$ and let $\t:=\matr{0}{1}{t^d}{0}$. Then $\FF=(\F_i,\P_i,\t_i)$ is an abelian $\tau$-sheaf of rank $2$, dimension $d$, and characteristic $\chr=V(t)\in\PP^1$ over $\Fq$. Hence the Frobenius endomorphism $\pi$ equals $\t$. If $d$ is odd then $\FF$ is primitive (that means $(d,r)=1$) and therefore simple by Proposition~\ref{PROP.5}. In particular, $\pi$ is semisimple. We have
\[
\mu_\pi\,=\,\chi_v \,=\, x^2-t^d \,=\, (x-\sqrt{t^d})(x+\sqrt{t^d})
\]
which means that $\pi_v$ is \emph{not} absolutely semisimple in characteristic $2$. Moreover, we calculate $r_\Qv(\chi_v,\chi_v)=1\cdot1\cdot2=2$ whereas in the field extension $\Qv(\sqrt{t})\,/\,\Qv$ we have
\[
r_{\Qv(\sqrt{t})}(\chi_v,\chi_v) \,=\, \left\{
\begin{array}{l@{\quad}l} 
2\cdot2\cdot1=4 & \mbox{in characteristic 2,} \\ 
1\cdot1\cdot1+1\cdot1\cdot1=2 & \mbox{in characteristic different from 2.} 
\end{array}
\right.
\]
Although the later has no further significance it illustrates the remark after Definition~\ref{Def3.3}. By Theorem \ref{THEOREM-2}/\ref{THEOREM-2.4}. we have $E=F=Q(\pi)$ commutative and $\II{E:Q}=2=r$. Moreover, $|\,\pi\,|_\infty=|\,\sqrt{t^d}\,|_\infty=q^{d/2}$ and $\chi_v$ is irreducible. But $\chi_v$ is not separable in characteristic 2.

If $d=2n$ is even then the minimal polynomial of $\pi$ is
\[
\mu_\pi\,=\,\chi_v \,=\, x^2 -t^d\,=\, (x-t^{d/2})(x+t^{d/2})\,.
\]
So $\pi$ is semisimple if and only if $\charakt(\Fq)\ne2$. In this case $\FF$ is quasi-isogenous to the abelian $\tau$-sheaf $\FF'$ with $\F'_i=\O_{C_L}(in\cdot\infty)^{\oplus 2}$ and $\t'_i=\matr{-t^n}{0}{0}{t^n}$. The quasi-isogeny $f:\FF'\to\FF$ is given by $f_{0,\eta}=\matr{-t^n}{1}{t^n}{1}:\F'_{0,\eta}\isoto\F_{0,\eta}$. The abelian $\tau$-sheaf $\FF'$ equals the direct sum $\FF^{(1)}\oplus\FF^{(2)}$ where $\F_i^{(j)}=\O_{C_L}(in\cdot\infty)$ and $\t_i^{(j)}=(-1)^j t^n$. Note that $\FF^{(1)}$ and $\FF^{(2)}$ are not isogenous over $\Fq$, since the equation $-t^n\cdot\s(g)=g\cdot t^n$ has no solution $g\in Q$ for $\charakt(\Fq)\ne2$. Therefore 
\[
Q\oplus Q\,=\,\bigoplus_{j=1}^2 \QEnd(\FF^{(j)})\,\cong\,E\,=\,F\,=\,Q[x]/(x^2-t^{2n})\,\cong\,Q\oplus Q\,.
\]

Now we consider the same abelian $\tau$-sheaf over $L=\Ff_{q^2}$. This means $\pi=\t^2=t^d\in Q$ and therefore $\chi_v=(x-t^d)^2$. Thus $\pi$ is semisimple. By Theorem \ref{THEOREM-2}/\ref{THEOREM-2.5} we have $F=Q(\pi)=Q$ and $E$ is central simple over $Q$ with $\II{E:Q}=4$ and $\inv_\infty E=\inv_\chr E=\frac{d}{2}$. Moreover, $|\,\pi\,|_\infty=|\,t^d\,|_\infty=q^d$. In this case, $\pi_v$ is absolutely semisimple. Note that if $d$ is even and $\charakt(\Fq)=2$ this is another example for Theorem~\ref{Thm3.8b}. 

If $d$ is odd then $\FF$ is still primitive, whence simple and $E$ is a division algebra. If $d=2n$ is even then the abelian $\tau$-sheaves $\FF^{(1)}$ and $\FF^{(2)}$ defined above are isomorphic $\FF^{(1)}\isoto\FF^{(2)},1\mapsto\lambda$ where $\lambda\in\BF_{q^2}$ satisfies $\lambda^{q-1}=-1$. Therefore $M_2(Q)=M_2\bigl(\QEnd(\FF^{(1)})\bigr)\cong E$ in accordance with the Hasse invariants just computed.
\end{Example}

\begin{Example}\label{Ex3.15}
We compute another example which displays other phenomena. Let $C=\PP^1_\Fq$ and let $C\setminus\{\infty\}=\Spec\Fq[t]$. Let $\F_i=\O_{C_L}(\lceil\frac{i-1}{2}\rceil\cdot\infty)^{\oplus2}\oplus\O_{C_L}(\lceil\frac{i}{2}\rceil\cdot\infty)^{\oplus2}$, let $\P_i$ be the natural inclusion, and let $\t_i$ be given by the matrix
\[
T\,:=\,\left( \begin{array}{cccc}0&0&0&a\\0&b&1&0\\t&0&-b&0\\0&t&0&0 \end{array}\right)\qquad\text{with}\quad a,b\in\Fq\setminus\{0\}\,.
\]
Then $\FF$ is an abelian $\tau$-sheaf of rank $4$ and dimension $2$ with $l=2,k=1$ and characteristic $\chr=V(t)\in\PP^1$. One checks that the minimal polynomial of the matrix $T$ is $x^4-b^2x^2-at^2$ which is irreducible over $Q$ if $\charakt(\Fq)\ne2$, since it has neither zeroes in $\Fq[t]$ nor quadratic factors in $Q[x]$. If $\charakt(\Fq)=2$ then the minimal polynomial is a square and $\FF$ is not semisimple.

For $L=\Fq$ and $2\nmid q$ we obtain $\pi=\t$ semisimple and $E=F=Q(\pi)=Q[x]/(x^4-b^2x^2-at^2)$.

For $L=\BF_{q^2}$ we have $\pi=\t^2$ and the minimal polynomial of $\pi$ over $Q$ is $x^2-b^2x-at^2$, which is irreducible also in characteristic $2$ since it has no zeroes in $\Fq[t]$. Hence $\pi$ is semisimple, $F$ is a field with $[F:Q]=2$ and $[E:F]=4$ by Corollary~\ref{Cor3.11b}. This again illustrates Theorem~\ref{Thm3.8b}. We compute the decomposition of $\infty$ and $\chr$ in $F$.

\smallskip

\noindent
\emph{Decomposition of $\chr$:} Modulo $t$ the polynomial $x^2-b^2x-at^2$ has two zeroes $x=b^2$ and $x=0$ in $\Fq$. So by Hensel's lemma $F\otimes_Q Q_\chr\cong F_v\oplus F_{v'}$ splits with $F_v\cong F_{v'}\cong Q_\chr$ and $v(\pi)=0$ and $v'(\pi)=v'(at^2)=2$. Thus the Hasse invariants of $E$ are $\inv_v E=\inv_{v'}E=0$.

\smallskip

\noindent
\emph{Decomposition of $\infty$:} Set $y=\pi/t$. Then $y^2-\frac{b^2}{t}y-a=0$.\\[1mm]
\emph{Case (a).} If $2|q$ then $(y-a^{q/2})^2-\frac{b^2}{t}(y-a^{q/2})-\frac{b^2}{t}a^{q/2}=0$, that is, $\infty$ ramifies in $F$, $F\otimes_Q Q_\infty=F_w$ with $w(\frac{\pi}{t}-a^{q/2})=1$ and $w(\frac{1}{t})=2\cdot\infty(\frac{1}{t})=2$. So $[F_w:Q_\infty]=2$ and $\inv_w E=0$.

\smallskip

\noindent
\emph{Case (b).} If $2\nmid q$ and $\sqrt{a}\in\Fq$ then the polynomial $y^2-\frac{b^2}{t}y-a$ has two zeroes $y=\pm\sqrt{a}$ modulo $\frac{1}{t}$. So by Hensel's lemma $F\otimes_Q Q_\infty\cong F_w\oplus F_{w'}$ splits with $[F_w:Q_\infty]=[F_{w'}:Q_\infty]=1$. Thus the local Hasse invariants of $E$ are $\inv_w E=\inv_{w'}E=\frac{1}{2}$. As was remarked in \ref{Rem3.14b} such a distribution of the Hasse invariants can occur only if $d\ge2$.

\smallskip

\noindent
\emph{Case (c).} If $2\nmid q$ and $\sqrt{a}\notin\Fq$ then $y^2-\frac{b^2}{t}y-a$ is irreducible modulo $\frac{1}{t}$ and $\infty$ is inert in $F$, $F\otimes_Q Q_\infty=F_w$ with $[F_w:Q_\infty]=2$. Thus the Hasse invariant of $E$ is $\inv_w E=0$.

\medskip

In case (b) $E$ is a division algebra and $\FF$ is simple. In cases (a) and (c) $E\cong M_2(F)$ and $\FF$ is quasi-isogenous to $(\FF')^{\oplus2}$ for an abelian $\tau$-sheaf $\FF'$ of rank $2$, dimension $1$ and $\QEnd(\FF')=F$. This surprising result is due to the fact that $\FF'$, being of dimension $1$, is associated with a Drinfeld module and thus of the form $\F'_i=\O_{C_L}(\lceil\frac{i}{2}\rceil\cdot\infty)\oplus\O_{C_L}(\lceil\frac{i-1}{2}\rceil\cdot\infty)$ with $\t'_i=\matr{c}{d}{t}{0}$ and $c,d\in\BF_{q^2}$. Then $\pi'=(\t')^2=\matr{c^{q+1}+d^qt}{c^qd}{ct}{dt}$ has minimal polynomial $x^2-(c^{q+1}+(d+d^q)t)x+d^{q+1}t^2$ which must be equal to $x^2-b^2x-at^2$. This is possible only if $d+d^q=0$ and $d^{q+1}=-a$. So either $d\in\Fq$ and $2|q$ and we are in case (a), or $d\in\BF_{q^2}\setminus\Fq$, $d^q=-d$, and $a=d^2$. The later implies $2\nmid q$ and $\sqrt{a}=d\notin\Fq$ and we are in case (c). If we choose $c=b$ in case (c) a quasi-isogeny $f:\FF\to(\FF')^{\oplus2}$ over $\BF_{q^2}$ is given for instance by
\[
\left(\begin{array}{cccc} 
d & a & -bd/t & 0 \\
0 & 0 & -d & a \\
0 & 0 & d/t & a/t \\
1 & -d & 0 & bd/t 
\end{array}\right)\,.
\]
\end{Example}

% -----------------------------------------------------------------------------

\bigskip

\subsection{Kernel Ideals for Pure Anderson Motives} \label{Sect1.7}

In this section we investigate which orders of $E$ can arise as endomorphism rings $\End(\ulM)$ for pure Anderson motives $\ulM$. For this purpose we define for each right ideal of the endomorphism ring $\End(\ulM)$ an isogeny with target $\ulM$ and discuss its properties. This generalizes Gekeler's results for Drinfeld modules \cite[\S 3]{Gekeler} and translates the theory of Waterhouse \cite[\S 3]{Wat} for abelian varieties to the function field case. These two sources are themselves the translation, respectively the higher dimensional generalization of Deuring's work on elliptic curves~\cite{Deuring}.

Let $\ulM$ be a pure Anderson motive over $L$ and abbreviate $R:=\End(\ulM)$. Let $I\subset R$ be a right ideal which is an $A$-lattice in $E:=R\otimes_A Q$. This is equivalent to saying that $I$ contains an isogeny, since every lattice contains some isogeny $a\cdot\id_\ulM$ for $a\in A$ and conversely the existence of an isogeny $f\in I$ implies that the lattice $f\cdot\dual{f}\cdot R$ is contained in $I$.

\begin{Definition}\label{Def1.7.1}
\begin{suchthat}
\item 
Let $\ulM^I$ be the pure Anderson sub-motive of $\ulM$ whose underlying $A_L$-module is $\sum_{g\in I}\im(g)$. This is indeed a pure Anderson motive, since if $I=f_1R+\ldots+f_nR$ are arbitrary generators, then $\ulM^I$ equals the image of the morphism 
\[
(f_1,\ldots,f_n):\ulM\oplus\ldots\oplus\ulM\es\longto\es\ulM\,.
\]
As $I$ contains an isogeny, $\ulM^I$ has the same rank as $\ulM$ and the natural inclusion is an isogeny which we denote $f_I:\ulM^I\to \ulM$.
\item 
If $I=\{\,f\in R:\im(f)\subset\ulM^I\,\}$ then $I$ is called a \emph{kernel ideal for $\ulM$}.
\end{suchthat}
\end{Definition}

The later terminology is borrowed from Waterhouse~\cite[\S 3]{Wat}. Since $\{\,f\in R:\im(f)\subset\ulM^I\,\}$ is the right ideal annihilating $\coker f_I$ one should maybe use the name ``cokernel ideal'' instead.

\begin{Proposition}\label{Prop3.21b}
Let $I\subset R$ be a right ideal which is a lattice, and consider the right ideal $J:=\{\,f\in R: \im(f)\subset\ulM^I\,\}\subset R$ containing $I$. Then $\ulM^{J}=\ulM^I$. In particular, $J$ is a kernel ideal for $\ulM$. We call $J$ the \emph{kernel ideal for $\ulM$ associated with $I$}.
\end{Proposition}

\begin{proof}
Obviously $J$ is a right ideal and $\ulM^{J}\subset\ulM^I$ by definition of $J$. Conversely $\ulM^I\subset\ulM^{J}$ since $I\subset J$.
\end{proof}

\begin{Lemma}\label{Lemma1.7.2}
\begin{suchthat}
\item 
For any $g\in I$, $f_I^{-1}\circ g:\ulM\to\ulM^I$ is a morphism and $g=f_I\circ(f_I^{-1}\circ g)$.
\item 
If $I=gR$ is principal, $g$ an isogeny, then $f_I^{-1}\circ g:\ulM\to\ulM^I$ is an isomorphism and $I$ is a kernel ideal.
\end{suchthat}
\end{Lemma}

\begin{proof}
1 is obvious since the image of $g$ lies inside $\ulM^I$.\\
2. Clearly $f_I^{-1}\circ g$ is injective since $g$ is an isogeny and surjective by construction, hence an isomorphism. To show that $I$ is a kernel ideal let $f\in R$ satisfy $\im(f)\subset\ulM^I$. Consider the diagram
\[
\xymatrix @C+2pc {\ulM \ar@{-->}[d]_h \ar[r]^{f_I^{-1}\circ f} & \ulM^I \ar[r]^{f_I} & \ulM\\
\ulM \ar[ur]_{f_I^{-1}\circ g}
}
\]
and let $h:=(f_I^{-1}\circ g)^{-1}\circ(f_I^{-1}\circ f)$. Then $f=gh\in I$ as desired.
\end{proof}

\noindent
{\it Example.} If $a\in A$ and $I=aR$, then $\ulM^I=a\ulM$ and $\coker f_I=\ulM/a\ulM$. More generally if $\Fa\subset A$ is an ideal and $I=\Fa R$ then $\ulM^I=\Fa\ulM$ and $\coker f_I=\ulM/\Fa\ulM$.

\begin{Proposition}\label{Prop1.7.3}
Let $I\subset R$ and $J\subset \End(\ulM^I)$ be right ideals which are lattices in $E$. Then also the product $K:=f_I\cdot J\cdot f_I^{-1}\cdot I$ is a right ideal of $R$ and a lattice in $E$ and $f_K^{-1}\circ f_I\circ f_J$ is an isomorphism of $(\ulM^I)^J$ with $\ulM^K$
\[
(\ulM^I)^J \xrightarrow{\es f_J\;} \ulM^I \xrightarrow{\es f_I\;} \ulM \xleftarrow{\;f_K\es} \ulM^K\,.
\]
\end{Proposition}

\begin{proof}
  If $f\in I$ and $g\in J$ then the morphism $f_I^{-1}\circ f:\ulM\to\ulM^I$ can be composed with $f_I\circ g$ to yield an element of $R$. Since $I$ and $J$ contain isogenies, $K$ is a right ideal and contains an isogeny. Clearly the images of $f_I\circ f_J$ and $f_K$ in $\ulM$ coincide since they equal the sum $\sum_{i,j}f_I\circ g_j\circ f_I^{-1} \circ f_i(\ulM)$ for sets of generators $\{f_i\}$ of $I$ and $\{g_j\}$ of $J$.
\end{proof}

\begin{Theorem}\label{Thm1.7.4}
Let $I,J\subset\End(\ulM)=:R$ be right ideals which are lattices in $E:=R\otimes_A Q$ and consider the following assertions:
\begin{suchthat}
\item 
$I$ and $J$ are isomorphic $R$-modules,
\item 
the pure Anderson motives $\ulM^I$ and $\ulM^J$ are isomorphic.
\end{suchthat}
Then 1 implies 2 and if moreover $I$ and $J$ are kernel ideals, also 2 implies 1.
\end{Theorem}

\begin{proof}
$1\Rightarrow 2$. Since $I$ and $J$ are lattices, the $R$-isomorphism $I\to J$ extends to an $E$-isomorphism of $E$ and is thus given by left multiplication with a unit $g\in E^\times$, that is, $J=gI$. There is an $a\in A$ such that $ag\in I\subset R$. Then $\im(ag)\subset\ulM^I$, that is, $f_I^{-1}\circ ag:\ulM\to\ulM^I$ is an isogeny.

Let $K$ be the right ideal $f_I\cdot\bigl(f_I^{-1}\circ ag\circ f_I\cdot\End(\ulM^I)\bigr)\cdot f_I^{-1}\cdot I$ of $R$. We claim that $\ulM^K\cong\ulM^{(ag)I}$. Namely, $\ulM^{(ag)I}\subset\ulM^K$ since $agI\subset K$. Conversely if $f\in I$, $h\in\End(\ulM^I)$, and $m\in\ulM$, then we find $m':=f_I\circ h\circ f_I^{-1}\circ f(m)\in\ulM^I$, that is, $m'=\sum_i f_i(m_i)$ for suitable $f_i\in I$ and $m_i\in \ulM$. It follows that $ag(m')=\sum_i agf_i(m_i)\in\ulM^{(ag)I}$ and therefore $\ulM^{(ag)I}=\ulM^K$. 

Applying Lemma~\ref{Lemma1.7.2} and Proposition~\ref{Prop1.7.3} now yields an isomorphisms $\ulM^I\cong\ulM^K=\ulM^{(ag)I}$. Likewise we obtain $\ulM^J\cong\ulM^{aJ}$ and the equality $aJ=agI$ then implies $\ulM^J\cong\ulM^I$ as desired.

\smallskip
\noindent
$2\Rightarrow 1$. Let $I$ and $J$ be kernel ideals and let $u:\ulM^I\to\ulM^J$ be an isomorphism. There is an $a\in A$ with $a\ulM\subset\ulM^I$. Therefore $g:=f_J\circ u\circ(f_I^{-1}\circ a):\ulM\to \ulM$ is an isogeny. 

We claim that $gI=aJ$, that is, left multiplication by $a^{-1}g$ is an isomorphism of $I$ with $J$. Let $f\in I$, then $h:=f_J\circ u\circ(f_I^{-1}\circ f)\in R$ has $\im(h)\subset\ulM^J$. So $h\in J$ since $J$ is a kernel ideal, and $g f=ah\in aJ$, since $a$ commutes with all morphisms. Conversely let $h\in J$, then $f:=f_I\circ u^{-1}\circ(f_J^{-1}\circ h)\in R$ has $\im(f)\subset \ulM^I$. So $f\in I$ since $I$ is a kernel ideal, and $ah=gf\in gI$ as desired.
\end{proof}

\begin{Proposition}\label{Prop1.7.5}
Let $I\subset R$ be a right ideal which is a lattice in $E$. Then $f_I\cdot\End(\ulM^I)\cdot f_I^{-1}$ contains the left order $\O=\{\,f\in E:fI\subset I\,\}$ of $I$ and equals it if $I$ is a kernel ideal.
\end{Proposition}

\noindent
{\it Remark.} Recall that $\End(\ulM^I)\otimes_A Q$ is identified with $E$ by mapping $h\in \End(\ulM^I)$ to $f_I\circ h\circ f_I^{-1}$.

\begin{proof}
Let $f\in \O$ and $g\in I$. Then $f g\in I$ and $f_I^{-1}\circ f\circ f_I\circ(f_I^{-1}\circ g)=f_I^{-1}\circ f g$ is a morphism from $\ulM$ to $\ulM^I$. If $g$ varies, the images of $f_I^{-1}\circ g$ exhaust all of $\ulM^I$. Hence $f_I^{-1}\circ f\circ f_I$ is indeed an endomorphism of $\ulM^I$. 
Conversely let $I$ be a kernel ideal and let $f=f_I\circ h\circ f_I^{-1}\in f_I\cdot\End(\ulM^I)\cdot f_I^{-1}$. If $g\in I$ then $f\circ g=f_I\circ h\circ(f_I^{-1}\circ g)\in R$ has $\im(f\circ g)\subset \ulM^I$. So $fg\in I$ as desired.
\end{proof}

We will now draw conclusions about the endomorphism ring $R$ similar to Waterhouse' results \cite{Wat} on abelian varieties by simply translating his arguments.

\begin{Theorem}\label{ThmW3.13}
Every maximal order in $E$ occurs as the endomorphism ring $f\cdot\End(\ulM')\cdot f^{-1}\subset E$ of a pure Anderson motive $\ulM'$ isogenous to $\ulM$ via an isogeny $f:\ulM'\to \ulM$.
\end{Theorem}

\begin{proof}
Let $S$ be a maximal order of $E$. Then the lattice $R$ contains $aS$ for some $a\in A$. Consider the right ideal $I=aS\cdot R$ whose left order contains $S$. By Proposition~\ref{Prop1.7.5}, $f_I\cdot\End(\ulM^I)\cdot f_I^{-1}$ contains the left order of $I$. Since $S$ is maximal we find $S=f_I\cdot\End(\ulM^I)\cdot f_I^{-1}$.
\end{proof}

\begin{Theorem}\label{ThmW3.14}
If $E$ is semisimple and $\End(\ulM)$ is a maximal order in $E$, so is $f_I\cdot\End(\ulM^I)\cdot f_I^{-1}$ for any right ideal $I\subset R$.
\end{Theorem}

\begin{proof}
By \cite[Theorem 21.2]{Re} the left order of $I$ is also maximal and then Proposition~\ref{Prop1.7.5} yields the result.
\end{proof}

%\bigskip

From now on we assume that $L$ is a finite field and we set $e:=[L:\Fq]$. Let $\pi$ be the Frobenius endomorphism of $\ulM$.

\begin{Proposition}\label{PropW3.5}
The order $R$ in $E$ contains $\pi$ and $\deg(\pi)/\pi$.
\end{Proposition}

\begin{proof}
Clearly the isogeny $\pi$ belongs to $R$. Let now $a\in\deg(\pi)$. Then $a$ annihilates $\coker\pi$ by \ref{Prop3.28a} and so there is an isogeny $f:\ulM\to\ulM$ with $\pi\circ f=a$. The image $a/\pi$ of $f$ in $E$ belongs to $R$.
\end{proof}

\begin{Proposition}\label{ThmW3.15}
If $\ulM$ is a semisimple pure Anderson motive over a finite field and $\End(\ulM)$ is a maximal order in $E=\End(\ulM)\otimes_A Q$, then every right ideal $I\subset\End(\ulM)$, which is a lattice, is a kernel ideal for $\ulM$, and $\deg(f_I)=N(I):=\bigl(N(f):f\in I\bigr)$.
\end{Proposition}

\begin{proof}
(cf.\ \cite[Theorem 3.15]{Wat}) Let $f\in I$, then $f=f_I\circ f_I^{-1}f$ and $N(f)\in\deg(f)\subset\deg(f_I)$ by Lemma~\ref{Lemma1.7.7}. Therefore $N(I)\subset\deg(f_I)$. Let $R'$ be the left order of $I$. It is maximal by \cite[Theorem 21.2]{Re}. For a suitable $a\in A$ the set $J':=\{\,x\in E: xI\subset aR\,\}$ is a right ideal in $R'$ and a lattice in $E$ and satisfies $J'\cdot I=aR$ by \cite[Theorem 22.7]{Re}. Let $J:=f_I^{-1}J'f_I\subset\End(\ulM^I)$ be the induced right ideal of $\End(\ulM^I)=f_I^{-1} R'f_I$; see \ref{Prop1.7.5}. Then $\coker f_I\circ f_J=\coker f_{J'I}=\coker a$ by Proposition~\ref{Prop1.7.3}. Therefore Theorem~\ref{Prop3.4.1} and \cite[24.12 and 24.11]{Re} imply
\[
N(a)\cdot A=N(J')\cdot N(I)\subset(\deg f_J)(\deg f_I)=\deg(a)=N(a)\cdot A\,.
\]
By the above we must have $N(I)=\deg(f_I)$ since $A$ is a Dedekind domain. If $I$ were not a kernel ideal its associated kernel ideal would be a larger ideal with the same norm. But this is impossible by \cite[24.11]{Re}.
\end{proof}

Like for abelian varieties there is a strong relation between the ideal theory of orders of $E$ and the investigation of isomorphy classes of pure Anderson motives isogenous to $\ulM$. We content ourselves with the following result which is analogous to Waterhouse~\cite[Theorem 6.1]{Wat}. The interested reader will find many other results without much difficulty.

\begin{Theorem}\label{ThmW6.1}
Let $\ulM$ be a simple pure Anderson motive of rank $r$ and dimension $d$ over the smallest possible field $\Fq$. Then
\begin{suchthat}
\item 
$\End(\ulM)$ is commutative and $E:=\End(\ulM)\otimes_A Q=Q(\pi)$.
\item 
All orders $R$ in $Q(\pi)$ containing $\pi$ are endomorphism rings of pure Anderson motives isogenous to $\ulM$. Any such order automatically contains $N(\pi)/\pi=N_{Q(\pi)/Q}(\pi)/\pi$.
\item 
For each such $R$ the isomorphism classes of pure Anderson motives isogenous to $\ulM$ with endomorphism ring $R$ correspond bijectively to the isomorphism classes of $A$-lattices in $E$ with order $R$.
\end{suchthat}
\end{Theorem}

\begin{proof} 
1 follows from \ref{Thm3.8} and \ref{Cor3.13b}.

\noindent
2. Let $R$ be an order in $Q(\pi)$ containing $\pi$ and let $v\ne\chr$ be a maximal ideal of $A$. Since $[E:Q]=r$ and $E_v$ is semisimple, there is by Lemma~\ref{Lemma3.4.1b} an isomorphism $E_v\isoto V_v\ulM$ of (left) $E_v$-modules given by $f\mapsto f(x)$ for a suitable $x\in V_v\ulM$. It identifies $R_v:=R\otimes_A A_v$ with a $\pi$-stable lattice $\Lambda_v=R_v\cdot x$ in $V_v\ulM$, which without loss of generality is contained in $T_v\ulM$. By Proposition~\ref{Prop2.7b} 
there is an isogeny $f:\ulM'\to \ulM$ of pure Anderson motives with $T_vf(T_v\ulM')=\Lambda_v$. By Theorem~\ref{TATE-CONJECTURE-MODULES} we conclude
\[
\End(\ulM')\otimes_A A_v =\End_{A_v[\pi]}(\Lambda_v)=R_v\,.
\]
For $v=\chr$ note that $Q_{\chr,L}=Q_\chr$ since $L=\Fq$. In particular $\BF_\chr=\BF_q$. Since $\dim_{Q_\chr}\ulN_\chr(\ulM)=r=[E:Q]$, Theorem~\ref{ThmLS2} together with Lemma~\ref{Lemma3.4.1b} show that $E_\chr$ is isomorphic to $\ulN_\chr(\ulM)$ as left $E_\chr$-modules. Since $R$ contains $\pi$, the image of $R_\chr:=R\otimes_A A_\chr$ in $\ulN_\chr(\ulM)$ is a local $\sigma$-subshtuka $\ulHM{}'$ of $\ulM_\chr(\ulM)$ of the same rank. (If it is not contained in $\ulM_\chr(\ulM)$, multiply it with a suitable $a\in A$.) Then Proposition~\ref{Prop2.18b} yields an isogeny of pure Anderson motives $f:\ulM'\to \ulM$ such that $\ulM_\chr(f)\bigl(\ulM_\chr(\ulM')\bigr)=\ulHM{}'$ and
\[
\End(\ulM')\otimes_A A_\chr=\End_{A_\chr[\phi]}(\ulHM{}')=R_\chr
\]
by Theorem~\ref{ThmLS2}.
Since each of these operations only modifies $\End(\ulM)$ at the respective place $v$, this shows that we may modify $\ulM$ at all places to obtain a pure Anderson motive $\ulM'$ with $\End(\ulM')=R$. Now the last statement follows from Proposition~\ref{PropW3.5} and Theorem~\ref{Prop3.4.1}.

\smallskip
\noindent
3. Let $R$ be such an order. By what we proved in 2 there is a pure Anderson motive $\ulTM$ for which all $T_v\ulTM\cong R_v$ and $\ulM_\chr(\ulTM)\cong R_\chr$. Let $I\subset R$ be a (right) ideal which is an $A$-lattice in $E$ and consider the isogeny $f_I:\ulTM^I\to\ulTM$. Under the above isomorphisms $T_v f_I(T_v\ulTM^I)\cong I\otimes_A A_v=:I_v$ and $\ulM_\chr f_I(\ulM_\chr\ulTM^I)\cong I\otimes_A A_\chr=:I_\chr$. Conversely if $f:\ulM'\to\ulTM$ is an isogeny then $\ulM_v f(\ulM_v\ulM')$ is a (left) $R_v$-module because $R=\End(\ulTM)$, hence isomorphic to an $R_v$-ideal $I_v$. This shows that any isogeny $f:\ulM'\to\ulTM$ is of the form $f_I:\ulTM^I\to\ulTM$.

If now $f\in R$ satisfies $\im(f)\subset\ulTM^I$ then $f\in I_\chr$ and $f\in I_v$ for all $v$ and therefore $f\in I$. This shows that every $I$ is a kernel ideal for $\ulTM$. By Proposition~\ref{Prop1.7.5}, $\End(\ulTM^I)$ is the (left) order of $I$. Since every lattice with order $R$ in $E$ is isomorphic to an ideal of $R$, we have
\[
\xymatrix @R-1pc {\{\,A\text{-lattices in }E\text{ with order }R\,\}/_\sim \ar@{=}[d] \\
 \{\,I\subset R\text{ Ideals with order }R\,\}/_\sim \es \ar[r]^{\sim\qquad\quad} &
\es \bigl\{\,\ulTM^I\xrightarrow{f_I}\ulTM\to\ulM\text{ with }\End(\ulTM^I)=R\,\bigr\}/_\sim
}
\]
and the assertion now follows from Theorem~\ref{Thm1.7.4}.
\end{proof}

% =============================================================================

\vfill
 
\noindent
\parbox[t]{8cm}{
Matthias Bornhofen  \\ 
Kolleg St. Sebastian\\
Hauptstr. 4\\
D -- 79252 Stegen \\
Germany  \\[0.1cm] }
\parbox[t]{11cm}{ 
Urs Hartl  \\ 
Institute of Mathematics  \\ 
University of Muenster\\
Einsteinstr.\ 62\\
D--48149 Muenster\\
Germany  \\[0.1cm] 
http:/\!/www.math.uni-muenster.de/u/urs.hartl/
}

% =============================================================================

{\small
\begin{thebibliography}{XXX}
\addcontentsline{toc}{section}{References}

\bibitem[An1]{Anderson} G.\ Anderson: $t$-Motives, \emph{Duke Math.\ J.\/} {\bfseries 53} (1986), 457--502.


\bibitem[An2]{Anderson2} G.\ Anderson: On Tate Modules of Formal $t$-Modules, \emph{Internat.\ Math.\ Res.\ Notices} {\bfseries 2} (1993), 41--52.




% FIX VVVVVVVVVVVVVVVVVVVVVVVVVVVVVVVVVVVVVVVVVVV
\bibitem[BH1]{BH_A} M.\ Bornhofen, U.\ Hartl: \emph{Pure Anderson Motives and Abelian $\tau$-Sheaves}, \emph{Math.\ Z.} {\bfseries } (2010); see also {\tt arXiv:0709.2809}.



\bibitem[BH2]{BH_B} M.\ Bornhofen, U.\ Hartl: Pure Anderson Motives over Finite Fields, \emph{J.\ Number Th.\/} {\bfseries 129}, n.\ 2 (2009), 247--283; see also {\tt arXiv:0709.2815}.


\bibitem[BS]{BS}A.\ Blum, U.\ Stuhler: Drinfeld Modules and Elliptic Sheaves, in: \emph{Vector Bundles on Curves: New Directions}, S.\ Kumar, G.\ Laumon, U.\ Stuhler, M.\ S.\ Narasimhan, eds., pp.\ 110--188, Lecture Notes in Mathematics 1649, Springer-Verlag, Berlin, etc.\ 1991.


\bibitem[Boe]{Boeckle02} G.\ B\"ockle: Global $L$-functions over function fields,  \emph{Math.\ Ann.} {\bfseries 323} (2002),  no.\ 4, 737--795.


%\bibitem[BLR]{BLR}S.\ Bosch, W.\ L\"utkebohmert, M.\ Raynand: \emph{ N\'eron Models}, Springer-Verlag, Berlin, etc (1990).


%\bibitem[Bo1]{Bo}S.\ Bosch: \emph{ Algebra}, Springer-Verlag, Berlin, etc (2001).


%\bibitem[Bo2]{Bo2}S.\ Bosch: \emph{ Lineare Algebra}, Springer-Verlag, Berlin, etc (2003).


\bibitem[Bou]{Bou} N.\ Bourbaki: \emph{ Alg\`ebre, chapitre 8}, Hermann, Paris (1958).


%\bibitem[Co]{Co}P.\ M.\ Cohn: \emph{ Algebra, volume 2}, John Wiley \& Sons, Ltd., Chichester (1991).


\bibitem[Deu]{Deuring} M.\ Deuring: Die Typen der Multiplikatorenringe elliptischer Funktionenk\"orper, \emph{Abh.\ Math.\ Sem.\ Hansischen Univ.\/} {\bfseries 14} (1941), 197--272.


\bibitem[Dr1]{Drinfeld} V.G.\ Drinfeld: Elliptic Modules, \emph{Math.\ USSR-Sb.\/} {\bfseries 23} (1976), 561--592.


\bibitem[Dr2]{Drinfeld3} V.G.\ Drinfeld: Commutative Subrings of Certain Noncommutative Rings, \emph{Funct.\ Anal.\ Appl.\/} {\bfseries 11} (1977), 9--12.


\bibitem[Dr3]{Drinfeld5} V.G.\ Drinfeld: Moduli variety of $F$-sheaves, \emph{Funct.\ Anal.\ Appl.\/} {\bfseries 21} (1987), no.\ 2, 107--122.


\bibitem[EGA]{EGA} A.\ Grothendieck: \emph{{\'E}lements de G{\'e}o\-m{\'e}trie Alg{\'e}\-brique},
Publ.\ Math.\ IHES {\bfseries 4}, {\bfseries 8}, {\bfseries 11}, {\bfseries 17}, {\bfseries 20}, {\bfseries 24}, {\bfseries 28}, {\bfseries 32}, Bures-Sur-Yvette, 1960--1967; see also
Grundlehren {\bfseries 166}, Springer-Verlag, Berlin etc.\ 1971.



%\bibitem[Eis]{Eis}D.\ Eisenbud: \emph{ Commutative Algebra --- with a view toward algebraic geometry}, Graduate Texts in Mathematics, Springer-Verlag, New York, etc (1995).


\bibitem[Fal]{Fal}G.\ Faltings: Endlichkeitss\"atze für abelsche Variet\"aten \"uber Zahlk\"orpern, \emph{Invent.\ Math.} {\bfseries 73} (1983), 349--366.


%\bibitem[Ful]{Ful}W.\ Fulton: \emph{ Algebraic Curves}, W.\ A.\ Benjamin, Inc., New York, Amsterdam (1969).


\bibitem[Gar]{Gardeyn2} F.\ Gardeyn: A Galois criterion for good reduction of $\tau$-sheaves, \emph{J.\ Number Theory\/} {\bfseries 97} (2002), 447--471.


\bibitem[Gek]{Gekeler} E.-U.\ Gekeler: On finite Drinfel'd modules, \emph{J.\ Algebra} {\bfseries 141} (1991), no.\ 1, 187--203. 


\bibitem[Ge]{Genestier} A.\ Genestier: \emph{Espaces sym{\'e}triques de Drinfeld\/}, Ast{\'e}risque {\bfseries 234}, Soc.\ Math.\ France, Paris 1996.


\bibitem[Gos]{Go}D.\ Goss: \emph{ Basic Structures of Function Field Arithmetic}, Springer-Verlag, Berlin, etc (1996).


%\bibitem[Ha]{Ha}R.\ Hartshorne: \emph{ Algebraic Geometry}, Springer-Verlag, New-York, Heidelberg (1977).


\bibitem[Ha1]{Hl} U.\ Hartl: Uniformizing the Stacks of Abelian Sheaves, in \emph{Number Fields and Function fields - Two Parallel Worlds, Papers from the 4th Conference held on Texel Island, April 2004}, G.\ van der Geer, B.\ Moonen, R.\ Schoof, Editors, Progress in Math.\ 239, Birkh\"auser-Verlag, Basel 2005, pp.\ 167--222.


\bibitem[Ha2]{HartlPSp} U.\ Hartl: Period Spaces for Hodge Structures in Equal Characteristic, to appear in \emph{Ann.\ of Math.} (2010), see also {\tt arXiv:math.NT/0511686}.





\bibitem[HH]{HH}U.\ Hartl, M.\ Hendler: \emph{ Change of Coefficients for Drinfeld Modules, Shtukas, and Abelian $\tau$-Sheaves}, Preprint on {\tt arXiv:math.NT/0608256} (2006).


\bibitem[Hei]{Heiden} G.-J.\ van der Heiden: \emph{Weil Pairing and the Drinfeld Modular Curve}, PhD thesis, University of Groningen, 2003.


\bibitem[Jac]{Jacobson} N.\ Jacobson: \emph{The Theory of Rings}, AMS Mathematical Surveys, vol.\ II, American Mathematical Society, New York, 1943.


\bibitem[Kat]{Katz1} N.\ Katz: $p$-adic properties of modular schemes and modular forms, \emph{Modular functions of one variable, III (Proc.\ Internat.\ Summer School, Univ.\ Antwerp, Antwerp, 1972)},  pp.\ 69--190. LNM {\bfseries 350}, Springer, Berlin, 1973. 


\bibitem[Kim]{Kim} W.~ Kim: \emph{Galois deformation theory for norm fields and its arithmetic applications}, PhD thesis, University of Michigan, 2009.

\bibitem[Lan]{La}S.\ Lang: \emph{ Abelian Varieties}, Springer, New York (1983).


\bibitem[Lau]{Laumon} G.\ Laumon: \emph{Cohomology of Drinfeld Modular Varieties I\/}, Cambridge Studies in Advanced Mathematics {\bfseries 41}, Cambridge University Press, Cambridge, 1996.


\bibitem[Pap]{Papanikolas} M.\ Papanikolas: Tannakian duality for Anderson-Drinfeld motives and algebraic independence of Carlitz logarithms, \emph{Invent.\ Math.} {\bfseries 171} (2008), 123--174.


\bibitem[PT]{PT} R.\ Pink, M.\ Traulsen: The isogeny conjecture for $t$-motives
  associated to direct sums of Drinfeld modules, \emph{J.\ Number Theory} {\bfseries 117}, no.\ 2 (2006), 355--375.


\bibitem[Rei]{Re}I.\ Reiner: \emph{ Maximal Orders}, Academic Press, London, etc (1975).


\bibitem[Ses]{Se}C.\ S.\ Seshadri: \emph{ Fibr\'es vectoriels sur les courbes alg\'ebriques}, Ast\'erisque 96, Soci\'et\'e math\'ematique de France (1982).


\bibitem[Sta]{Stalder} N.\ Stalder: \emph{Algebraic Monodromy Groups of $A$-Motives}, PhD thesis, ETH Zuerich, 2007.


\bibitem[Tae]{Taelman} L.\ Taelman: \emph{On $t$-Motifs}, PhD thesis, University of Groningen, 2007.


\bibitem[Tag]{Taguchi95b}Y.\ Taguchi: The Tate conjecture for t-motives, \emph{Proc.\ Amer.\ Math.\ Soc.\/} {\bfseries 123} (1995), n.\ 11, 3285--3287.


\bibitem[Tam]{Tam}A. Tamagawa: Generalization of Anderson's $t$-motives and Tate conjecture, in: \emph{Moduli Spaces, Galois Representations and $L$-Functions}, S\=urikaisekikenky\=usho K\=oky\=uroku, n.\ 884, Kyoto (1994), 154--159.


\bibitem[Tat]{Tat}J.\ Tate: Endomorphisms of Abelian Varieties over Finite Fields, \emph{Invent.\ Math.} {\bfseries 2} (1966), 134--144.


\bibitem[TW]{TW}Y.\ Taguchi, D.\ Wan: $L$-Functions of $\varphi$-Sheaves and Drinfeld Modules, \emph{J.\ Amer.\ Math.\ Soc.} {\bfseries 9}, $\mbox{n}^o$ 3 (1996), 755--781.


\bibitem[Wat]{Wat}W.\ C.\ Waterhouse: Abelian Varieties over Finite Fields, \emph{Ann.\ Sci.\ \'Ecole Norm.\ Sup.\ (4)}, {\bfseries 2} (1969), 521--560.


\bibitem[WM]{WM} W.\ Waterhouse, J.\ Milne: Abelian Varieties over Finite Fields, \emph{1969 Number Theory Institute (Proc.\ Sympos.\ Pure Math., Vol.\ XX, Stony Brook, N.Y., 1969)},  pp.\ 53--64, Amer.\ Math.\ Soc., Providence, R.I.\ 1971. 

%\bibitem[Yu]{Yu} J.-K.\ Yu: Isogenies of Drinfel'd modules over finite fields, \emph{J.\ Number Theory} {\bfseries 54} (1995),  no.\ 1, 161--171. 


\end{thebibliography}
}
\end{document}